\documentclass[11pt]{article}
\usepackage{amssymb, url}
\usepackage{amsmath,amsthm, xpatch}
\usepackage{mathtools}
\usepackage{dsfont}
\usepackage{verbatim}
\usepackage{graphicx}
\usepackage{epstopdf}
\usepackage{fullpage}
\usepackage[hang,flushmargin]{footmisc}

\usepackage{tikz}
\usetikzlibrary{calc}
\usetikzlibrary{decorations.pathreplacing,calligraphy}

\RequirePackage{doi}
\usepackage{hyperref}

\usepackage{enumitem}
\usepackage{datetime}
\usepackage{bbm}

\usepackage{color}

\newtheorem{definition}{Definition}
\newtheorem{theorem}{Theorem}[]
\newtheorem{corollary}[theorem]{Corollary}
\newtheorem{lemma}{Lemma}
\newtheorem{conjecture}{Conjecture}

\newtheorem{problem}{Problem}

\newtheorem*{remark}{Remark}

\newtheorem{claim}{Claim}[theorem]
\newtheorem{subclaim}{Claim}[claim]
\newtheorem{subsubclaim}{Claim}[subclaim]

\theoremstyle{plain}

\newenvironment{theoremlateproof}[1]
 {%
  \renewcommand{\thetheorem}{\ref{#1}}%
  \begin{proof}[Proof of Theorem~\ref{#1}]%
 }
 {\end{proof}}

\newenvironment{claimlateproof}[1]
 {%
  \renewcommand{\theclaim}{\ref{#1}}%
  \begin{proof}[Proof of Claim~\ref{#1}]%
 }
 {\end{proof}}

\newcounter{caseCount}
\AtBeginEnvironment{claim}{\setcounter{caseCount}{0}}

\providecommand{\keywords}[1]
{
  \small	
  \textbf{\textit{Keywords---}} #1
}

\begin{document}
\title{Packing a Degree Sequence Realization With A Graph}
\author{James M.\ Shook$^{3,4}$}

\footnotetext[3]{National Institute of Standards and Technology, Computer Security Division, Gaithersburg, MD; {\tt james.shook@nist.gov}.}
\footnotetext[4]{Official Contribution of the National Institute of Standards and Technology; Not subject to copyright in the United States.}
\maketitle
\begin{abstract}Two simple $n$-vertex graphs $G_{1}$ and $G_{2}$, with respective maximum degrees $\Delta_{1}$ and $\Delta_{2}$, are said to pack if $G_{1}$ is isomorphic to a subgraph of the complement of $G_{2}$. The BEC conjecture by Bollob\'{a}s, Eldridge, and Catlin, states that if $(\Delta_{1}+1)(\Delta_{2}+1)\leq n+1$, then $G_{1}$ and $G_{2}$ pack. The BEC conjecture is true when $\Delta_{1}=2$ and has been confirmed for a few other classes of graphs with various conditions on $\Delta_{1}$, $\Delta_{2}$, or $n$. We show that if \[(\Delta_{1}+1)(\Delta_{2}+1)\leq n+\min\{\Delta_{1},\Delta_{2}\},\] then there exists a simple graph with an identical degree sequence as $G_{1}$ that packs with $G_{2}$.  However, except for a few cases, we show that this bound is not sharp. As a consequence of our work, we confirm the BEC conjecture if $G_{1}$ is the vertex disjoint union of a unigraph and a forest $F$ such that either $F$ has at least $\Delta_{2}+1$ components or at most $2\Delta_{2}-1$ edges.
\end{abstract}

\section{Introduction}
All graphs in this paper are finite and have no loops. Let $G=(V(G),E(G))$ be a graph with vertex set $V(G)$ and edge set $E(G)$. We denote the degree of a vertex $x\in V(G)$ by $deg_{G}(x)$, and the maximum and minimum degree of a graph will be represented by $\Delta(G)$ and $\delta(G)$, respectively. For $X\subseteq V(G)$, we let $G[X]$ denote the graph induced by the vertices in $X$. We let $G^{+}$ represent the subgraph of $G$ induced by all the vertices in $G$ with degree at least one, and let $\delta^{+}(G)=\delta(G^{+})$. The notation $\overline{G}$ denotes the complement of $G$. 

Let $k$ be a non-negative integer. The disjoint vertex union of $k$ copies of a graph $G$ will be denoted by $kG$. We use $K^{k}$ for the complete graph on $k$ vertices, and $K^{k,k}$ for the complete balanced bipartite graph. The independent set with $k$ vertices is denoted by $I^{k}$, and the path and cycle with $k$ edges are denoted by $P^{k}$ and $C^{k}$, respectively.

We have organized our paper as follows. First, we present the BEC Conjecture in Section~\ref{BECConjecture} since our work was motivated by it. In Section~\ref{sec:degSequence} we define degree sequence packing and discuss known results. We present our main result at the beginning 
 of Section~\ref{sec:mainTheorem}, but defer the proof until Section~\ref{sec:mainTheoremProof}. In Section~\ref{sec:kfactors}, we demonstrate how our main result can be improved when the graph under consideration has a $k$-factor. We propose a possible next step towards a proof of the BEC Conjecture in Section~\ref{sec:componentWisePacking}, and in Section~\ref{sec:uniGraphs} we confirm the BEC conjecture for unigraphs.
 
\subsection{The BEC Conjecture}\label{BECConjecture}
Let $G_{1}$ and $G_{2}$ be two simple graphs with vertex set $V=\{v_{1},\ldots, v_{n}\}$. If $G_{1}$ is a subgraph of $\overline{G_{2}}$, then we either say $G_{1}$ can be \textit{embedded} into $\overline{G_{2}}$ or $G_{1}$ \textit{packs} with $G_{2}$. Both perspectives can be found in the literature, but in this paper we will favor the idea of graph packing. Specifically, if there is a graph $F$ on $V$ that is isomorphic to $G_{1}$ such that $E(F)\cap E(G_{2})=\emptyset$, then we call $F + G_{2}$ a packing of $G_{1}$ with $G_{2}$.

Graph packing has been extensively studied \cite{Kierstead2009,Yap1988}. This includes a large focus on extremal conditions for when a pair of graphs pack \cite{Sauer1978, Kaul2007, Kaul2022}. Our current work is motivated by the following conjecture of Bollob\'{a}s and Eldridge \cite{Bollobas1978} and  Catlin \cite{Catlin1976}.

\begin{conjecture}[BEC conjecture]\label{conjecture:bollobas_eldridge}
Let $G_{1}$ and $G_{2}$ be simple graphs of order $n$ with $\Delta(G_{i})=\Delta_{i}$ for $i=1,2$.  If \begin{equation}\label{eq:BEC}(\Delta_{1}+1)(\Delta_{2}+1)\leq n+1,\end{equation} then $G_{1}$ and $G_{2}$ pack. 
\end{conjecture}
If true, the BEC conjecture would be a considerable generalization of the Hajnal-Szemer\'{e}di Theorem (See \cite{Hajnal1970}, and see \cite{Kierstead2008} for a short proof). Recall that the Hajnal-Szemer\'{e}di Theorem states that for every positive integer $k$, every simple graph $G$ with $\Delta(G)\leq k$ has a vertex $k$ coloring in which any two color classes differ in size by at most one. Thus, to prove the Hajnal-Szemer\'{e}di Theorem with the BEC Conjecture simply pack $G$ with a $|G|$-vertex graph whose vertices are equitably distributed among $k+1$ complete graphs.



In general, (\ref{eq:BEC}) cannot be improved since there are pairs of graphs with $(\Delta_{1}+1)(\Delta_{2}+1)=n+2$. One such example is when $G_{1}=\Delta_{2}K^{2}$  and $G_{2}=K^{\Delta_{2},\Delta_{2}}$ with $\Delta_{2}$ odd. Another example is when $G_{1}=\Delta_{2}K^{2}$ and $K^{\Delta_{2}+1}\subseteq G_{2}$, or the more general $G_{1}=\Delta_{2}K^{\Delta_{1}+1}\cup K^{\Delta_{1}-1}$ and $G_{2}=\Delta_{1}K^{\Delta_{2}+1}\cup K^{\Delta_{2}-1}$.
 
Aigner and Brandt \cite{Aigner1993} confirmed Conjecture~\ref{conjecture:bollobas_eldridge}  when $\Delta_{1}\leq 2$. For $n$ sufficiently large, Alon and Fisher \cite{Alon1996} gave another proof for $\Delta_{1}\leq 2$ and Csaba\footnote{Csaba states in \cite{Csaba2007} that he has a manuscript of a proof for $\Delta_{1}=4$ when $n$ sufficiently large.}, Shokoufandeh, and Szemer\'{e}di \cite{Csaba2003} proved it when $\Delta_{1}=3$. There are results for specific class of graphs. Wang \cite{Wang1994} confirmed the conjecture when $G_{1}$ is a forest. Csaba proved it when $G_{1}$ is bipartite and $n$ is sufficiently large \cite{Csaba2007}. Bollab\'{a}s, Kostochka, and Nakprasit \cite{Bollobas2008} confirmed Conjecture~\ref{conjecture:bollobas_eldridge} when $G_{1}$ is $d$-degenerate, $n>\Delta_{1}\Delta_{2}$, $\Delta_{1}\geq 40d$, and $\Delta_{2}\geq 215$. Batenburg and Kang \cite{CamesVanBatenburg2018} showed that $G_{1}$ and $G_{2}$ pack if $G_{1}$ has no cycles of length four, six, or eight and $\Delta_{1}\geq 940060$ or $\Delta_{1}\geq \Delta_{2}\geq 27620$. Other work related to Conjecture~\ref{conjecture:bollobas_eldridge} can be found in \cite{Kaul2008, Kostochka2007, Eaton2000}. 

\subsection{Degree Sequence Packing}\label{sec:degSequence}
A finite  sequence of non-negative integers  $\pi=(d_{1},\ldots,d_{n})$ is said to be a degree sequence of a graph $G$ if and only if there is a bijection $\phi:[n]\rightarrow [n]$ such that $deg_{G}(v_{i})=d_{\phi(i)}$ for all $v_{i}\in V(G)$. We call $G$ a realization of $\pi$, and let $\mathcal{R}(\pi)$ be the set of all realizations of $\pi$ that are simple graphs. A degree sequence is positive if all of its terms are at least one. Given a graph $G$, we denote by $\pi(G)$ an arbitrary degree sequence of $G$. There is a large breath of literature on degree sequences \cite{Hakimi1962}.  However, our focus will be on packing realizations of degree sequences.

Given a vertex set $V=\{v_{1},\dots, v_{n}\}$ with degree sequences $\pi=(a_{1},\ldots, a_{n})$ and $\pi'=(b_{1},\dots,b_{n})$, several authors asked if there is a realization $G$ of $\pi$ that packs with a realization $H$ of $\pi'$ such that $deg_{G+ H}(v_{i})=a_{i}+b_{i}$. Work on this question can be found in \cite{Busch2012, Diemunsch2015, Yin2016} and \cite{Gollakota2020} for a special case with trees. 

For a simple graph $G$ and function $f:V(G)\rightarrow\{1,\ldots, n\}$, a spanning subgraph $H$ of $G$ is said to be an $f$-factor if $deg_{H}(x)=f(x)$.  One can see that the $f$-factor $H$ packs with $\overline{G}$ and is a realization of the sequence $(f(v_{1}),\ldots, f(v_{n}))$. On the other hand, if $H$ is a realization of a positive degree sequence $\pi=(d_{1},\ldots,d_{n})$ that packs with $\overline{G}$, then $H$ is an $f$-factor of $G$ where $f(v)=deg_{H}(v)$ for all $v\in V$. There is a substantial amount of literature \cite{Plummer2007} on $f$-factors, but of particular interest is the following theorem of Katerinis and Tsikopoulos that we have recast as a theorem about degree sequence packing. 

\begin{theorem}[\cite{Katerinis2000}]\label{thm:ffactor}If $G_{1}$ and $G_{2}$ are simple graphs of order $n$ with $\Delta(G_i)=\Delta_i$ and $\delta_{1}=\delta(G_{1})$ for $i=1,2$ such that
    \begin{itemize}
    \item $\delta_{1}\geq 1$,
    \item $n\geq \frac{\delta_{1}+\Delta_{1}}{\delta_{1}}(\Delta_{2}+1)$, and
    \item $n>\frac{\delta_{1}+\Delta_{1}}{\delta_{1}}(\delta_{1}+\Delta_{1}-3)$,
\end{itemize} then for any bijection $\phi:[n]\rightarrow [n]$, there is a realization $H$ of $\pi(G_{1})=(d_{1},\ldots,d_{n})$ that packs with $G_{2}$ such that $deg_{H}(v_{i})=d_{\phi(i)}$.
\end{theorem}
Theorem~\ref{thm:ffactor} is very close to a result by Kano and Tokushige in \cite{Kano1992}. However, the theorems have slight differences, and we suspect there may be a way to unify them. We will make use of Theorem~\ref{thm:ffactor} in Section~\ref{sec:kfactors}.

In this paper, we will focus on a less restrictive degree sequence packing definition that resembles the packing definition used in the BEC conjecture. 
\begin{definition} Given simple graphs $G_{1}$ and $G_{2}$ with vertex set $V$, we say a degree sequence $\pi(G_{1})=(d_{1},\ldots,d_{n})$ packs with $G_{2}$ if there is a realization $H$ of $\pi(G_{1})$ that packs with $G_{2}$. 
\end{definition}
Note this definition allows both $G_{1}$ and $G_{2}$ to have trivial components and has no restrictions on the degree of the vertices of $H+G_{2}$. Csaba and Vasarhelyi \cite{Csaba2019}  studied this less strict definition from the perspective of graph embedding, but required every degree of the sequence to be a positive integer. We recast one of their theorems on embedding from the perspective of degree sequence packing below.
\begin{theorem}[Csaba and Vasarhelyi \cite{Csaba2019}]\label{thm:CsabaEmbed} For every $\eta>0$ and $D\in \mathbb{N}$, there exists an $n_{0}=n_{0}(\eta, D)$ such that for all $n>n_{0}$ if $G$ is a graph on $n$ vertices with $\Delta(G)\leq n\bigg(\frac{1}{2}-\eta \bigg )-1$ and $\pi=(d_{1},\ldots,d_{n})$ is a positive degree sequence of length $n$ with $d_{i}\leq D$ for all $i\in [n]$, then $\pi$ packs with $G$.
\end{theorem}
Theorem~\ref{thm:CsabaEmbed} says that for a fixed $\Delta(G_{1})$, if $n$ is large enough, then $\pi(G_{1})$ packs with $G_{2}$; even as $\Delta(G_{2})$ approaches $\frac{n}{2}-1$. Since their proof relies on Szemer\'{e}di’s graph regularity theorem, $n$ is quite large. One will notice that Theorem~\ref{thm:ffactor} gives better bounds in many instances and allows for a more restrictive definition of packing. We will show that $n$ can be made considerably smaller when compared to Theorem~\ref{thm:CsabaEmbed} and still pack at the cost of limiting the size of $\Delta(G_{2})$. Moreover, we allow $|G_{1}^{+}|<|G_{1}|$; which can make packing realizations more difficult since there are less realizations to choose from. We state all of this in our main theorem in Section~\ref{sec:mainTheorem}.

\section{Our Contributions}\label{sec:mainTheorem}
We begin this section with our main result.
\begin{theorem}\label{theorem:BEC-half}If $G_{1}$ and $G_{2}$ are simple graphs of order $n$ with $\Delta(G_i)=\Delta_i$ for $i=1,2$ such that \begin{equation}\label{eq:main}
(\Delta_{1}+1)(\Delta_{2}+1)\leq n+1+\min\{\Delta_{1},\Delta_{2}\}+g(G_{1},G_{2})\end{equation} where $g(G_{1},G_{2})=0$ when $\Delta_{1}=\delta_{+}(G_{1})$ and $g(G_{1},G_{2})=\Delta_{2}-1$ when $\Delta_{1}>\delta_{+}(G_{1})$, then there is a realization of $\pi(G_{1})$ that packs with $G_{2}$ unless one of the following conditions are met:
\begin{enumerate}[label=(F\arabic*),ref=(F\arabic*)]
     \item\label{main:fivecycle} $G_{1}=C^{5}\cup K^{1}$ and $G_{2}=2K^{3}$.
     \item\label{main:complete} $\Delta_{2}\leq\Delta_{1}$, $G_{1}=K^{\Delta_{1}+1}\cup I^{\Delta_{1}\Delta_{2}-1}$, and $G_{2}=\Delta_{1}K^{\Delta_{2}+1}$.
    \item\label{main:matchCompleteBipart} $G_{1}=\Delta_{2}K^{2}$  and $G_{2}=K^{\Delta_{2},\Delta_{2}}$ with $\Delta_{2}$ odd. 
    \item\label{main:matchingCompleteGraph} $G_{1}=\Delta_{2}K^{2}$ and $K^{\Delta_{2}+1}\subseteq G_{2}$.
\end{enumerate}
\begin{proof}
    See Section~\ref{sec:mainTheoremProof}
\end{proof}
\end{theorem}

There are no surprises in \ref{main:fivecycle}, \ref{main:complete}, \ref{main:matchCompleteBipart}, and \ref{main:matchingCompleteGraph}. The graphs in \ref{main:matchCompleteBipart} and \ref{main:matchingCompleteGraph} are often stated as examples of sharpness for Conjecture~\ref{conjecture:bollobas_eldridge}, and the graphs in \ref{main:fivecycle} and \ref{main:complete} are extremal examples for Brooks Theorem \cite{Brooks1941}.

The function $g(G_{1},G_{2})$ is sharp. For instance, $G_{1}$ is regular in \ref{main:fivecycle}, \ref{main:complete}, \ref{main:matchCompleteBipart}, and \ref{main:matchingCompleteGraph}. Moreover, if we allow $\Delta_{1}=\Delta_{2}-1$ in  \ref{main:complete}, then we see that $g(G_{1},G_{2})$ cannot be increased by one when $G_{1}$ is regular. It also cannot be increased by one when $G_{1}$ is not regular since $P^{4}$ does not pack with $K^{3}\cup K^{1}$. There are larger examples in which $g(G_{1},G_{2})$ cannot be increased by two. Take as an example the situation where $G_{1}$ is the complete graph minus an edge and $G_{2}$ is the vertex disjoint union of $\Delta(G_{1})-1$ copies of $K^{\Delta(G_{1})+1}$. Clearly, $G_{1}$ and $G_{2}$ do not pack. 

As noted before in \cite{Bollobas2008}, it is easier to pack two graphs of order $n$ if their maximum degrees are much smaller than $n$. This is also true with degree sequence packing. However, degree sequence packing is also easier when $|G^{+}_{1}|$ is close to $n$. This idea is evident with Theorem~\ref{thm:ffactor} and Theorem~\ref{thm:CsabaEmbed}. Since all of our sharpness examples are small and in some sense extreme, it seems reasonable to suspect that $g(G_{1},G_{2})$ could be larger when $\Delta(G_{2})>\Delta(G_{1})$ and $|G^{+}_{1}|$ is close to $n$. 

To avoid the forbidden pairs of graphs listed in \ref{main:fivecycle}, \ref{main:complete}, \ref{main:matchCompleteBipart}, and \ref{main:matchingCompleteGraph}, Theorem~\ref{theorem:BEC-half} implies the following theorem that resembles Conjecture~\ref{conjecture:bollobas_eldridge}.

\begin{theorem}\label{cor:BECdegree}For simple graphs $G_{1}$ and $G_{2}$ of order $n$ with $\Delta(G_i)=\Delta_i$ for $i=1,2$, if  
\begin{equation}\label{eq:BECdegreeSequence}
    (\Delta_{1}+1)(\Delta_{2}+1)\leq n+\min\{\Delta_{1},\Delta_{2}\}
\end{equation} or \begin{equation}
    (\Delta_{1}+1)(\Delta_{2}+1)\leq n+\min\{\Delta_{1}+\Delta_{2},2\Delta_{2}\}
\end{equation} when $G_{1}^{+}$ is not regular, then $\pi(G_{1})$ packs with $G_{2}$. 
\end{theorem}

\subsection{$k$-factors}\label{sec:kfactors}
For some $H\in \mathcal{R}(\pi(G_{1}))$, if $H^{+}$ has a $k$-factor (a $k$-regular spanning subgraph), then $n$ can be smaller and still imply that $\pi(G_{1})$ packs with $G_{2}$. Moreover, we can require the realization that packs with $G_{2}$ to have a large $k$-regular subgraph. To prove this we use a combination of Theorem~\ref{thm:ffactor}, Theorem~\ref{theorem:BEC-half}, and Theorem~\ref{cor:BECdegree}.

\begin{theorem}\label{thm:regularKReduction}For simple graphs $G_{1}$ and $G_{2}$ of order $n$ with $\Delta(G_i)=\Delta_i$ for $i=1,2$ such that $G^{+}_{1}$ is $k$-regular, if $|G_{1}^{+}|\geq 2\Delta_{2}+2(k-1)$,
then for any $X\subseteq V$ with $|X|=|G_{1}^{+}|$ there is a realization $H$ of $\pi(G^{+}_{1})$ that packs with $G_{2}[X]$ unless $G_{1}^{+}$ and $G_{2}[X]$ satisfy \ref{main:matchCompleteBipart} or \ref{main:matchingCompleteGraph} when $|G_{1}^{+}|=2\Delta_{2}$.
\begin{proof}If $k=1$, then the theorem follows from a direct application of Theorem~\ref{theorem:BEC-half} since in this case  $G_{1}^{+}$ and $G_{2}[X]$ do not satisfy \ref{main:fivecycle}, \ref{main:complete}, \ref{main:matchCompleteBipart}, nor \ref{main:matchingCompleteGraph}. Thus, we assume $k\geq 2$. 

For $\Delta_{2}=1$, $|G_{1}^{+}|\geq 2\Delta_{2}+2(k-1)\geq 2k$ and $G_{1}^{+}$ and $G_{2}[X]$ do not satisfy \ref{main:fivecycle}, \ref{main:complete}, \ref{main:matchCompleteBipart}, nor \ref{main:matchingCompleteGraph}. Therefore, since $|G^{+}_{1}|$ satisfies (\ref{eq:main}), Theorem~\ref{theorem:BEC-half} says $\pi(G_{1}^{+})$ packs with $G_{2}[X]$. Thus, we may assume $\Delta_{2}\geq 2$. This implies $|G_{1}^{+}|\geq 2\Delta_{2}+2(k-1)\geq 4k-5$ when $k\leq \Delta_{2}+1$. Therefore, Theorem~\ref{thm:ffactor} says $\pi(G_{1}^{+})$ packs with $G_{2}[X]$ when $k\leq \Delta_{2}+1$. Let us now assume $k\geq \Delta_{2}+2\geq 4$. Let $r$ be the largest non-negative integer such that there is a graph $H_{r}$ that is the edge-disjoint union of $G_{2}[X]$ and an $r$-factor. Suppose $r<k$, and let $k'=\min\{3,k-r\}$. Thus, $|G_{1}^{+}|\geq 2\Delta_{2}+2(k-1)\geq 2(\Delta_{2}+r)$, and since $k'\leq 3\leq \Delta_{2}+1$, we see that $|G_{1}^{+}|\geq 4k'-5$. Therefore, by Theorem~\ref{thm:ffactor}, a $k'$-factor with $|G_{1}^{+}|$ vertices packs with $H_{r}[X]$. However, this contradicts our choice of $r$. Thus, $\pi(G^{+}_{1})$ packs with $G_{2}[X]$.
\end{proof}
\end{theorem}

If $|G^{+}_{1}|$ is even, then it is straight forward to modify the proof of Theorem~\ref{thm:regularKReduction} to show that there is a realization of $\pi(G^{+}_{1})$ with $k$ edge-disjoint $1$-factors that packs with $G_{2}$. 

\begin{corollary}\label{cor:disjoint1Factors}For simple $n$-vertex graphs $G_{1}$ and $G_{2}$ with $\Delta(G_i)=\Delta_i$ for $i=1,2$ such that $G^{+}_{1}$ is $k$-regular, if $|G^{+}_{1}|$ is even such that $|G_{1}^{+}|\geq 2\Delta_{2}+2(k-1)$,
then for any $X\subseteq V$ with $|X|=|G_{1}^{+}|$ there is a realization $H$ of $\pi(G^{+}_{1})$ with $k$ edge-disjoint $1$-factors that packs with $G_{2}[X]$ unless $G_{1}^{+}$ and $G_{2}[X]$ satisfy \ref{main:matchCompleteBipart} or \ref{main:matchingCompleteGraph} when $|G_{1}^{+}|=2\Delta_{2}$.
\end{corollary}

\begin{theorem}\label{thm:notRegularKReduction}For simple $n$-vertex graphs $G_{1}$ and $G_{2}$ with $\Delta(G_i)=\Delta_i$ for $i=1,2$ such that $|G^{+}_{1}|\geq 2(\Delta_{1}+\Delta_{2})-1$, if \[(\Delta_{1}+1)(\Delta_{2}+1)\leq n+\min\{\Delta_{1}+k\Delta_{2}, k+(k+1)\Delta_{2}\}\] and some realization of $\pi(G^{+}_{1})$ has a $k$-factor, then there is a realization $H$ of $\pi(G_{1})$ that packs with $G_{2}$ such that has a $H^{+}$ has a $k$-factor.
\begin{proof}Let $W$ be a realization of $\pi(G^{+}_{1})$ with a $k$-factor $F$. Let $H_{1}=W-F$. Clearly, $\Delta(H_{1})=\Delta_{1}-k$. Since $(\Delta_{1}-k+1)(\Delta_{2}+1)\leq n+\min\{\Delta_{1}-k,\Delta_{2}\}$, Theorem~\ref{cor:BECdegree} says there is a realization $H_{2}$ of $\pi(H_{1})$ that packs with $G_{2}$. Let $Z=H_{2}+ G_{2}$, and let $X=V(H^{+}_{2})$. Since $\Delta(Z)\geq \Delta_{2}+\Delta_{1}-k$ and $|F|\geq 2(\Delta_{1}+\Delta_{2})-1$, Theorem~\ref{thm:regularKReduction} says $\pi(F)$ packs with $Z[X]$. Thus, there is a realization $H$ of $\pi(G_{1})$ that packs with $G_{2}$ such that $H^{+}$ has a $k$-factor with at least $2(\Delta_{1}+\Delta_{2})-1$ vertices.
\end{proof}
\end{theorem}

In \cite{Shook2022}, we showed that if $G_{1}$ has at least $\Delta_{1}-\delta^{+}+1$ vertices with degree at least $\Delta_{1}-\delta^{+}+k-1$, then some realization of $\pi(G_{1})$ has a $k$-factor. Thus, if $|G^{+}|\geq 2(\Delta_{1}+\Delta_{2})-1$ and $\pi(G_{1})$ doesn't have a $k$-factor, then $G_{1}$ has at most $\Delta_{1}-\delta^{+}$ vertices with degree at least $\Delta_{1}-\delta^{+}+k-1$. Perhaps with this information and Theorem~\ref{thm:notRegularKReduction} one could improve Theorem~\ref{theorem:BEC-half} when $|G^{+}|\geq 2(\Delta_{1}+\Delta_{2})-1$.

\subsection{Component-Wise Packing}\label{sec:componentWisePacking}
We denote by $\omega(G)$ the number of components in a graph $G$.

\begin{definition}Given graphs $G_{1}$ and $G_{2}$ such that $G_{1}$ has components $C_{1},\ldots, C_{\omega(G_{1})}$, we say $\pi(G_{1})$ component-wise packs with $G_{2}$ if there is a realization $H$ of $\pi(G_{1})$ that packs with $G_{2}$ such that $H$ has components $D_{1},\ldots, D_{\omega(G_{1})}$ where $D_{i}\in \mathcal{R}(C_{i})$.
\end{definition}

Considering the gap between Theorem~\ref{theorem:BEC-half} and the BEC conjecture it seems like a reasonable next step to pose the following problem.

\begin{problem}\label{problem1}If simple graphs $G_{1}$ and $G_{2}$ satisfy $(\Delta(G_{1})+1)(\Delta(G_{2})+1)\leq n+1$, then does $G_{1}$ component-wise pack with $G_{2}$?
\end{problem}

Answering Problem~\ref{problem1} in the negative would disprove the BEC conjecture. The example $G_{1}=\Delta_{2}K^{\Delta_{1}+1}\cup K^{\Delta_{1}-1}$ and $G_{2}=\Delta_{1}K^{\Delta_{2}+1}\cup K^{\Delta_{2}-1}$ shows that the bound in Problem~\ref{problem1} is sharp. 

\begin{theorem}[\cite{Aigner1993}]\label{thm:BECmax2} Let $G_{1}$ and $G_{2}$ be simple graphs of order $n$ with $\Delta(G_{i})=\Delta_{i}$ for $i=1,2$ such that $\Delta_{1}\leq 2$.  If $(\Delta_{1}+1)(\Delta_{2}+1)\leq n+1$, then $G_{1}$ and $G_{2}$ pack.
\end{theorem}

Thus, we only need to consider the problem for when both $\Delta_{1}$ and $\Delta_{2}$ are at least three. Applying the Hajnal-Semeredi Theorem to $G_{2}$ gives us a positive answer to the question when $|C_{i}|\leq \lfloor n/(\Delta_{2}+1)\rfloor$ for each component $C_{i}$ of $G_{1}$. We also know that $\pi(G_{1})$ component-wise packs with $G_{2}$ when $G_{1}$ is a forest.

\begin{theorem}[\cite{Wang1994}]\label{thm:BECForest} Let $G_{1}$ and $G_{2}$ be graphs with $n$ vertices such that $G_{1}$ is a forest. If $\Delta_{2}(\Delta(G_{1})+1)\leq n$ where $\Delta_{2}=\Delta(G_{2})$, then there is a packing of $G_{1}$ and $G_{2}$ unless $G_{1}=K^{1,n-1}$ and $G_{2}=\frac{n}{2}K^{2}$, $G_{1}=\Delta_{2}K^{2}$ and  $K^{\Delta_{2}+1}\subseteq G_{2}$, or $G_{1}=\Delta_{2}K^{2}$ and $G_{2}=K^{\Delta_{2},\Delta_{2}}$ with $\Delta_{2}$ odd.
\end{theorem}

So, what can we say if $G_{1}$ is the disjoint union of a graph and a forest? We start by stating a useful theorem of Golberg and Magdon-Ismail.

\begin{theorem}[\cite{Goldberg2011}]\label{thm:ForestPacking}
    If $F=T_{1}\cup \ldots \cup T_{p}$ is a forest, then every graph $G$ with at least $|F|$ vertices and minimum degree at least $|E(F)|$ contains $F$ as a subgraph.
\end{theorem}

\begin{lemma}\label{lem:ForestComponentWise}Let $G_{1}$ and $G_{2}$ be graphs on $V$ such that $\pi(G_{1})$ component-wise packs with $G_{2}$. If $F$ is a forest with $\omega(F)\geq \Delta(G_{2})+1$ such that $|F|+|G_{1}|=|V|$, then $\pi(G_{1}\cup F)$ component-wise packs with $G_{2}$.
\begin{proof}Since $\pi(G_{1})$ component-wise packs with $G_{2}$, there is a realization $H$ of $\pi(G_{1})$ whose components are realizations of components of $G_{1}$. Let $X$ be the set of vertices in $V$ such that $deg_{H}(x)=0$ for all $x\in X$. Thus, $|X|=|F|$. Since \[\delta(\overline{G_{2}[X]})\geq |F|-1-\Delta(G_{2}[x])\geq |F|-1-\Delta(G_{2})=|E(F)|+\omega(F)-(\Delta_{2}+1)\geq |E(F)|,\] Theorem~\ref{thm:ForestPacking} says $F$ is a subgraph of $\overline{G_{2}[X]}$. Thus, $F$ packs with $G_{2}[X]$, and therefore, $\pi(G_{1}\cup F)$ component-wise packs with $G_{2}$. 
\end{proof}
\end{lemma}

Perhaps a good starting point for Problem~\ref{problem1} is to prove Lemma~\ref{lem:ForestComponentWise} without the lower bound on $\omega(F)$.  Next we prove a partial result of Problem~\ref{problem1} that we will use in the next section. 

Let $G$ be a graph.  For $A\subseteq V(G)$ and $B\subseteq V(G)$, we let $E_{G}(A, B)$ be the set of edges with one end in $A$ and the other end in $B$, and let $e_{G}(A,B)=|E_{G}(A, B)|$. A set of vertices $X\subseteq V(G)$ is said to be \textit{dominating} if every vertex in $V-X$ is adjacent to a vertex in $X$. We say $G$ has a dominating clique if there is a dominating set $X$ such that $G[X]$ is a clique.


\begin{theorem}\label{thm:DominatingAndForest}Let $G_{1}$ and $G_{2}$ be graphs with vertex set $V$ with $\Delta(G_i)=\Delta_i$ for $i=1,2$ such that $G_{1}$ can be partitioned into a graph $C_{1}$ that has a dominating clique and a forest $F$ such that either $\omega(F)\geq \Delta(G_{2})+1$ or $|E(F)|\leq 2\Delta(G_{2})-1$. If $G_{1}$ and $G_{2}$ satisfy $(\Delta_{1}+1)(\Delta_{2}+1)\leq n+1$, then $\pi(G_{1})$ component-wise packs with $G_{2}$.
\begin{proof}By Theorem~\ref{thm:BECmax2}, we may assume both $\Delta(G_{1})$ and $\Delta(G_{2})$ are at least three. By Theorem~\ref{cor:BECdegree}, $\pi(C_{1})$ packs with $G_{2}$. Lemma~\ref{lem:ForestComponentWise} says $\pi(G_{1})$ component-wise packs with $G_{2}$ if $\omega(F)\geq \Delta(G_{2})+1$. Thus, we will assume $\omega(F)\leq \Delta(G_{2})$ for the rest of the proof. This implies $|F|=\omega(F)+|E(F)|\leq  3\Delta(G_{2})-1$. Thus, \[|C_{1}|=n-|F|\geq (\Delta(G_{1})+1)(\Delta(G_{2})+1)-1-(3\Delta(G_{2})-1)=(\Delta(G_{1})-2)\Delta(G_{2})+\Delta(G_{1})+1.\]

By Theorem~\ref{cor:BECdegree}, $K^{\Delta_{1}+1}$ packs with $G_{2}$. Thus, there is an independent set $X$ of $G_{2}$ with $\Delta(G_{1})+1$ vertices. If $|F|\leq \Delta_{1}+1$, then we may pack $F$ with $G_{2}[X]$. This implies $|C_{1}|=n-|F|\geq (\Delta(G_{1})+1)\Delta_{2}-1$. Since $C_{1}$ is connected and has more than $\Delta(G_{1})+1$ vertices, we may conclude that $C_{1}$ and $G_{2}-X$ do not satisfy \ref{main:fivecycle}, \ref{main:complete},  \ref{main:matchCompleteBipart}, nor \ref{main:matchingCompleteGraph}. Thus, Theorem~\ref{theorem:BEC-half} says $\pi(C_{1})$ packs with $G_{2}-X$. Therefore, $\pi(G_{1})$ component-wise packs with $G_{2}$. Thus, we are left with the case $3\Delta(G_{2})-1\geq |F|\geq \Delta(G_{1})+2$. 

Let $A$ be a dominating clique of $C_{1}$. Since every vertex in $V-F-|A|$ is adjacent in $C_{1}$ to a vertex in $A$, we see that $n-|F|-|A|\leq e_{C_{1}}(A,V-A-V(F))\leq |A|(\Delta(G_{1})+1-|A|)$.  Thus, \[(\Delta(G_{1})-2)\Delta(G_{2})+\Delta(G_{1})+1\leq n-|F|\leq |A|(\Delta(G_{1})+2-|A|).\] Since $|A|(\Delta(G_{1})+2-|A|)$ is maximized when $|A|=\frac{\Delta_{1}+2}{2}$, we deduce that $(\Delta(G_{1})-2)\Delta(G_{2})+\Delta(G_{1})+1\leq \frac{(\Delta_{1}+2)^{2}}{4}$. Expanding the right hand side of the last inequality and removing like terms we have the arrangement $4(\Delta(G_{1})-2)\Delta(G_{2})\leq \Delta(G_{1})^{2}.$ Since $\min\{\Delta(G_{1}),\Delta(G_{2})\}\geq 3$, we see a contradiction when $\Delta(G_{1})=3$. Therefore, we may assume $\Delta(G_{1})\geq 4$. Since $3\Delta(G_{2})-1\geq |F|\geq \Delta(G_{1})+2$, we deduce that $4\Delta(G_{2})\geq \Delta(G_{1})+3+\Delta(G_{2})\geq \Delta(G_{1})+6$. Therefore, we have the contradiction \begin{align}\Delta(G_{1})^{2}&\geq 4(\Delta(G_{1})-2)\Delta(G_{2})\notag\\ &\geq (\Delta(G_{1})-2)(\Delta(G_{1})+6)\notag\\
&=\Delta(G_{1})^{2}+4\Delta(G_{1})-12\geq \Delta(G_{1})^{2}+4\qedhere \notag
\end{align}
\end{proof}
\end{theorem}

The bounds and arguments we use in the proof of Theorem~\ref{thm:DominatingAndForest} are almost certainly not optimal, but the proof demonstrates how one can use our results when thinking about Problem~\ref{problem1}.

\subsection{Unigraphs}\label{sec:uniGraphs}

\begin{figure}[h]
    \centering

    \begin{tikzpicture}
\draw[color=black] (1.25,.5) -- (2.25,2.5);
\draw[color=black] (1.25,.5) -- (.25,1.25);
\draw[color=black] (.25,1.25) -- (2.25,2.5);

\draw[color=black] (3.25,.5) -- (2.25,2.5);
\draw[color=black] (3.25,.5) -- (4.25,1.25);
\draw[color=black] (4.25,1.25) -- (2.25,2.5);

\draw[color=black] (2.25,2.5) -- (3.25,3.5);
\draw[color=black] (2.25,2.5) -- (1.25,3.5);
\draw[color=black] (1.25,3.5) -- (2.25,4.5);
\draw[color=black] (3.25,3.5) -- (2.25,4.5);

\draw[fill=black] (.25,1.25) circle (.1cm);
\draw[fill=black] (4.25,1.25) circle (.1cm);
\draw[fill=black] (1.25,.5) circle (.1cm);
\draw[fill=black] (2,.5) circle (.025cm);
\draw[fill=black] (2.25,.5) circle (.025cm);
\draw[fill=black] (2.5,.5) circle (.025cm);

\draw[fill=black] (1.25,3.5) circle (.1cm);
\draw[fill=black] (3.25,3.5) circle (.1cm);
\draw[fill=black] (2.25,4.5) circle (.1cm);
\draw[fill=black] (2.25,2.5) circle (.1cm);
\draw[fill=black] (3.25,.5) circle (.1cm);

\draw [decorate,
    decoration = {calligraphic brace,mirror,amplitude=5pt}] (1,.25) --  (3.5,.25) node[pos=0.5,below=5pt,black]{$l$};

\end{tikzpicture}
    \caption{$U_{3}(l)$ consists of a $C^{4}$ sharing exactly one vertex with $l$ cycles with three vertices.}
    \label{fig:U3l}
\end{figure}
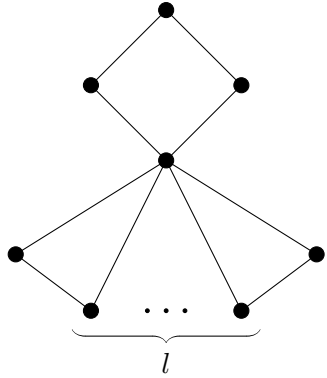

A graph $H$ is said to be a unigraph if $H$ is isomorphic to every realization of $\pi(H)$. This definition implies the complement of a unigraph is a unigraph. Some small examples of unigraphs include  $K^{n}$, $P^{4}$, $C^{4}$, and $C^{5}$.  More complex examples are $U_{2}(l,t)$ and $U_{3}(l)$. Where $U_{2}(l,t)=lK^{2}\cup K^{1,t}$ with $l\geq 1$ and $t\geq 1$, and $U_{3}(l)$ with $l\geq 0$ (See Figure~\ref{fig:U3l}) is a graph consisting of $l$ copies of $K^{3}$ that all mutually share exactly one vertex with $C_{4}$. All these examples along with a complete characterization of unigraphs can be found in a paper of  Tyshkevich \cite{Tyshkevich2000}. Tyshkevich's characterization shows that if $G$ is a connected unigraph that does not have a dominating clique, then $G$ is isomorphic to either $C^{5}$ or $U_{3}(l)$. Otherwise, if $G$ is not connected, then $G$ is isomorphic to $U_{2}(l,t)$ or the disjoint union of a unigraph and an independent set.

We are now ready to prove our main theorem of this section.

\begin{theorem}\label{thm:BECUnigraph}For simple graphs $G_{1}$ and $G_{2}$ of order $n$ with $\Delta(G_i)=\Delta_i$ for $i=1,2$, if $(\Delta_{1}+1)(\Delta_{2}+1)\leq n+1$ and $G_{1}$ is the vertex disjoint union of a unigraph $H$ and a forest $F$ with either $\omega(F)\geq \Delta_{2}+1$ or $|E(F)|\leq 2\Delta_{2}-1$, then $G_{1}$ packs with $G_{2}$.
\begin{proof}If $G_{1}$ is a forest, then Theorem~\ref{thm:BECForest} says $G_{1}$ and $G_{2}$ pack.
Since disconnected unigraphs can be partitioned into a connected unigraph and a forest, we may assume $H$ is a connected unigraph. If $H$ has a dominating clique, then the theorem follows from Theorem~\ref{thm:DominatingAndForest}. If $H$ does not have a dominating clique, then $H$ is isomorphic to either $C^{5}$ or $U_{3}(l)$. Therefore, $|H|\leq \Delta(H)+3$. This implies $|F|=n-|H|\geq (\Delta_{1}+1)(\Delta_{2}+1)-1-(\Delta(H)+3)\geq (\Delta_{1}+1)\Delta_{2}-2$. If $\omega(F)\leq \Delta_{2}$, then we have the contradiction \[|E(F)|=|F|-\omega(F)\geq (\Delta_{1}+1)\Delta_{2}-2-\Delta(H)=\Delta_{1}\Delta_{2}-2=2\Delta_{2}+\Delta(H)-2\geq 2\Delta_{2}+1.\]
If $\omega(F)\geq \Delta_{2}+1$, then Lemma~\ref{lem:ForestComponentWise} says $\pi(G_{1})$ component wise packs with $G_{2}$, and therefore, $G_{1}$ packs with $G_{2}$ since $H$ is a unigraph.  Thus, $G_{1}$ packs with $G_{2}$. 
\end{proof}
\end{theorem}

We finish with a simple observation about Theorem~\ref{thm:BECUnigraph}.

\begin{corollary}For simple graphs $G_{1}$ and $G_{2}$ of order $n$ with $\Delta(G_i)=\Delta_i$ for $i=1,2$, if $G_{1}$ is the vertex disjoint union of a unigraph $H$ and a matching $M$ such that  $(\Delta_{1}+1)(\Delta_{2}+1)\leq n+1$, then $G_{1}$ packs with $G_{2}$.
\begin{proof}The corollary follows directly from Theorem~\ref{thm:BECUnigraph} since $|E(M)|=|M|/2\leq \Delta_{2}\leq 2\Delta_{2}-1$ when $|M|/2=\omega(F)\leq \Delta_{2}$.   
\end{proof}
\end{corollary}
\section{\label{sec:mainTheoremProof}Proof of Theorem~\ref{theorem:BEC-half}}
Before we present the proof we need to take a moment to define some notation and two important concepts in Section~\ref{ExchInter}. Refer to \cite{Diestel2016} for any definitions and notation not explicitly defined in the proof. Let $G$ be a simple graph with vertex set $V$. For a vertex $u\in V$, we let $N_{G}(u)$ be the neighbors of $u$ in $G$ and $N_{G}[u]=N_{G}(u)\cup \{u\}$. If $uv$ is an edge of $\overline{G}$, then we say $uv$ is a non-edge of $G$. Moreover, we let $\alpha(G)$ be the order of a largest independent set in $G$ and $\chi(G)$ be the chromatic number of $G$. Later we will need the well-known inequality $\alpha(G)\geq \frac{n}{\chi(G)}$ whose proof can be often found as an exercise when discussing Brook's Theorem. 

\begin{theorem}[Brook's Theorem \cite{Brooks1941}]Every graph $G$ can be colored with $\Delta(G)$ colors unless $G$ contains a $K^{\Delta(G)+1}$ or $\Delta(G)=2$ and $G$ contains an odd cycle.
\end{theorem}

Given an additional simple graph $H$ with vertex set $V$, if there is an edge $xy\in E(G)\cap E(H)$, then we call $xy$ a $(G,H)$-bad pair. We denote the number of $(G,H)$-bad pairs by $b(G,H)$.

\subsection{Edge Exchanges and Vertex Interchanges}\label{ExchInter}
\begin{figure}[ht]
    \centering
    \begin{tikzpicture}
\draw[color=red] (1,3) -- (2,4);
\draw[color=red] (3,3) -- (2,2);
\draw[color=red] (7,3) -- (8,2);
\draw[color=red] (9,3) -- (8,4);
\draw[fill=black] (1,3) circle (.1cm);
\draw[fill=black] (3,3) circle (.1cm);
\draw[fill=black] (2,2) circle (.1cm);
\draw[fill=black] (2,4) circle (.1cm);

\draw[fill=black] (7,3) circle (.1cm);
\draw[fill=black] (9,3) circle (.1cm);
\draw[fill=black] (8,2) circle (.1cm);
\draw[fill=black] (8,4) circle (.1cm);
\node[] at (.5,3) {$x$};
\node[] at (3.5,3) {$y$};
\node[] at (2.5,4) {$v$};
\node[] at (2.5,2) {$u$};

\node[] at (6.5,3) {$x$};
\node[] at (9.5,3) {$y$};
\node[] at (8.5,4) {$v$};
\node[] at (8.5,2) {$u$};

\end{tikzpicture}
    \caption{An edge-exchange of the edge $vx$ and $uy$ with the non-edges $vy$ and $ux$.}
    \label{fig:DefEdgeExchanges}
\end{figure}
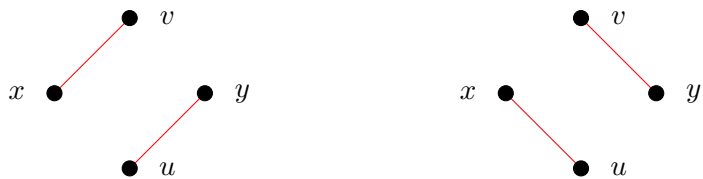

\begin{figure}[ht]
    \centering
    \begin{tikzpicture}
\draw[color=red] (3,4) -- (2,2);
\draw[color=red] (3,4) -- (3,2);
\draw[color=red] (3,4) -- (4,2);
\draw[color=red] (5,4) -- (4,2);
\draw[color=red] (5,4) -- (6,2);

\draw[color=red] (9,4) -- (10,2);
\draw[color=red] (9,4) -- (12,2);
\draw[color=red] (11,4) -- (8,2);
\draw[color=red] (11,4) -- (9,2);
\draw[color=red] (11,4) -- (10,2);

\draw[fill=black] (3,4) circle (.1cm);
\draw[fill=black] (3,2) circle (.1cm);
\draw[fill=black] (5,4) circle (.1cm);
\draw[fill=black] (2,2) circle (.1cm);
\draw[fill=black] (4,2) circle (.1cm);
\draw[fill=black] (6,2) circle (.1cm);

\draw[fill=black] (9,4) circle (.1cm);
\draw[fill=black] (9,2) circle (.1cm);
\draw[fill=black] (11,4) circle (.1cm);
\draw[fill=black] (8,2) circle (.1cm);
\draw[fill=black] (10,2) circle (.1cm);
\draw[fill=black] (12,2) circle (.1cm);

\node[] at (3.4,4) {$u$};
\node[] at (5.4,4) {$v$};
\node[] at (9.4,4) {$u$};
\node[] at (11.4,4) {$v$};

\end{tikzpicture}
    \caption{An interchange of $u$ with $v$.}
    \label{fig:DefInterchange}
\end{figure}
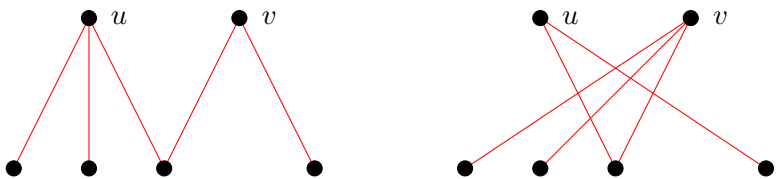

We will often need to construct a new realization of $\pi(G)$ from a given one that respects the labeling of the vertices in $V$. We do this in two ways. The first method is the traditional {\it edge-exchange} (See Figure~\ref{fig:DefEdgeExchanges}). This simply involves exchanging the edges $vx$ and $uy$ of $G$ with the non-edges $vy$ and $ux$ of $G$ to construct a new realization in $\mathcal{R}(\pi)$. The second method is thru vertex interchanges (See Figure~\ref{fig:DefInterchange}). Given two vertices $u$ and $v$ of a $G$ we {\it interchange} $u$ with $v$ by constructing a new graph $H\in \mathcal{R}(\pi)$ such that $N_{H}(u)=N_{G}(v)$, $N_{H}(v)=N_{G}(u)$, and $N_{H}(x)=N_{G}(x)$ for all $x\in V -\{u,v\}$. Note that $H=G$ when $u=v$.

\subsection{Beginning of the Proof}
\begin{theoremlateproof}{theorem:BEC-half}
Let $G_{2}$ be an arbitrary graph with vertex set $V=\{v_{1},\ldots,v_{n}\}$. Suppose there is a graph $G_{1}$ on $V$ with the fewest positive degree vertices such that $\pi(G_{1})$ does not pack with $G_{2}$, but $G_{1}$ and $G_{2}$ satisfy (\ref{eq:main}) and do not satisfy \ref{main:fivecycle},\ref{main:complete}, \ref{main:matchCompleteBipart}, nor \ref{main:matchingCompleteGraph}. 

\begin{remark}When considering the multi-graph $H+G_{2}$ for some $H\in \mathcal{R}(\pi(G_{1}))$ we like to color the edges of $H$ red and the edges of $G_{2}$ blue.
\end{remark}

We let $\Delta_{1}=\Delta(G_{1})$, $\Delta_{2}=\Delta_{2}(G_{2})$, and $\delta_{+}=\delta_{+}(G_{1})$.

For the rest of this section, we tackle the two critical cases of when $|G_{1}^{+}|\leq \Delta_{1}+1$ or $|G_{2}^{+}|\leq \Delta_{2}+1$. 

\begin{remark}Claim~\ref{cl:independenceBound} below is redundant since it is proved implicitly throughout the rest of the proof. However, we feel the reader would benefit by seeing such a critical case directly addressed first.
\end{remark}

\begin{claim}\label{cl:independenceBound}There exists a realization $H$ of $\pi(G_{1})$ such that $\alpha(H)\geq  \Delta_{2}+1$.
\begin{proof}Let $H$ be a realization of $\pi(G_{1})$ with the fewest components that are either a cycle or a $K^{\Delta_{1}+1}$.  From (\ref{eq:main}), we see that \begin{align}
n&\geq (\Delta_{1}+1)(\Delta_{2}+1)-1-\min\{\Delta_{1},\Delta_{2}\}-g(G_{1},G_{2})\notag \\
&= \Delta_{1}\Delta_{2}+\max\{\Delta_{1},\Delta_{2}\}-g(G_{1},G_{2}).\notag
\end{align}
If $H$ does not contain a $K^{\Delta_{1}+1}$ nor an odd cycle when $\Delta_{1}=2$, then by Brooks Theorem  $\chi(H)\leq \Delta_{1}$. This implies $\alpha(H)\geq \lceil n/\chi(G) \rceil \geq \lceil n/\Delta_{1} \rceil\geq \lceil(\Delta_{1}\Delta_{2}+1)/\Delta_{1}\rceil\geq \Delta_{2}+1$.

Suppose there is a $C\subseteq V$ such that $H[C]$ is a $K^{\Delta_{1}+1}$ or an odd cycle. If there is an edge $uu'$ such that $\{u,u'\}\subseteq V-C$, then we choose an edge $vv'$ of $H[C]$ and then exchange the edges $uu'$ and $vv'$ with the non-edges $uv$ and $u'v'$ of $H$ to combine two components of $H$. This creates a realization of $\pi(G_{1})$ with fewer components than $H$ that are either a $K^{\Delta_{1}+1}$ or a cycle. Thus, it must be the case that $V-C$ is an independent set in $H$. Therefore, $H^{+}$ is regular and $g(H,G_{2})=0$. If $n=\Delta_{1}+1$, then $\Delta_{1}=\Delta_{2}=1$. However, this implies both $G_{1}$ and $G_{2}$ have exactly one edge, and therefore, contradicts our claim that $G_{1}$ and $G_{2}$ do not satisfy \ref{main:matchingCompleteGraph}. Thus, $n>\Delta_{1}+1$, and therefore, $\alpha(H)\geq 2$. Thus, we may assume $\Delta_{2}\geq 2$ since any non-edge of $H$ takes care of the $\Delta_{2}=1$ case. If $H[C]=K^{\Delta_{1}+1}$, then 
 \begin{align}\alpha(H)&= \alpha(H[C])+n-(\Delta_{1}+1)\notag\\
 &\geq 1+\Delta_{1}\Delta_{2}+\max\{\Delta_{1},\Delta_{2}\}-(\Delta_{1}+1)\notag\\
 &=\Delta_{1}(\Delta_{2}-1)+\max\{\Delta_{1},\Delta_{2}\}\notag\\
 &\geq \Delta_{2}+1.\notag
 \end{align}
So, $C$ must be an odd cycle with at least five vertices, and therefore, $\alpha(H[C])=(|C|-1)/2$. Since $\Delta_{1}=2$ and $\Delta_{2}\geq 2$, we may deduce that $n\geq 2(\Delta_{2}+1)$. Therefore, we can conclude that
\begin{align}
    \alpha(H)&=\alpha(H[C])+n-|C|\notag\\
    &= (|C|-1)/2+n-|C|= n-(|C|+1)/2\notag\\
    &\geq \big\lceil(n-1)/2\big\rceil\notag\\
    &\geq \big\lceil(2\Delta_{2}+1)/2\big\rceil=\Delta_{2}+1.\notag
\end{align}
This completes the proof of Claim~\ref{cl:independenceBound}.
\end{proof}
\end{claim}

From Claim~\ref{cl:independenceBound}, we can reason that if $G^{+}_{2}$ has at most $\Delta_{2}+1$ vertices, then we simply find a realization $H$ of $\pi(G_{1})$ with an independent set of size $|G^{+}_{2}|$ on the same vertices as $G^{+}_{2}$ to construct a packing of $\pi(G_{1})$ with $G_{2}$.

\begin{claim}\label{cl:minOrder}$|V(G_{1}^{+})|\geq \Delta_{1}+2$.
\begin{proof} 
Suppose $|V(G_{1}^{+})|\leq \Delta_{1}+1$. Let $S$ be a maximum independent set in $G_{2}$. If $\alpha(G_{2})\geq \Delta_{1}+1$, then there is a contradiction since we can simply assign the vertices of $G_{1}^{+}$ into $S$ and arbitrarily assign the zero degree vertices of $G_{1}$ to the vertices of $V-S$ to construct a packing of $G_{1}$ with $G_{2}$. Suppose $\alpha(G_{2})\leq \Delta_{1}$. Let $J=\{u_{1},\ldots, u_{|J|}\}$ be a largest independent set in $G_{2}$ such that $N_{G_{2}}(u_{i})\cap N_{G_{2}}(u_{j})=\emptyset$ for any two distinct vertices $u_{i}$ and $u_{j}$ of $J$. Let $W=\bigcup_{u_{i}\in J}N_{G_{2}}(u_{i})$ and $X=V-J-W$. Since $J$ is maximal, we may conclude that $N_{G_{2}}(x)\cap W\neq \emptyset$ for every $x\in X$. Note that $\alpha(G_{2})\geq |J|+\alpha(G_{2}[X])$. From (\ref{eq:main}), we deduce that
\begin{align}
|J|(\Delta_{2}+1)+|X|&\geq \sum_{v_{i}\in J}(deg_{G_{2}}(v_{i})+1)+|X|=n\notag\\
&\geq (\Delta_{1}+1)(\Delta_{2}+1)-g(G_{1},G_{2})-(\min\{\Delta_{1},\Delta_{2}\}+1).\label{eq:minG1}
\end{align}

Consider the case $\Delta_{1}=\delta_{+}$, and therefore, $g(G_{1},G_{2})=0$. First, suppose $X$ is not empty. In this case $\Delta(G_{2}[X])\leq \Delta_{2}-1$ and $\alpha(G_{2}[X])\geq 1$. Therefore, if $|J|=\Delta_{1}$, then $\alpha(G_{2})\geq \Delta_{1}+1$. We now assume $|J|<\Delta_{1}$. From (\ref{eq:minG1}), we see that $|X|\geq (\Delta_{1}+1-|J|)(\Delta_{2}+1)-(\Delta_{2}+1)\geq (\Delta_{1}-|J|)(\Delta_{2}+1)$. Therefore, $\alpha(G_{2}[X])\geq \frac{|X|}{\Delta_{2}}>\Delta_{1}-|J|$. However, this gives us a contradiction since $\alpha(G_{2})\geq |J|+\alpha(G_{2}[X])\geq \Delta_{1}+1$. Now suppose $X$ is empty. Since $\Delta_{1}\geq |J|$, we see from (\ref{eq:minG1}) that $\Delta_{2}\leq \Delta_{1}$, $|J|=\Delta_{1}$, $n=\Delta_{1}(\Delta_{2}+1)$, and $deg_{G_{2}}(u_{i})=\Delta_{2}$ for all $u_{i}\in J$. If there are vertices $v$ and $v'$ not adjacent in $G_{1}$, then we may assign the vertices of $V(G^{+}_{1})-\{v'\}$ to those in $J$ such that $v$ is mapped to $u_{1}$ and assign $v'$ to some vertex in $N_{G_{2}}(u_{1})$ to construct a packing of $G_{1}$ with $G_{2}$. Thus, $G_{1}$ must be a complete graph. Since $G_{1}$ is complete and $\Delta_{2}\leq \Delta_{1}$, we know that $G_{2}\neq \Delta_{1}K^{\Delta_{2}+1}$ since $G_{1}$ and $G_{2}$ do not satisfy \ref{main:complete}. Therefore, there is a $u_{i}\in J$ with two distinct vertices $v$ and $v'$ in $N_{G_{2}}(u_{i})$ that are not adjacent. However, this implies $J-\{u_{i}\}+\{v,v'\}$ is an independent set with $\Delta_{1}+1$ vertices. Thus, $\Delta_{1}\neq \delta_{+}$.
 
Since $\Delta_{1}>\delta_{+}$, we know that $g(G_{1},G_{2})=\Delta_{2}-1$. By (\ref{eq:minG1}), we know that $|X|\geq (\Delta_{1}-1-|J|)(\Delta_{2}+1)+2$. If $\Delta_{1}-1\geq |J|$, then  $\alpha(G_{2}[X])\geq \frac{|X|}{\Delta_{2}}>\Delta_{1}-1-|J|$, and therefore, $\alpha(G_{2})=\Delta_{1}$. If $|J|=\Delta_{1}$, then $\alpha(G_{2})\geq |J|\geq \Delta_{1}$. Returning to $S$, we choose a $y$ in $G_{1}^{+}$ with minimum degree and assign the vertices of $G_{1}^{+}-\{y\}$ to vertices in $S$. Let $Z\subseteq S$ be those vertices identified with the vertices in $N_{G_{+}}(y)$. If there is some vertex $u\in V-S$ not adjacent in $G_{2}$ to a vertex in $Z$, then we can identify $y$ with $u$ to achieve a packing of $G_{1}$ and $G_{2}$. We now show such a $u$ exists. Since $y$ has degree $\delta_{+}$ in $G_{1}$ and the vertices in $Z$ are adjacent in $G_{2}$ to at most $\delta_{+}\Delta_{2}$ vertices in $V-S$, we know that \begin{equation} \label{eq:minOrder1} |V-(S\cup N_{G_{2}}(Z))|\geq n-(\Delta_{1}+\delta_{+}\Delta_{2}).
\end{equation}
Recognizing that $\min\{\Delta_{1},\Delta_{2}\}\leq \Delta_{2}$, we learn from (\ref{eq:main}) that \[n\geq \Delta_{1}\Delta_{2}+\Delta_{1}+\Delta_{2}+1-1-\Delta_{2}-(\Delta_{2}-1)=\Delta_{2}(\Delta_{1}-1)+\Delta_{1}+1.\]
When we incorporate this last inequality with (\ref{eq:minOrder1}) we see that \[|V-(S\cup N_{G_{2}}(Z))|\geq \Delta_{2}(\Delta_{1}-1)+\Delta_{1}+1-(\Delta_{1}+\delta_{+}\Delta_{2})
=\Delta_{2}(\Delta_{1}-\delta_{+}-1)+1.\]
Since $\Delta_{1}>\delta_{+}$, we see that our required $u$ exists, and we have completed the proof of Claim~\ref{cl:minOrder}.
\end{proof}
\end{claim}

\subsubsection{$G_{1}$ and $G_{2}$ $\delta_{+}$-nearly-pack}
For $x\in V$ and $H\in \mathcal{R}(\pi(G_{1}))$, if $x$ is in every $(H,G_{2})$-bad pair we say $H+G_{2}$ is a near-packing with respect to $x$. We say $G_{1}$ and $G_{2}$ $k$-nearly-pack if there exists an $H\in \mathcal{R}(\pi(G_{1}))$ and an $x\in V$ with $deg_{H}(x)=k$ such that $H+G_{2}$ is a near-packing with respect to $x$.

In this section, we show there are realizations of $G_{1}$, with useful properties, that $\delta_{+}$-nearly pack with $G_{2}$. We then finish the section with some helpful claims.

\begin{claim}\label{claim:exists}$G_{1}$ and $G_{2}$ $\delta_{+}$-nearly-pack.

\begin{proof}Let $y$ be a vertex with minimum positive degree in $G_{1}$.  Havel and Hakimi \cite{Hakimi1962,Havel1955} showed that there is a realization $F$ of $\pi(G_{1})$ such that $y$ is adjacent in $F$ to the $\delta_{+}$ vertices with highest degree in $F$. Let $H$ be the graph formed by removing in $F$ every edge incident with $y$. Clearly, $H$ has fewer vertices of positive degree than $G_{1}$, and therefore, we can say $\pi(H)$ packs with $G_{2}$ once we verify that $H$ and $G_{2}$ satisfy (\ref{eq:main}) and do not satisfy \ref{main:fivecycle}--\ref{main:matchingCompleteGraph}. 

Observe that $\Delta(G_{1})=\Delta(H)+i$ for some $i\in\{0,1\}$, and therefore, \[\min\{\Delta(G_{1}),\Delta(G_{2})\}\leq \min\{\Delta(H),\Delta(G_{2})\}+i.\]
Combining this with (\ref{eq:main}) we see that
\begin{align}
    n&\ge(\Delta(H)+i+1)(\Delta(G_{2})+1)-\min\{\Delta(H)+i,\Delta(G_{2})\}-1-g(G_{1},G_{2})\notag \\
    &\geq  (\Delta(H)+1)(\Delta(G_{2})+1)-\min\{\Delta(H),\Delta(G_{2})\}-1+i\Delta(G_{2})-g(G_{1},G_{2}).\label{eq:induction}
\end{align}
If we show that $i\Delta(G_{2})-g(G_{1},G_{2})\geq -g(H,G_{2})$, then $H$ and $G_{2}$ satisfy (\ref{eq:main}). If $H^{+}$ is not regular, then $g(H,G_{2})\geq g(G_{1},G_{2})$ and $H$ does not satisfy \ref{main:fivecycle}--\ref{main:matchingCompleteGraph}. If $H^{+}$ is regular, then $g(H,G_{2})=0$ and $G_{1}$ has at most $\delta_{+}$ vertices with maximum degree. This implies $i=1$, and therefore, $i\Delta(G_{2})-g(G_{1},G_{2})\geq 1>0 -g(H,G_{2})$. Consequently, $\Delta(H)$ and $\Delta(G_{2})$ do not obtain equality in (\ref{eq:main}) like the pairs of graphs in \ref{main:fivecycle}--\ref{main:matchingCompleteGraph}. Thus, there is a realization $H'$ of $\pi(H)$ that packs with $G_{2}$. Once we add back the edges between $y$ and the $deg_{G_{1}}(y)$ vertices in $H'$ with smaller degree than in $G_{1}$, we construct a realization $G'_{1}$ of $\pi(G_{1})$ such that  $G'_{1}+ G_{2}$ is a near-packing with respect to $y$.
\end{proof}
\end{claim}

For the rest of the proof, we let $m = \min\{\Delta(G_{1}),\Delta(G_{2})\}+1$. For an $H\in \mathcal{R}(\pi(G_{1}))$, we let $L(H)$ be the set of vertices in $H$ with degree zero. When it is clear which $H$ we are working with, we will simply use $L$ and $|L|$ for $L(H)$ and $|L(H)|$, respectively.

We define $\mathcal{\hat{H}}$ such that $(H,y)\in \mathcal{\hat{H}}$ if and only if $H\in \mathcal{R}(\pi(G_{1}))$, $y\in V$ with $deg_{G_{1}}(y)=\delta_{+}$, and $H+ G_{2}$ is a near-packing with respect to $y$. Note that Claim~\ref{claim:exists} implies $\hat{\mathcal{H}}$ is not empty.

\begin{claim}\label{cl:MainObservation}For $(H,y)\in \mathcal{\hat{H}}$, if $l\in L(H)$, then $l$ is adjacent in $G_{2}$ to a vertex in $N_{H}(y)$.
\begin{proof}If there is an $l\in L(H)$ not adjacent in $G_{2}$ to a vertex in $N_{H}(y)$, then we simply interchange $y$ with $l$ in $H$ to construct a packing of $\pi(G_{1})$ with $G_{2}$. 
\end{proof}
\end{claim}

For $(H,y)\in \mathcal{\hat{H}}$, we define the partition $N_{H}(y)=S_{2}(H,y)\cup S_{1}(H,y)\cup S_{0}(H,y)$ such that every vertex in $S_{2}(H,y)$ is adjacent in $G_{2}$ to at most $\Delta_{2}-2$ vertices in $L(H)\cup \{y\}$ and for $i\in \{0,1\}$, every vertex in $S_{i}(H,y)$ is adjacent in $G_{2}$ to exactly $\Delta_{2}-i$ vertices in $L(H)\cup \{y\}$.

\begin{remark}When $L(G_{1})\neq \emptyset$ and given an $(H,y)\in \hat{\mathcal{H}}$, the following Claim~\ref{cl:modify} allows us to modify $H$ in a useful way. The claim is critical to the proof of Theorem~\ref{theorem:BEC-half}. However, its proof is long and involves a different strategy than the main proof. It may benefit the reader to see how the claim is used in Claim~\ref{cl:desiredNearPacking} before tackling the proof of Claim~\ref{cl:modify}.
\end{remark}

\begin{claim}\label{cl:modify}For $(H,y)\in \hat{\mathcal{H}}$, if $L(G_{1})\neq \emptyset$ and there is a smallest $t\in \{0,1\}$ such that there is an $l\in L(H)\cup \{y\}$ with \begin{equation}\label{eq:modify}
    N_{G_{2}}(l)\cap N_{H}(y)\subseteq \bigcup_{0\leq j\leq t}S_{j}(H,y),
\end{equation}
then for any $s\in N_{G_{2}}(l)\cap S_{t'}(H,y)$ with $t'\leq t$, there is a set of edges $\{u_{1}w_{1},\ldots,u_{q}w_{q}\}$ in $H$ with the following properties:
\begin{enumerate}[label=\Roman*.,ref=(\Roman*)]
    \item\label{claim2con1} $u_{1}=y$, $w_{1}=s$, and $w_{i}u_{i+1}\notin E(H)\cup E(G_{2})$  for $i\leq q-1$.
    \item\label{claim2con2}$w_{q}\in V-L(H)-N_{H}[y]$.
    \item\label{claim2con3}If $w_{q}$ is adjacent in $G_{2}$ to $l$, then $w_{q}$ is adjacent in $G_{2}$ to at most $\Delta_{2}-(t'+1)$ vertices in $L(H)\cup \{y\}$.
\end{enumerate}
\begin{proof}
See Section~\ref{sec:modify}.
\end{proof}
\end{claim}

With Claim~\ref{cl:modify} we can show the existence of a desirable  $(H,y)\in \hat{\mathcal{H}}$ when $L(G_{1})\neq \emptyset$. 

\begin{claim}\label{cl:desiredNearPacking}If $L(G_{1})\neq \emptyset$, then there exist an $(H,y)\in \hat{\mathcal{H}}$ such that $N_{G_{2}}(l)\cap S_{2}(H,y)\neq \emptyset$ for every $l\in L(H)\cup \{y\}$. Consequently, $\Delta_{2}\geq 3$.
\begin{proof}Choose an arbitrary $(H,y)\in\hat{\mathcal{H}}$. We fix $y$, and define the set $X$ such that $((H',y),l)\in X$ if and only if $(H',y)\in \hat{\mathcal{H}}$ and $l\in L(H')\cup \{y\}$. We choose an $((H_{1},y),l)\in X$ such that  
\begin{enumerate}[label=\roman*.,ref=(\roman*)]
\item \label{choice1:desiredNearPacking} we maximize  $|S_{2}(H_{1},y)|$,
\item \label{choice2:desiredNearPacking} subject to \ref{choice1:desiredNearPacking}, we minimize $|N_{G_{2}}(l)\cap S_{2}(H_{1},y)|$,
\item \label{choice3:desiredNearPacking} subject to \ref{choice2:desiredNearPacking}, we minimize  $|N_{G_{2}}(l)\cap S_{1}(H_{1},y)|$,
\item \label{choice4:desiredNearPacking} subject to \ref{choice3:desiredNearPacking}, we minimize $|S_{0}(H_{1},y)|$, and
\item \label{choice5:desiredNearPacking} subject to \ref{choice4:desiredNearPacking} we minimize  $|N_{G_{2}}(l)\cap S_{0}(H_{1},y)|$.
\end{enumerate}

Suppose $N_{G_{2}}(l)\cap S_{2}(H_{1},y)=\emptyset$. By Claim~\ref{cl:MainObservation}, every vertex in $L(H_{1})$ is adjacent in $G_{2}$ to a vertex in $N_{H}(y)$. Thus, there is a largest $t\in \{0,1\}$ such that $N_{G_{2}}(l)\cap S_{t}(H_{1},y)\neq \emptyset$. By \ref{choice2:desiredNearPacking} and \ref{choice3:desiredNearPacking}, $\bigcup_{j\geq t}(N_{G_{2}}(l')\cap S_{j}(H,y))\neq \emptyset$ for all $l'\in L(H_{1})$. Therefore, $H$, $y$, $l$, and $t$ satisfy the condition of Claim~\ref{cl:modify}. Let $s\in N_{G_{2}}(l)\cap S_{t}(H_{1},y)$.  By Claim~\ref{cl:modify}, there is a set of $q$ edges
$\{u_{1}w_{1},\ldots,u_{q}w_{q}\}$ of $H_{1}$ with $u_{1}=y$, $w_{1}=s$, $w_{q}\in V-L(H_{1})-N_{H_{1}}[y]$, and $w_{j}u_{j+1}\notin E(H_{1})\cup E(G_{2})$ for $j\leq q-1$ such that $w_{q}$ satisfies \ref{claim2con3}. If we exchange in $H_{1}$ the edges $\{u_{1}w_{1},\ldots,u_{q}w_{q}\}$ with the non-edges $\{w_{1}u_{2},\ldots,w_{q}u_{1}\}$ we create a realization $H_{2}$ of $\mathcal{R}(\pi(G_{1}))$. Observe that $N_{H_{2}}(y)=N_{H_{1}}(y)-\{s\}+\{w_{q}\}$, $L(H_{1})=L(H_{2})$, and $S_{2}(H_{1},y)\subseteq S_{2}(H_{2},y)$. Suppose $l$ is not adjacent in $G_{2}$ to $w_{q}$. If $t=1$, then $((H_{2},y),l)$ contradicts \ref{choice3:desiredNearPacking}. When $t=0$, then $|S_{0}(H_{2},y)|\leq |S_{0}(H_{1},y)|$. Thus, $((H_{2},y),l)$ contradicts \ref{choice4:desiredNearPacking} or \ref{choice5:desiredNearPacking}. Now suppose $l$ is adjacent in $G_{2}$ to $w_{q}$. By \ref{claim2con3}, $w_{q}$ is adjacent in $G_{2}$ to at most $\Delta_{2}-(t+1)$ vertices in $L(H_{2})\cup \{y\}$. However, $((H_{2},y),l)$ contradicts \ref{choice1:desiredNearPacking} when $t=1$. When $t=0$, we have two cases to consider. If there is an $l'\in (L(H_{2})\cup \{y\})-\{l\}$ such that $N_{G_{2}}(l')\cap S\subseteq S_{0}(H_{2},y)$, then $((H_{2},y),l')$ contradicts \ref{choice4:desiredNearPacking} since $|S_{0}(H_{2},y)|<|S_{0}(H_{1},y)|$. When no such $l'$ exists, $((H_{2},y),l)$ contradicts \ref{choice4:desiredNearPacking} since $|S_{0}(H_{2},y)|<|S_{0}(H_{1},y)|$ and $|N_{G_{2}}(l)\cap S_{1}(H_{2},y)|=1$. Therefore, $N_{G_{2}}(l)\cap S_{2}(H_{1},y)\neq \emptyset$, and thus, $\Delta_{2}\geq 3$. Since $|N_{G_{2}}(l)\cap S_{2}(H_{1},y)|$ is minimized by \ref{choice2:desiredNearPacking}, $((H,y),l)$ proves Claim~\ref{cl:desiredNearPacking}.
\end{proof}
\end{claim}

We define $\mathcal{H}\subseteq \hat{\mathcal{H}}$ such that $(H,y)\in \mathcal{H}$ if and only if $b(H,G_{2})\leq b(H',G_{2})$ for every $(H',y')\in \hat{\mathcal{H}}$, and we define $\acute{\mathcal{H}}\subseteq \hat{\mathcal{H}}$ such that $(H,y)\in \Acute{\mathcal{H}}$ if and only if $L(G_{1})\neq \emptyset$ and $N_{G_{2}}(l)\cap S_{2}(H,y)\neq \emptyset$ for every $l\in L(H)\cup \{y\}$. By Claim~\ref{cl:desiredNearPacking}, $\acute{\mathcal{H}}$ is not empty when $L(G_{1})\neq \emptyset$. Trivially, if $|L(G_{1})|+|\{y\}|\leq \Delta_{2}-2$, then $\Acute{\mathcal{H}}=\hat{\mathcal{H}}$. On the other hand, $\Acute{\mathcal{H}}\cap \mathcal{H}$ could be empty. Since we want to take advantage of the pairs in $\Acute{\mathcal{H}}$, we end this section by giving an adequate upper bound on $|L(G_{1})|$ when  $\Acute{\mathcal{H}}\cap \mathcal{H}$ is empty. 

\begin{claim}\label{cl:L(H)bound}If $\Acute{\mathcal{H}}\cap \mathcal{H}=\emptyset$, then $|L(G_{1})|< (\delta_{+}-1)(\Delta_{2}-2)-\Big\lfloor\frac{\Delta_{2}-2}{\delta_{+}}\Big\rfloor$.
\begin{proof}The claim is trivial if $L(G_{1})=\emptyset$; suppose $L(G_{1})\neq \emptyset$. Since $\Acute{\mathcal{H}}\cap \mathcal{H}=\emptyset$, $|L(G_{1})\cup \{y\}|\geq \Delta_{2}-1$. By Claim~\ref{cl:desiredNearPacking}, there is an $(H,y)\in \Acute{\mathcal{H}}$ and $\Delta_{2}\geq 3$. Suppose  $|L(H)|\geq (\delta_{+}-1)(\Delta_{2}-2)-\Big\lfloor\frac{\Delta_{2}-2}{\delta_{+}}\Big\rfloor$. Let $S'\subseteq S_{2}(H,y)$ be a smallest set such that every vertex in $L(H)\cup \{y\}$ is adjacent in $G_{2}$ to a vertex in $S'$. Since $|L(H)\cup \{y\}|\geq \Delta_{2}-1$, $|S'|\geq 2$, and therefore, $\delta_{+}\geq 2$. We see that \[(\delta_{+}-1)(\Delta_{2}-2)-\Big\lfloor\frac{\Delta_{2}-2}{\delta_{+}}\Big\rfloor\leq |L(H)|\leq e_{G_{2}}(S', L(H))\leq |S'|(\Delta_{2}-2)-|\{y\}|.\] Rearranging we deduce that \[(\delta_{+}-1-|S'|)(\Delta_{2}-2)+1\leq \Big\lfloor\frac{\Delta_{2}-2}{\delta_{+}}\Big\rfloor.\] This can only be true if $|S'|\geq \delta_{+}-1$ since $\delta_{+}\geq 2$ and $\Delta_{2}\geq 3$.  If there is an $l\in L(H)\cup \{y\}$ that is adjacent to exactly one vertex in $N_{H}(y)$, then we interchange $y$ with $l$ in $H$ to construct another realization $H'$ of $\pi(G_{1})$. Note that $(H',l)\in \Acute{\mathcal{H}}$ since $L(H')\cup \{l\}=L(H)-\{y\}+\{l\}$ and $S_{2}(H,y)=S_{2}(H',y)$. However, $H'+ G_{2}$ has exactly one bad pair, and therefore, we have the contradiction $(H',l)\in \mathcal{H}$. Thus, no vertex in $L(H)\cup \{y\}$ is adjacent in $G_{2}$ to exactly one vertex in $N_{H}(y)$. Thus, by our selection of $S'$, $|S'|=\delta_{+}-1$, and therefore, $\delta_{+}\geq 3$. Moreover, it must be the case that \[2|L(H)\cup \{y\}|\leq e_{G_{2}}(N_{H}(y), L(H)\cup \{y\})\leq \delta_{+}\Delta_{2}-2|S'|.\] This implies  
\[2((\delta_{+}-1)(\Delta_{2}-2)-\Big\lfloor\frac{\Delta_{2}-2}{\delta_{+}}\Big\rfloor+1)\leq \delta_{+}(\Delta_{2}-2)+2.\] Rearranging and simplifying we see that \[(\delta_{+}-2)(\Delta_{2}-2)\leq \frac{\Delta_{2}-2}{\delta_{+}}.\] However, this is a contradiction since $\delta_{+}\geq 3$ and $\Delta_{2}\geq 3$. Thus, $|L(H)|< (\delta_{+}-1)(\Delta_{2}-2)-\Big\lfloor\frac{\Delta_{2}-2}{\delta_{+}}\Big\rfloor$.\qedhere
\end{proof}
\end{claim}

The rest of the proof is broken up into the cases $\delta_{+}=1$ and $\delta_{+}\geq 2$. For both cases, we begin by choosing a pair in $\mathcal{H}$ and a related bad pair, but in different ways. However, once those choices are made, the setup for their respective arguments are similar. Therefore, in the next section we define the common sets we need and prove some supporting claims about them.
\subsubsection{\label{sec:setup} Preliminary Setup and Supporting Claims.}

 \begin{figure}[t]
    \centering
    \begin{tikzpicture}

\coordinate (x) at (2.75,4.5);
\coordinate (y) at (4.25,4.5);
\coordinate (u) at (2.5,2.75);
\coordinate (v) at (3.5,2.5);
\coordinate (z) at (4.75,1.5);

\draw[] (2,0) rectangle (5,.8);
\draw[] (2,1) rectangle (5,3.25);
\draw[] (0,1) rectangle (1.75,3.25);
\draw[] (5.25,1) rectangle (7,3.25);
\draw[] (2.25,1) rectangle (5,1.80);  
\draw[color=red, thick] (x) -- (.75,2.75);
\draw[color=blue, very thick] (x) -- (.75,2);
\draw[color=red, thick] (y) -- (6.5,2.75);
\draw[color=blue, very  thick] (y) -- (6.5,2);
\draw[color=blue, very  thick] (x) -- (u);
\draw[color=red, thick] (y) -- (u);
\draw[color=blue, very thick] (x) -- (v);
\draw[color=blue, very thick] (y) -- (v);
\draw[color=red, thick] (x) -- (z);
\draw[color=blue, very thick] (y) -- (z);
\draw[color=red, thick] (4.45,.15) -- (3.95,1.45);
\draw[color=red, thick] (4.45,.15) -- (4.45,1.45);
\draw[color=red, thick] (4.45,.15) -- (4.95,1.45);
\draw[color=red, thick] (3.75,.15) -- (4.45,1.45);
\draw[color=red, thick] (3.75,.15) -- (3.25,1.15);
\draw [red, thick]  (x) to[out=-20,in=200] (y);
\draw [blue, very thick] (x) to[out=20,in=160] (y);

  \draw[fill=black] (x) circle (.1cm);
  \draw[fill=black] (y) circle (.1cm);

\node[] at (2.5,4.75)  {$x$};
\node[] at (4.5,4.75) {$y$};
\node[] at (2.85,.35) {$I_{G}(x,y)$};
\node[] at (2.85,2.15) {$W_{G}(x,y)$};
\node[] at (.9,1.5) {$X_{G}-Y_{G}$};
\node[] at (6.15,1.5) {$Y_{G}-X_{G}$};
\node[] at (3.1,1.45) {$\hat{W}_{G}(x,y)$};
\end{tikzpicture}
    \caption{This is a typical setup for a given $(H,y)\in \mathcal{H}$ and a $(H,G_{2})$-bad pair $(x,y)$. Edges of $H$ are colored red, and edges of $G_{2}$ are a little thicker and colored blue.}
    \label{fig:TypicalSetup}
\end{figure}
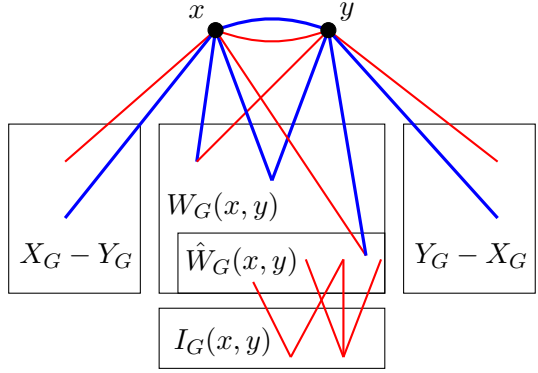
 
Given an $(H,y)\in \mathcal{H}$ and a $(H,G_{2})$-bad pair $(x,y)$. We denote $H\cup G_{2}$ by $G$. Let $W_{G}[x,y] = N_{G}[x]\cap N_{G}[y]$, $W_{G}(x,y)=W_{G}[x,y]-\{x,y\}$, $I_{G}(x,y)=V-N_{G}[\{x,y\}]$, and $\hat{W}_{G}(x,y)\subseteq W_{G}(x,y)$ be the set of vertices adjacent in $H$ to vertices in $I_{G}(x,y)$. We let $X_{G}=N_{G}[x]$ and $Y_{G}=N_{G}[y]$. See Figure~\ref{fig:TypicalSetup} for a visualization of these sets. Since every vertex in $\hat{W}_{G}(x,y)$ is adjacent in $H$ to at most $\Delta_{1}$ vertices in $I_{G}(x,y)$, we define \[q_{G}(x,y)=\sum_{u\in \hat{W}_{G}(x,y)}(\Delta_{1}-|N_{H}(u)\cap I_{G}(x,y)|).\] We also need to define $\epsilon_{G}(x,y)=|W_{G}(x,y)|-|\hat{W}_{G}(x,y)|$. When it is clear which $(H,y)$ and bad pair $(x,y)$ is under consideration, we may simplify our notation by dropping the subscripts in favor of $q$, $\epsilon$, $I$, $X$, $Y$, $\hat{W}$, and $W$ for $W(x,y)$. Similarly, we will use $g$ in place of $g(H, G_{2})$.

\begin{claim}\label{cl:setup}The following is true:
\begin{itemize}
\item $N_{H}(u)\subseteq W_{G}(x,y)$ for all $u\in I_{G}(x,y)$.
\item If $u\in Y_{G}-X_{G}$, then $N_{H}(u)\subseteq Y_{G}$.
\item If $u\in N_{H}(y)\cap (Y_{G}-X_{G})$ and $v\in X_{G}-Y_{G}$, then $N_{H}(v)-\{x\}\subseteq W_{G}(x,y)$. 
\item If $u\in N_{H}(y)\cap (Y_{G}-X_{G})$ and $v\in N_{H}(x)\cap (X_{G}-Y_{G})$, then $N_{H}(u)-\{y\}\subseteq W_{G}(x,y)$.
\end{itemize}
\begin{proof}Suppose $u\in Y\cup I$, and let $u'\in N_{H}(u)$. If $u'\in X\cup I$, then we exchange in $H$ the edges $xy$ and $uu'$ with the non-edges $xu$ and $yu'$ of $G$ to construct an $H'\in \mathcal{R}(\pi(G_{1}))$ such that  $(H',y)\in\hat{\mathcal{H}}$ and $b(H',G_{2})<b(H,G_{2})$. However, this contradicts $(H,y)\in \hat{\mathcal{H}}$, and therefore, $u'\notin X\cup I$. A similar argument holds if $u\in X\cup I$ and $u'\in Y\cup I$. Thus, we have proved the first two listed items. Now suppose there is a $u\in N_{H}(y)\cap (Y-X)$ and a $v\in X-Y$. If there is a $v'\in N_{H}(v)\cap (Y-X)$, then we can exchange in $H$ the edges $xy$, $yu$, and $vv'$ with the non-edges $xu$, $yv'$, and $yv$ of $G$ to construct an $H'\in \mathcal{R}(\pi(G_{1}))$ such that  $(H',y)\in\hat{\mathcal{H}}$ and $b(H',G_{2})<b(H,G_{2})$. However, this contradicts $(H,y)\in\hat{\mathcal{H}}$, and therefore, $v'\in W_{G}(x,y)$. Thus, $N_{H}(v)-\{x\}\subseteq W(x,y)$. If $v\in N_{H}(x)\cap (X-Y)$ and there is a $u'\in N_{H}(u)\cap (Y-X)$, then we can exchange in $H$ the edges $xy$, $xv$, and $uu'$ with the non-edges $xu$, $xu'$, and $yv$ of $G$ to construct an $H'\in \mathcal{R}(\pi(G_{1}))$ such that  $(H',y)\in\hat{\mathcal{H}}$ and $b(H',G_{2})<b(H,G_{2})$. Thus, $N_{H}(u)-\{y\}\subseteq W(x,y)$.
\end{proof}
\end{claim}

\begin{claim}\label{cl:setup2}If $u\in  V-Y_{G}$ and $v\in  V-X_{G}$, then every vertex in $N_{H}(u)$ is adjacent in $G$ to every vertex in $N_{H}(v)$.
\begin{proof}
    If there is a $u'\in N_{H}(u)$ and a $v'\in N_{H}(v)$ that are not adjacent in $G$, then we can exchange in $H$ the edges $uu'$, $v'v$, and $xy$ with the non-edges $u'v'$, $vx$, and $yu$ of $G$ to construct an $H'\in \mathcal{R}(\pi(G_{1}))$ such that $(H',y)\in \hat{\mathcal{H}}$ and $b(H',G_{2})<b(H,G_{2})$. However, this is a contradiction since $(H,y)\in \mathcal{H}$. Thus, every vertex in $N_{H}(u)$ is adjacent in $G$ to every vertex in $N_{H}(v)$.
\end{proof}
\end{claim}

\begin{claim}\label{cl:goodneighbors}If $deg_{H}(z)\leq 2$ for some $z\in \{x,y\}$, then $deg_{H}(u)\geq 2$ for all $u\in N_{H}(z)-W[x,y]$.
\begin{proof}If there is a $u\in N_{H}(z)-W[x,y]$ such that $deg_{H}(u)=1$, then we may interchange $u$ with $z$ to construct a realization of $\pi(G_{2})$ that packs with $G_{2}$. This contradiction implies $deg_{H}(u)\geq 2$ for all $u\in N_{H}(z)-W[x,y]$.
\end{proof}
\end{claim}

\begin{claim}\label{cl:minW}There is a $u\in N_{G}[x,y]-W_{G}(x,y)$ such that $|N_{H}(u)\cap W_{G}(x,y)|\geq \delta_{+}-1$.
\begin{proof}The claim is trivial when $\delta_{+}=1$. The claim also easily follows if $N_{H}(x)\subseteq W_{G}[x,y]$, or $N_{H}(y)\subseteq W_{G}[x,y]$. When $\delta_{+}\geq 2$ and there are vertices $u\in N_{H}(y)-W[x,y]$ and $v\in N_{H}(x)-W[x,y]$, the last item in Claim~\ref{cl:setup} says $N_{H}(u)-\{y\}\subseteq W_{G}(x,y)$. Thus, $|N_{H}(u)\cap W(x,y)|\geq \delta_{+}-1$.
\end{proof}
\end{claim}

We let $k_{y}=\Delta_{2}-deg_{G_{2}}(y)$ and $k_{x}=\Delta_{2}-deg_{G_{2}}(x)$. Since $(x,y)$ is a bad pair, we know $|X|\leq deg_{H}(x)+deg_{G_{2}}(x)\leq  \Delta_{1}+\Delta_{2}-k_{x}$ and $|Y|\leq deg_{H}(y)+deg_{G_{2}}(y)=\delta_{+}+\Delta_{2}-k_{y}$.  Furthermore,
\begin{align}
|N_{G}[x,y]|&=|X\cup Y|\notag \\
&\leq |X|+|Y|-|W[x,y]|\notag \\
&\leq deg_{H}(x)+deg_{G_{2}}(x)+deg_{H}(y)+deg_{G_{2}}(y)-|W[x,y]|\notag \\
&\leq \Delta_{1}+\Delta_{2}-k_{x}+\delta_{+}+\Delta_{2}-k_{y}-|W[x,y]|\notag \\
&= \Delta_{1}+\delta_{+}+2\Delta_{2}-|W(x,y)|-2-k_{x}-k_{y}.\notag
\end{align}
As $n\geq \Delta_{1}\Delta_{2}+\Delta_{1}+\Delta_{2}+1-m-g$, we can deduce that
\begin{align}
|I_{G}(x,y)|&=|V-N_{G}[x,y]|\notag\\
&\geq\Delta_{2}(\Delta_{1}-1)+\Delta_{1}-deg_{H}(x)+3-m-g-\delta_{+}+|W(x,y)|+k_{x}+k_{y}. \label{eq:I}
\end{align}
\begin{claim}\label{cl:InotEmpyt} $I_{G}(x,y)\neq \emptyset$.
\begin{proof}If $\Delta_{1}-deg_{H}(x)+|W(x,y)|+k_{x}(H)+k_{y}(H)\geq 
\delta_{+}$, then by (\ref{eq:I}), \[|I|\geq \Delta_{2}(\Delta_{1}-1)+3-m-g\geq \Delta_{2}(\Delta_{1}-2)+4-m\geq 1.\] Thus, we may assume $\delta_{+}-1\geq \Delta_{1}-deg_{H}(x)+|W(x,y)|+k_{x}+k_{y}$. However, Claim~\ref{cl:minW} implies that $\delta_{+}-1= |W(x,y)|$, $\Delta_{1}=deg_{H}(x)$, and $k_{x}=k_{y}=0$. If $\Delta_{1}\geq 3$, then from (\ref{eq:I}) we see that \[|I|\geq \Delta_{2}(\Delta_{1}-1)+2-m-g\geq \Delta_{2}(\Delta_{1}-2)+3-m\geq 2.\] So, we assume $\Delta_{1}\leq 2$. If $\delta_{+}=2$, then $\Delta_{1}=2$, and therefore, $g=0$ and $|W(x,y)|= 1$. Therefore, by (\ref{eq:I}), $|I|\geq 2(\Delta_{2}+1)-m\geq 1$. Thus, $\delta_{+}=1$. In this case, $|W(x,y)|=0$.
By Claim~\ref{scl:basicObservation}, $L(H)\subseteq N_{G}(x)$. 

Suppose $\Delta_{1}=2$. Thus, there is a $v\in N_{H}(x)\cap (X-Y)$. If $\Delta_{2}=1$, then (\ref{eq:I}) says $|I|\geq 1$. If $\Delta_{2}\geq 2$, then there must be a $u\in N_{G_{2}}(y)\cap (Y-X)$ since $k_{y}=0$. This implies there is a $u'\in N_{H}(u)\cap (Y-X)$. However, we have a contradiction by the last item in Claim~\ref{cl:setup}. 

We are left with the case $\Delta_{1}=1$. In this situation, $(H,x)\in \hat{\mathcal{H}}$ and therefore, by Claim~\ref{cl:MainObservation}, $L(H)=\emptyset$. Thus, $H$ is a perfect matching. By Claim~\ref{cl:setup}, there are no edges of $H$ between $X-Y$ and $Y-X$. Thus, $H[X-\{x,y\}]$ and $H[Y-\{x,y\}]$ are perfect matchings. This implies $\Delta_{2}$ is odd since $|X-\{x,y\}|=|Y-\{x,y\}|=\Delta_{2}-1$.  If $\Delta_{2}=1$, then $x$ and $y$ are the sole vertices in $H$ and together violate \ref{main:matchCompleteBipart}. If $\Delta_{2}\geq 2$, then since both $X-\{x,y\}$ and $Y-\{x,y\}$ contain edges of $H$, we see by Claim~\ref{cl:setup}, that every vertex in $Y-X$ is adjacent in $G_{2}$ to every vertex in $X-Y$. Thus, $(X-Y)\cup \{y\}$ and $(Y-X)\cup \{x\}$ are independent sets in $G_{2}$. Since $|(X-Y)\cup \{y\}|=|(Y-X)\cup \{x\}|=\Delta_{2}$, $G_{1}$ and $G_{2}$ satisfy \ref{main:matchCompleteBipart}. However, this is a contradiction, and thus, $I$ must not be empty.
\end{proof}
\end{claim}

\begin{claim}\label{cl:WIcomplete}$\hat{W}$ is complete in $G$.
\begin{proof}Since every vertex in $\hat{W}$ is adjacent in $H$ to a vertex in $I$, this claim easily follows from Claim~\ref{cl:setup2}.
\end{proof}
\end{claim}

We will often need the following inequality.
\begin{equation}\label{eq:bounds}
\delta_{+}|I-L|\leq e_{H}(I,\hat{W})= \Delta_{1}|\hat{W}|-q.
\end{equation}

\begin{claim}\label{cl:epAnddelta}$\epsilon+q\geq \delta_{+}-1$.
\begin{proof}By (\ref{cl:minW}), there is a $u\in N_{G}[x,y]-W_{G}(x,y)$ with at least $\delta_{+}-1$ neighbors in $W(x,y)$. Thus, $\epsilon \geq e_{H}(u, W(x,y)-\hat{W})$ and $q\geq e_{H}(u, \hat{W})$. Since \[\epsilon+q\geq e_{H}(u, W(x,y)-\hat{W})+e_{H}(u, \hat{W})\geq \delta_{+}-1,\] we have proved the claim.
\end{proof}
\end{claim}

\subsubsection{\label{sec:d=1}$\delta_{+}=1$}

For this section, we choose an $(H,y)\in \mathcal{H}$ and an $(H,G_{2})$-bad pair $(x,y)$ such that  
\begin{enumerate}[label=(C\arabic*),ref=(C\arabic*)]
    \item \label{choice:1}$|N_{H}(x)\cap W_{H+ G_{2}}(x,y)|$ is minimized,
    \item \label{choice:2} subject to \ref{choice:1}, $e_{H}(N_{H}(x)-W_{H+ G_{2}}(x,y),W_{H+ G_{2}}(x,y))$ is minimized, and
    \item \label{choice:3} subject to \ref{choice:2}, $deg_{H}(x)$ is maximized.
\end{enumerate}
Let $G=H+ G_{2}$. By Claim~\ref{cl:MainObservation}, we know that $L(H)\subseteq X$, and by Claim~\ref{cl:InotEmpyt}, we know that $I_{G}(x,y)$ is not empty. Therefore, $W_{G}(x,y)$ is not empty. Since $deg_{H}(y)=1$, $W_{G}(x,y)\subseteq N_{G_{2}}(y)$, and therefore, $\Delta_{2}\geq 2$ and $|W_{G}(x,y)|\leq \Delta_{2}-1$. 

\begin{claim}\label{cl:LisAdjacenttohat{W}}Every vertex in $L(H)$ is adjacent in $G_{2}$ to every vertex in $\hat{W}_{G}(x,y)$.
\begin{proof}
    Suppose there is a $u\in \hat{W}$ that is not adjacent in $G$ to some $l\in L$. Given a $w\in N_{H}(u)\cap I$ we interchange $y$ with $l$ in $H$ to construct a realization $H'$ of $\pi(G_{1})$. We then exchange in $H'$ the edges $xl$ and $wu$ with the non-edges $ul$ and $xw$ of $H'+ G_{2}$ to construct a realization that packs with $G_{2}$. Thus, every vertex in $\hat{W}$ is adjacent in $G_{2}$ to every vertex in $L$.
\end{proof}
\end{claim}

\begin{claim}\label{cl:d1neq1}$\Delta_{1}>1$
\begin{proof}Suppose $\Delta_{1}=1$. From (\ref{eq:I}), we see that $|I|=|\hat{W}|+\epsilon+k_{y}+k_{x}$. On the other hand, (\ref{eq:bounds}) implies $|I|\leq |\hat{W}|-q$. Combining these bounds implies $\epsilon+k_{y}+k_{x}+q=0$. Thus, $Y-X$ and $X-Y$ are empty, and therefore, $L$ is empty. Thus, $G_{1}$ is a complete matching and $G_{2}[Y]$ induces a $K^{\Delta_{2}+1}$. However, this contradicts our assumption that $G_{1}$ and $G_{2}$ do not satisfy \ref{main:matchingCompleteGraph}. Therefore, $\Delta_{1}\geq 2$.
\end{proof}
\end{claim}

Since $g(H,G_{2})=\Delta_{2}-1$, we can simplify the right hand side of (\ref{eq:I}) to 
\begin{equation}\label{eq:Idel1}
    |I|\geq \Delta_{2}(\Delta_{1}-2)+\Delta_{1}-deg_{H}(x)+3-m+|\hat{W}|+\epsilon+k_{y}+k_{x}.
\end{equation}

We let $I_{1}$ be the set of vertices in $I$ with degree one in $H$. This implies \[2|I|-|I_{1}|\leq e_{H}(I,\hat{W})=\Delta_{1}|\hat{W}|-q.\] After solving for $|I_{1}|$, we deduce that \begin{equation}\label{eq:I1}
    |I_{1}|\geq (\Delta_{1}-2)(2(\Delta_{2}-1)-|\hat{W}|)+2(\Delta_{1}-deg_{H}(x))+2(\epsilon+k_{y}+k_{x})+q.
\end{equation}

Let $Z\subseteq \hat{W}$ be the set of vertices adjacent in $H$ to a vertex in $I_{1}$. We define \[q'=\sum_{u\in Z}(\Delta_{1}-|N_{H}(u)\cap I_{1}|).\] Thus, $e_{H}(I_{1}, Z)= \Delta_{1}|Z|-q'$. Note that $|I_{1}|=e_{H}(I_{1}, Z)$. 

\begin{claim}\label{cl:q'Bound} $q+q'\leq 2|Z|$. Moreover, if $q+q'= 2|Z|$, then $deg_{H}(x)=\Delta_{1}$, $\Delta_{1}\geq 3$, $\epsilon+k_{y}+k_{x}=0$, and $\Delta_{2}-1=|Z|$.
    \begin{proof}Suppose $q+q'\geq 2|Z|+c$ for some non-negative integer $c$. Using $\Delta_{1}|Z|-q'= |I_{1}|$ and (\ref{eq:I1}) we deduce that
  \begin{align}0&= |I_{1}|-\Delta_{1}|Z|+q'\notag\\
  &\geq (\Delta_{1}-2)(2(\Delta_{2}-1)-|\hat{W}|)+2(\Delta_{1}-deg_{H}(x))+2(\epsilon+k_{y}+k_{x})-\Delta_{1}|Z|+q+q'\notag\\
  &= (\Delta_{1}-2)(2(\Delta_{2}-1)-|\hat{W}|)+2(\Delta_{1}-deg_{H}(x))+2(\epsilon+k_{y}+k_{x})-\Delta_{1}|Z|+2|Z|+c\notag\\
  &= (\Delta_{1}-2)(2(\Delta_{2}-1)-|\hat{W}|-|Z|)+2(\Delta_{1}-deg_{H}(x))+2(\epsilon+k_{y}+k_{x})+c.\notag
  \end{align}
This last inequality is only true if $deg_{H}(x)=\Delta_{1}$, $\epsilon+k_{y}+k_{x}=0$, $c=0$, and either $\Delta_{1}=2$ or $\hat{W}=Z$ and $\Delta_{2}-1=|Z|$. Since $c=0$, we may conclude that $q+q'= 2|Z|$. Suppose $\Delta_{1}=2$. Since $deg_{H}(x)=2$, there is a $u\in N_{H}(x)-\{y\}$. If $u\in X-Y$, then we may interchange $x$ with $u$ in $H$ to construct a realization of $\pi(G_{1})$ that packs with $G_{2}$. Since $\epsilon=0$, $W_{G}(x,y)=\hat{W}$. If $u\in \hat{W}$, then $q\geq 1$.  Let $w\in I_{1}$ and $v\in N_{H}(w)$. If $v=u$, then we may interchange $x$ with $w$ in $H$ to construct a realization of $\pi(G_{1})$ that packs with $G_{2}$. If $u\neq v$, then we exchange the edges $xy$, $xu$, and $wv$ of $H$ with the non-edges $uy$, $xv$ and $xw$ of $G$ to construct a realization $H'$ of $\pi(G_{1})$. Next, we interchange $x$ with $w$ in $H'$, and finally, we interchange $x$ with $y$ in the resulting realization to construct a realization of $\pi(G_{1})$ that packs with $G_{2}$. Therefore, $\Delta_{1}\geq 3$.
\end{proof}
\end{claim}

\begin{claim}\label{cl:NoGoodI1}If $|N_{H}(x)\cap W_{G}(x,y)|\geq 1$, $e_{H}(N_{H}(x)-W_{G}(x,y),W_{G}(x,y))\geq 1$, or $deg_{H}(x)=1$, then there does not exist a $z\in Z$ with $deg_{H}(z)\geq 2$ such that $N_{H}(z)\cap (Y\cup N_{H}(x))=\emptyset$ and $|N_{H}(z)\cap I_{1}|\geq deg_{H}(z)-1$.
\begin{proof}Suppose such a $z$ exists. We choose a $w\in N_{H}(z)$ while giving preference to one not in $I_{1}$ if it exists. We first interchange $y$ with $x$ in $H$ to construct an $H'\in \mathcal{R}(\pi(G_{1}))$. In $H'$ we exchange the edges $wz$ and $xy$ with the non-edges $wy$ and $xz$ of $G$ and to create an $H''\in \mathcal{R}(\pi(G_{1}))$ such that $(H'',y)\in \hat{H}$ and has the bad pair $(y,z)$. We now have a contradiction to \ref{choice:1}, \ref{choice:2}, and \ref{choice:3} since $deg_{H''}(z)\geq 2$, $|N_{H''}(z)\cap W_{H''+ G_{2}}(z,y)|=0$, and every vertex in $N_{H''}(z)$ has degree one in $H''$.
\end{proof}
\end{claim}

\begin{claim}\label{cl:I1stats}If $|N_{H}(x)\cap W_{G}(x,y)|\geq 1$, $e_{H}(N_{H}(x)-W_{G}(x,y),W_{G}(x,y))\geq 1$, or $deg_{H}(x)=1$, then $q+q'=2|Z|$.
\begin{proof}By Claim~\ref{cl:NoGoodI1}, every $z\in Z$ either has degree one in $H$, $|N_{H}(z)\cap N_{G}[x,y]|\geq 1$, or $|N_{H}(z)\cap (I-I_{1})|\geq 2$. Let $P\subseteq Z$ be the set of all vertices adjacent in $H$ to $\Delta_{1}$ vertices in $I$. We deduce that \[q\geq \sum_{z\in Z-P}\Delta_{1}-|N_{H}(z)\cap I|\geq |Z-P|.\] By Claim~\ref{cl:NoGoodI1}, every vertex in $P$ must be adjacent in $H$ to at least two vertices in $I-I_{1}$. Therefore, \[q'\geq \sum_{z\in Z-P}\Delta_{1}-|N_{H}(z)\cap I|+\sum_{z\in P} \Delta_{1}-|N_{H}(z)\cap I_{1}|\geq |Z-P|+2|P|.\] Combining the last two inequalities we see that $q+q'\geq 2|Z-P|+2|P|=2|Z|$. Thus, $q+q'=2|Z|$ by Claim~\ref{cl:q'Bound}.
\end{proof}
\end{claim}

By Claim~\ref{cl:I1stats} and Claim~\ref{cl:q'Bound}, if $|N_{H}(x)\cap W_{G}(x,y)|\geq 1$, $e_{H}(N_{H}(x)-W_{G}(x,y),W_{G}(x,y))\geq 1$, or $deg_{H}(x)=1$, then $deg_{H}(x)=\Delta_{1}$, $\Delta_{1}\geq 3$, $\epsilon+k_{y}+k_{x}=0$, and $\Delta_{2}-1=|Z|$. Thus, $deg_{H}(x)\geq 2$.

\begin{claim}\label{cl:exactlyOne}If $|N_{H}(x)\cap W_{G}(x,y)|\geq 1$ or $e_{H}(N_{H}(x)-W_{G}(x,y),W_{G}(x,y))\geq 1$, then every vertex in $Z$ is adjacent in $H$ to exactly one vertex in $Y\cup N_{H}(x)$ and  $\Delta_{1}-1$ vertices in $I_{1}$.
\begin{proof} Suppose there is a $z\in Z$ such that $N_{H}(z)\cap (Y\cup N_{H}(x))=\emptyset$. Previously we showed $Z=W_{G}(x,y)$ and $|Z|=\Delta_{2}-1$, and therefore, $z$ is adjacent in $G_{2}$ to every vertex in $W_{G}[x,y]$. Thus, $z$ is not adjacent in $G_{2}$ to any vertex in $N_{H}(x)\cap (X-Y)$. Let $w\in N_{H}(z)\cap I_{1}$. If there is a $u\in N_{H}(x)\cap W_{G}(x,y)$, then in $H$ we exchange the edges $wz$ and $xy$ with the non-edges $wx$ and $yz$ of $G$ to construct an $H'\in \mathcal{R}(\pi(G_{1}))$ such that $(H',y)\in \hat{H}$ and has the bad pair $(y,z)$. We now have a contradiction to \ref{choice:1} since $|N_{H'}(z)\cap W_{H'+ G_{2}}(z,y)|=0$. Thus, $N_{H}(x)\cap W_{G}(x,y)=\emptyset$. If there is a $u\in N_{H}(x)$ that is adjacent in $H$ to a vertex in $W_{G}(x,y)$, then in $H$ we exchange the edges $wz$ and $ux$ with the non-edges $wx$ and $uz$ of $G$ to construct an $H'\in \mathcal{R}(\pi(G_{1}))$ such that $(H',y)\in \hat{H}$ and has the bad pair $(x,y)$. We now have a contradiction to \ref{choice:2} since $|e_{H}(N_{H}(x)-W_{G}(x,y),W_{G}(x,y))|>|e_{H'}(N_{H'}(x)-W_{H'+ G_{2}}(x,y)),W_{H'+ G_{2}}(x,y))|$. Thus, every vertex in $Z$ is adjacent in $H$ to a vertex in $Y\cup N_{H}(x)$. This implies $q\geq\sum_{v\in Z} \big (\Delta_{1}-|N_{H}(v)\cap I|\big )\geq |Z|$. If there is a $v\in Z$ with $|N_{H}(v)\cap I_{1}|\leq \Delta_{1}-2$, then $q'\geq \sum_{v\in Z}\big (\Delta_{1}-|N_{H}(v)\cap I_{1}|\big )\geq |Z|+1$. However, this gives us a contradiction since $q+q'\leq 2|Z|$. Thus, every vertex in $Z$ is adjacent to exactly one vertex in $Y\cup N_{H}(x)$ and $\Delta_{1}-1$ vertices in $I_{1}$.
\end{proof}
\end{claim}

\begin{claim}\label{cl:uNotAdjacenttoI}If there is a $u\in N_{H}(x)\cap W_{G}(x,y)$, then $\{u\}=N_{H}(x)\cap W_{G}(x,y)$, and $N_{G_{2}}(u)=N_{G_{2}}(x)$.
\begin{proof}We showed previously that $Z=W_{G}(x,y)$ since $N_{H}(x)\cap W_{G}(x,y)\neq \emptyset$. By Claim~\ref{cl:LisAdjacenttohat{W}}, $u$ is adjacent in $G_{2}$ to every vertex in $L$. By Claim~\ref{cl:exactlyOne}, no vertex in $N_{G_{2}}(x)\cap (X-Y)$ is adjacent in $H$ to a vertex in $Z$. Therefore, Claim~\ref{cl:setup2} implies every vertex in $N_{G_{2}}(x)\cap (X-Y-L)$ is adjacent in $G_{2}$ to $u$. Moreover, since $u$ is adjacent in $H$ to $x$ and is not adjacent in $H$ to any other vertex in $Y\cup N_{H}(x)$, we deduce that $u$ is adjacent in $G_{2}$ to every vertex in $N_{G_{2}}(x)\cup (N_{H}(x)\cap (Z-\{u\}))$. Since $k_{x}=0$ implies $deg_{G_{2}}(x)=\Delta_{2}$, $u$ must be the only vertex in $N_{H}(x)\cap W_{G}(x,y)$. This completes the proof of Claim~\ref{cl:uNotAdjacenttoI}.
\end{proof}
\end{claim}

Since $deg_{H}(x)\geq 3$ when $N_{H}(x)\cap W_{G}(x,y)\neq \emptyset$, Claim~\ref{cl:uNotAdjacenttoI} implies $N_{H}(x)-Y$ is not empty. Therefore, $N_{H}(x)-Y$ is not empty since $deg_{H}(x)\geq 2$ in general. We partition $N_{H}(x)-Y=D_{1}\cup D_{2}$ where $D_{1}$ are all vertices $u\in N_{H}(x)-Y$ such that $N_{H}(u)\subseteq N_{H}[x]$. If $u\in D_{1}$ is not adjacent in $G_{2}$ to a vertex in $N_{H}(x)$, then we may interchange in $H$ $u$ with $x$ to find a realization of $\pi(G_{1})$ that packs with $G_{2}$. Therefore, every vertex in $D_{1}$ is adjacent in $G_{2}$ to a vertex in $N_{H}(x)$.

Suppose there is a $u\in D_{2}$ that is not adjacent in $G$ to a $v\in \hat{W}$. Let $w\in N_{H}(v)\cap I_{1}$, and $u'\in N_{H}(u)-D_{1}$. By Claim~\ref{cl:setup2}, we can reason that $u'\in W_{G}(x,y)$. If we exchange the edges $vw$ and $ux$ with the non-edges $uv$ and $xw$ of $G$, then we create a $H'\in \mathcal{R}(\pi(G_{1})$ such that $(H',y)\in \mathcal{H}$ and $e_{H}(N_{H}(x)-W_{G}(x,y),W_{G}(x,y))>e_{H'}(N_{H'}(x)-W_{H'+ G_{2}}(x,y),W_{H'+ G_{2}}(x,y))$. However, this contradicts \ref{choice:2}. Therefore, every vertex in $D_{2}$ must be adjacent in $G$ to every vertex in $Z$. 

We choose a $w\in I_{1}$ that is adjacent in $G_{2}$ to the fewest vertices in $N_{H}(x)$. Let $k=|N_{G_{2}}(w)\cap N_{H}(x)|$. Suppose $k\leq 1$. We choose a $v\in N_{H}(x)$ such that we give preference to one in $N_{G_{2}}(w)$ if it exists. Let $u$ be the unique vertex in $N_{H}(w)\cap Z$. If $u\notin N_{H}(x)$, then we exchange the edges $xv$, $xy$, and $wu$ with the non-edges $yv$, $xu$, and $xw$ of $G$. Then in the resulting realization we interchange $x$ with $w$ to create a realization that packs with $G_{2}$. If $u\in N_{H}(x)$, then since Claim~\ref{cl:exactlyOne} and Claim~\ref{cl:uNotAdjacenttoI} imply $u$ is not adjacent in $G$ to $v$, we exchange the edges $xv$ and $uw$ with the non-edges $xw$ and $uv$ of $G$. We then interchange in the resulting realization $x$ with $w$ to create a realization that packs with $G_{2}$. Therefore, every vertex in $I_{1}$ is adjacent in $G_{2}$ to at least two vertices in $N_{H}(x)-\{y\}$. As a result, $k\geq 2$ and $\Delta_{1}\geq 3$.

If there is some $u\in D_{1}$ not adjacent in $G_{2}$ to a vertex in $N_{H}(x)$, then we have a contradiction since we may create a packing of $\pi(G_{1})$ and $G_{2}$ by interchanging $u$ with $x$ in $H$. By Claim~\ref{cl:uNotAdjacenttoI}, we may conclude that $u$ must be adjacent in $G_{2}$ to a vertex in $D_{1}\cup D_{2}$. Therefore, every vertex in $D_{1}$ is adjacent in $G_{2}$ to a vertex in $D_{1}\cup D_{2}$.

Let $D'_{2}\subseteq D_{2}$ be the set of vertices that are adjacent in $G_{2}$ to $\Delta_{2}$ vertices in $I_{1}$. This implies every vertex in $D'_{2}$ is adjacent in $H$ to every vertex in $Z$. Thus, $q\geq |D'_{2}||Z|$.

Observe that \[ k|I_{1}|\leq e_{G_{2}}(I_{1},D_{1}\cup D_{2})\leq (deg_{H}(x)-1)(\Delta_{2}-1)+|D'_{2}|\leq (\Delta_{1}-1)(\Delta_{2}-1)+|D'_{2}|.\]
Using (\ref{eq:I1}) we deduce that 
\[ k|I_{1}|\geq k(\Delta_{1}-2)(2(\Delta_{2}-1)-|\hat{W}|)+2k(\Delta_{1}-deg_{H}(x))+2k(\epsilon+k_{y}+k_{x})+kq.\] To simplify things we let $c=2k(\Delta_{1}-deg_{H}(x))+2k(\epsilon+k_{y}+k_{x})$.
Combining the bounds of $e_{G_{2}}(I_{1},D_{1}\cup D_{2})$ and grouping all the terms onto the right hand side, we deduce that 
\begin{align}
    0&\geq  k(\Delta_{1}-2)(2(\Delta_{2}-1)-|\hat{W}|)+kq-(\Delta_{1}-1)(\Delta_{2}-1)-|D'_{2}|+c\notag\\
    &\geq (\Delta_{1}-2)((2k-1)(\Delta_{2}-1)-2|\hat{W}|)-(\Delta_{2}-1)+kq-|D'_{2}|+c.\notag
\end{align}
Recall that $k\geq 2$ and $q\geq |D'_{2}||Z|$. Therefore, to avoid a contradiction the last inequality can only be true if $\Delta_{1}=3$, $|\hat{W}_{I_{1}}|=\Delta_{2}-1$, $k=2$, $D_{2}=\emptyset$, $q=0$, and $c=0$. We showed that every vertex in $D_{1}$ is adjacent in $G_{2}$ to a vertex in $D_{1}\cup D_{2}$. Since $D_{2}=\emptyset$ and $\Delta_{1}=3$, $D_{1}$ has exactly two vertices $u$ and $u'$ with degree one in $H$ that are adjacent in $G_{2}$ to themselves and every vertex in $I_{1}$. This implies $|I_{1}|=\Delta_{2}-2$ since (\ref{eq:I1}) says $|I_{1}|\geq \Delta_{2}-2$. Since $q=0$ and $|Z|=\Delta_{2}-1$, $W[x,y]$ induces a $K^{\Delta_{2}+1}$ in $G_{2}$. Thus, neither $u$ nor $u'$ are adjacent in $G$ to vertices in $Z$. By (\ref{eq:Idel1}), $|I|\geq \Delta_{2}$. Therefore, there is a $w'\in I-I_{1}$ not adjacent in $G_{2}$ to $u$. Let $v\in N_{H}(w')$ and $w\in N_{H}(v)\cap I_{1}$. We can now finish the proof of this section by constructing a packing of $\pi(G_{1})$ and $G_{2}$. We begin by interchanging $w$ with $w'$ in $H'$ to construct a realization $H'$ of $\pi(G_{1})$. Following this, we exchange in $H'$ the edges $xu$, $xy$, and $w'u$ with the non-edges $yu$, $xv$, and $xw'$ of $G$. Then we complete our construction by interchanging $x$ with $w'$ in the resulting realization.
\subsubsection{\label{sec:deltageq2}$\delta_{+}\geq 2$}
We choose an $(H,y)\in \mathcal{H}$ such that we give preference to one in $\Acute{\mathcal{H}}\cap \mathcal{H}$ if it exists. Furthermore, we choose any $x\in V$ that forms a bad pair with $y$.

\begin{claim}\label{cl:IcapL}$|I_{G}(x,y)\cap L(H)|=(\delta_{+}-1)(\Delta_{2}-2)-\eta$ for some non-negative $\eta\geq 0$.
\begin{proof}If $\Acute{\mathcal{H}}\cap \mathcal{H}=\emptyset$, then Claim~\ref{cl:L(H)bound}, implies $|I_{G}(x,y)\cap L(H)|=(\delta_{+}-1)(\Delta_{2}-2)-\eta$ for some $\eta>\Big\lfloor\frac{\Delta_{2}-2}{\delta_{+}}\Big\rfloor$. If $(H,y)\in \Acute{\mathcal{H}}\cap \mathcal{H}$, then every vertex in $I\cap L(H)$ must be adjacent in $G_{2}$ to a vertex in $S_{2}(H,y)$, but not adjacent in $G$ to $x$ nor $y$. Thus, $|I\cap L(H)|\leq e_{G_{2}}(I\cap L(H), S_{2}(H,y))\leq (\delta_{+}-1)(\Delta_{2}-2)$. 
\end{proof}
\end{claim}

We can use (\ref{cl:IcapL}) and $|W(x,y)|=|\hat{W}|+\epsilon$ to show that 
\[\begin{split}
|I-L|=|I|-|I\cap L|\geq \Delta_{2}(\Delta_{1}-1)+\Delta_{1}-&deg_{H}(x)+3-m-g-\delta_{+}\\&+\hat{W}+\epsilon+k_{y}+k_{x}-((\delta_{+}-1)(\Delta_{2}-2)-\eta).\notag
\end{split}\]
This reduces to \begin{equation}|I-L|\geq \Delta_{2}(\Delta_{1}-\delta_{+})+\Delta_{1}-deg_{H}(x)-m-g+|\hat{W}|+1+\delta_{+}+\eta+\epsilon+k_{y}+k_{x}.\label{eq:I-Lprevious}\end{equation}

\begin{claim}\label{cl:I-Lnotempty}$I-L$ is not empty.
\begin{proof}By contradiction, suppose $I-L$ is empty. By Claim~\ref{cl:InotEmpyt},  $I$ is not empty, and therefore there is an  $l\in I\cap L$. By definition, $\hat{W}$ would be empty, and from (\ref{cl:minW}) we can see that $\epsilon\geq \delta_{+}-1$. We now take advantage of $\delta_{+}$ being at least two. We first consider the case $\Delta_{1}>\delta_{+}$. Thus, $g=\Delta_{2}-1$. By (\ref{eq:I-Lprevious}), we see that \[|I-L|\geq \Delta_{2}(\Delta_{1}-\delta_{+})-m-(\Delta_{2}-1)+1+\delta_{+}+\epsilon\geq \Delta_{2}(\Delta_{1}-\delta_{+}-1)+2\delta_{+}+1-m.\] If $\Delta_{1}\geq \delta_{+}+2$, then $I-L$ is not empty since $m<\Delta_{2}+1+2\delta_{+}$. If $\Delta_{1}=\delta_{+}+1$, then $I-L$ is not empty since $2\delta_{+}+1=\Delta_{1}+\delta_{+}>m$. Thus, $\Delta_{+}=\delta_{+}$, and therefore, by (\ref{eq:I-Lprevious}), $|I-L|\geq 1+\Delta_{1}-m+\epsilon\geq \Delta_{1}-1\geq 1$. This proves Claim~\ref{cl:I-Lnotempty}.\qedhere
\end{proof}
\end{claim}

\begin{claim}\label{cl:D2lower}$\Delta_{2}\geq \delta_{+}+1$
\begin{proof}By Claim~\ref{cl:IcapL}, $|I\cap L|=(\delta_{+}-1)(\Delta_{2}-2)-\eta$ for some non-negative $\eta\geq0$. By contradiction, suppose $\Delta_{2}\leq \delta_{+}$. Thus, $m=\Delta_{2}+1$. By Claim~\ref{cl:InotEmpyt} and Claim~\ref{cl:I-Lnotempty}, $|I-L|\neq \emptyset$. Since $N_{H}(w)\subseteq \hat{W}$ for all $w\in I-L$, we deduce that $|\hat{W}|\geq \delta_{+}\geq \Delta_{2}$. Since $\hat{W}\cup \{x,y\}$ is complete, we know that every vertex in $\hat{W}$ is adjacent in $H$ to at least one vertex in $\hat{W}\cup \{x,y\}$. Thus, $q\geq |\hat{W}|$ and $e_{H}(I,\hat{W})\leq (\Delta_{1}-1)|\hat{W}|$.
After reducing (\ref{eq:I-Lprevious}), we see that \[|I-L|\geq \Delta_{2}(\Delta_{1}-\delta_{+}-1)+\delta_{+}-g+|\hat{W}|.\] Combining this with the upper bound of $e_{H}(I, \hat{W})$ we see that \[(\Delta_{1}-1)|\hat{W}|\geq \delta_{+}(\Delta_{2}(\Delta_{1}-\delta_{+}-1)+\delta_{+}-g+|\hat{W}|).\] Combining all the terms onto the right hand side and simplifying we deduce that \[0\geq (\Delta_{1}-\delta_{+}-1)(\delta_{+}\Delta_{2}-|\hat{W}|)+\delta_{+}(\delta_{+}-g).\] If $\Delta_{1}>\delta_{+}$, then we have a contradiction since $\delta_{+}>g$. If $\Delta_{1}=\delta_{+}$, then since $g=0$ and $|\hat{W}|\geq \delta_{+}$, we have the contradiction \[0\geq \delta_{+}^{2}-\delta_{+}\Delta_{2}+|\hat{W}|\geq \delta_{+}(\delta_{+}-\Delta_{2}+1)\geq 1.\] This completes the proof of Claim~\ref{cl:D2lower}
\end{proof}
\end{claim}

Since $\delta_{+}\geq 2$, Claim~\ref{cl:D2lower} implies $\Delta_{2}\geq\delta_{+}+1\geq 3$.

Suppose $\Delta_{1}=\delta_{+}$. In this case, $g=0$ and \[|I-L|\geq |\hat{W}|+\Delta_{1}-deg_{H}(x)+\eta+\epsilon+k_{y}+k_{x}.\] Using (\ref{eq:bounds}) and arranging terms onto the left hand side we deduce that \[\delta_{+}(\Delta_{1}-deg_{H}(x)+\eta+\epsilon+k_{y}+k_{x})+q\leq 0.\] However, this is a contradiction since Claim~\ref{cl:epAnddelta} shows that $\epsilon+q\geq 1$. Thus, $\Delta_{1}>\delta_{+}$.

Since $\Delta_{1}>\delta_{+}$ and $g=\Delta_{2}-1$, we may further reduce (\ref{eq:I-Lprevious}) to
\begin{equation}\label{eq:I-L}
    |I-L|\geq \Delta_{2}(\Delta_{1}-\delta_{+}-1)+\Delta_{1}-deg_{H}(x)-m+|\hat{W}|+2+\delta_{+}+\eta+\epsilon+k_{y}+k_{x}.
\end{equation}

Combining (\ref{eq:bounds}) and (\ref{eq:I-L}) we can show that \[\Delta_{1}|\hat{W}|-q\geq \delta_{+}(\Delta_{2}(\Delta_{1}-\delta_{+}-1)+\Delta_{1}-deg_{H}(x)-m+|\hat{W}|+2+\delta_{+}+\eta+\epsilon+k_{y}+k_{x}).\] Combining the $|\hat{W}|$ terms onto the left hand side we have the useful arrangement \begin{equation}
    (\Delta_{1}-\delta_{+})|\hat{W}|-q\geq \delta_{+}(\Delta_{2}(\Delta_{1}-\delta_{+}-1)+\Delta_{1}-deg_{H}(x)-m+2+\delta_{+}+\eta+\epsilon+k_{y}+k_{x}). \label{eq:boundsAdjusted}
\end{equation} 

\begin{claim}\label{cl:WlessD2}$|\hat{W}|\leq \Delta_{2}-1$
\begin{proof}By contradiction, we suppose $|\hat{W}|=\Delta_{2}+r$ for some non-negative integer $r$. We can see that every vertex in $\hat{W}$ is adjacent in $G$ to $\Delta_{2}+r+1$ vertices in $\hat{W}\cup \{x,y\}$. This implies every vertex $u\in \hat{W}$ is adjacent in $H$ to at least $r+1+e_{G_{2}}(u, L\cup (Y-X))$ vertices in $\hat{W}\cup \{x,y\}$. Therefore, $q=(r+1)|\hat{W}|+c$ for some non-negative integer $c\geq e_{G_{2}}(\hat{W}, L\cup (Y-X))+e_{H}(\hat{W},Y-X)$. We need to show that $\epsilon+\eta+c\geq 1$. Suppose $\epsilon+\eta=0$. In this case, $|L\cap I|=(\delta_{+}-1)(\Delta_{2}-2)$. Since $\delta_{+}\geq 2$ and $\Delta_{2}\geq 3$, there is some $l\in L\cap I$. By Claim~\ref{cl:MainObservation}, $l$ is adjacent in $G_{2}$ to a $u\in N_{H}(y)-\{x\}$. If $u\in \hat{W}$, then $c\geq 1$. Suppose $u\in Y-X$, and let $u'\in N_{H}(u)-\{y\}$. If $u'\in \hat{W}$, then $c\geq 1$. If $u'\in Y-\{y\}-\hat{W}$, then by Claim~\ref{cl:setup2}, $c\geq 1$ since every vertex in $\hat{W}$ is adjacent in $G$ to $u$. Thus, $c\geq 1$ when $\epsilon+\eta=0$, and therefore, $\epsilon+\eta+c\geq 1$.

Returning to (\ref{eq:boundsAdjusted}) with our new bound for $q$, we deduce \begin{equation}\label{eq:WlessD2_2}
\delta_{+}(\Delta_{2}(\Delta_{1}-\delta_{+}-1)+\Delta_{1}-deg_{H}(x)-m+2+\delta_{+}+\eta+\epsilon+k_{y}+k_{x})\leq (\Delta_{1}-\delta_{+}-r-1)|\hat{W}|-c.\end{equation} Since $m\leq \Delta_{1}+1$ and $r\geq 0$,
we can rearrange the left hand side of (\ref{eq:WlessD2_2}) and bound the right hand side such that \begin{equation}\label{eq:WlessD2_3}\delta_{+}(\Delta_{2}-1)(\Delta_{1}-\delta_{+}-1)+\delta_{+}(\Delta_{1}-deg_{H}(x)+\eta+\epsilon+k_{y}+k_{x})\leq (\Delta_{1}-\delta_{+}-1)|\hat{W}|-c.\end{equation} Since $\Delta_{2}\geq 3$, $\delta_{+}\geq 2$, and $|\hat{W}|\leq \Delta_{2}+r\leq \Delta_{2}-1+\delta_{+}-1,$ we may conclude that $|\hat{W}|\leq \delta_{+}(\Delta_{2}-1)$. Substituting this new upper bound for $|\hat{W}|$ in (\ref{eq:WlessD2_3}), we clearly see a contradiction since $\epsilon+\eta+c\geq 1$. Thus, $|\hat{W}|\leq \Delta_{2}-1$, and we have proved Claim~\ref{cl:WlessD2}.
\end{proof}
\end{claim}

\begin{claim}$\Delta_{1}=\delta_{+}+1$
\begin{proof}By contradiction, we suppose $\Delta_{1}\geq \delta_{+}+2$. Since $m\leq \Delta_{1}+1$ and $|\hat{W}|\leq \Delta_{2}-1$, we may deduce from (\ref{eq:boundsAdjusted}) that \[\delta_{+}(\Delta_{2}-1)(\Delta_{1}-\delta_{+}-1)+\delta_{+}(\eta+\epsilon+k_{y}+k_{x})\leq (\Delta_{1}-\delta_{+})(\Delta_{2}-1)-q.\] Gathering all of the terms on to the left hand side we see that \[(\Delta_{2}-1)((\Delta_{1}-\delta_{+}-1)(\delta_{+}-1)-1)+\delta_{+}(\eta+\epsilon+k_{y}+k_{x})+q\leq 0.\] Since $\Delta_{2}\geq 3$ and $\Delta_{1}-\delta_{+}\geq 2$, we see that $2(\delta_{+}-2)+\delta_{+}(\eta+\epsilon+k_{y}+k_{x})+q\leq 0.$ Claim~\ref{cl:epAnddelta} says $\epsilon +q\geq \delta_{+}-1\geq 1$. Therefore, we can deduce the contradiction \[1\leq \delta_{+}-2+\delta_{+}(\eta+\epsilon+k_{y}+k_{x})+q\leq 0.\] Therefore, $\Delta_{1}=\delta_{+}+1$.
\end{proof}
\end{claim}

Since $\Delta_{1}=\delta_{+}+1$,  Claim~\ref{cl:D2lower} implies $\Delta_{2}\geq \Delta_{1}=\delta_{+}+1$, and therefore, $m=\Delta_{1}+1$. We may now simplify (\ref{eq:boundsAdjusted}) to 
 \begin{equation}
    \delta_{+}(\Delta_{1}-deg_{H}(x)+\eta+\epsilon+k_{y}+k_{x})\leq |\hat{W}|-q. \label{eq:boundsAdjustedFurther}
\end{equation} 

\begin{claim}\label{cl:etqBound}$\eta\leq \bigg\lfloor\frac{\Delta_{2}-2}{\delta_{+}}\bigg\rfloor$ and $(H,y)\in \Acute{\mathcal{H}}\cap \mathcal{H}$
\begin{proof} Since $|\hat{W}|\leq \Delta_{2}-1$ and $\epsilon+q\geq \delta_{+}-1\geq 1$, we may deduce from  (\ref{eq:boundsAdjustedFurther}) that 
$\delta_{+}\eta+1\leq \Delta_{2}-1$,  and therefore, $\eta\leq \bigg\lfloor\frac{\Delta_{2}-2}{\delta_{+}}\bigg\rfloor$. This implies \[(\delta_{+}-1)(\Delta_{2}-2)-\bigg\lfloor\frac{\Delta_{2}-2}{\delta_{+}}\bigg\rfloor\leq |I_{G}(x,y)\cap L(H)|\leq |L(G_{1})|.\] By Claim~\ref{cl:L(H)bound}, $\Acute{\mathcal{H}}\cap \mathcal{H}$ is not empty. Thus, 
$(H,y)\in \Acute{\mathcal{H}}\cap \mathcal{H}$ by our choice of $(H,y)$.
\end{proof}
\end{claim}

\begin{claim}\label{cl:eta2}Every vertex in $N_{H}(y)-\{x\}$ is adjacent in $G_{2}$ to at most $\eta+2$ vertices in $V-(I\cap L)$.
\begin{proof}Let $|S_{2}(H,y)-\{x\}|=\delta_{+}-c$ for some positive integer $c\geq 1$.Let $t_{y}$ be the number of vertices in $S_{2}(H,y)-\{x\}$ adjacent in $G_{2}$ to at least $\eta+3$ vertices in $V-(I\cap L)$. Those $t_{y}$ vertices are adjacent in $G_{2}$ to at most $\Delta_{2}-\eta-3$ vertices in $I\cap L$. Since $(H,y)\in \Acute{\mathcal{H}}\cap \mathcal{H}$, every vertex in $I\cap L$ is adjacent in $G_{2}$ to a vertex in $S_{2}(H,y)-\{x\}$. Therefore, we may deduce that
\begin{align}
    (\delta_{+}-1)(\Delta_{2}-2)-\eta&\leq e_{G_{2}}(I\cap L, S_{2}(H,y))\notag\\
    &\leq t_{y}(\Delta_{2}-\eta-3)+(|S_{2}(H,y)|-t_{y})(\Delta_{2}-2)\notag\\
    &= t_{y}(\Delta_{2}-\eta-3)+(\delta_{+}-c-t_{y})(\Delta_{2}-2)\notag \\
    &\leq (\delta_{+}-1)(\Delta_{2}-2)-t_{y}(\eta+1)-(c-1)(\Delta_{2}-2).\notag
\end{align} Since $\eta<\Delta_{2}-2$, we may conclude that $t_{y}=0$ and $c=1$. Since $c=1$, $S_{2}(H,y)-\{x\}=N_{H}(y)-\{x\}$. Thus, the proof is complete.
\end{proof}
\end{claim}

\begin{claim}\label{cl:WIG2Nei}Every vertex in $\hat{W}$ is adjacent in $G_{2}$ to at least $\eta\delta_{+}+2\delta_{+}-1$ vertices in $N_{G}[x,y]$.
\begin{proof}We choose a $u\in\hat{W}$. Observe that \[q=2e_{H}(\hat{W},\hat{W})+e_{H}(\hat{W},N[x,y]-\hat{W}).\] Since 
$e_{H}(\hat{W},\hat{W})\geq e_{H}(u,\hat{W})$, we see that \[\frac{q-e_{H}(\hat{W},N[x,y]-\hat{W})}{2}\geq e_{H}(u,\hat{W}).\] From (\ref{eq:boundsAdjustedFurther}), we see that
\begin{align}
    |N_{G_{2}}(u)\cap \hat{W}|&=|\hat{W}-\{u\}|-e_{H}(u,\hat{W})\notag\\
    &\geq \delta_{+}(\eta+\epsilon+k_{y}+k_{x})-1+q-\frac{q-e_{H}(\hat{W},N[x,y]-\hat{W})}{2}.\notag
\end{align}
Since $q\geq e_{H}(\hat{W},N[x,y]-\hat{W})\geq e_{H}(u,N[x,y]-\hat{W})$, we see that \[|N_{G_{2}}(u)\cap \hat{W}|\geq \delta_{+}(\eta+\epsilon+k_{y}+k_{x})+e_{H}(\hat{W},N[x,y]-\hat{W})-1.\] Applying the previous inequality to $|N_{G_{2}}(u)\cap N[x,y]|=|N_{G_{2}}(u)\cap \hat{W}|+e_{G_{2}}(u,N[x,y]-\hat{W})$, we establish
\begin{equation}|N_{G_{2}}(u)\cap N[x,y]|\geq\delta_{+}(\eta+\epsilon+k_{y}+k_{x})+e_{H}(\hat{W},N[x,y]-\hat{W})+e_{G_{2}}(u,N[x,y]-\hat{W})-1.\label{eq:blueNieghborsWI}
\end{equation}

To finish this proof we need to establish a lower bound for $e_{H}(\hat{W},N[x,y]-\hat{W})+e_{G_{2}}(u,N[x,y]-\hat{W})$. Since $|\hat{W}|\leq \Delta_{2}-1$, we see that $|X|+|Y|-2|W_{G}[x,y]|\geq 2(\Delta_{2}+\delta_{+})-2(\Delta_{2}-1+2))-2(k_{y}+k_{x}) \geq 2(\delta_{+}-1)-2(k_{y}+k_{x})$. Let $r$ be the number of vertices in $(X-Y)\cup (Y-X)$ adjacent in $H$ to a vertex in $\hat{W}$. If $v\in (X-Y)\cup (Y-X)$ is not adjacent in $H$ to vertices in $\hat{W}$, then since $\delta_{+}\geq 2$, there is a $v'\in N_{H}(u)$ such that $v'\notin W[x,y]\cup L$. Therefore, by Claim~\ref{cl:setup2}, $u$ is adjacent in $G_{2}$ to $v$. Thus, $u$ is adjacent in $G_{2}$ to $2(\delta_{+}-1)-2(k_{y}+k_{x})-r$ vertices in $(X-Y)\cup (Y-X)$. Since $u$ is adjacent in $G$ to both $x$ and $y$, we may conclude that \begin{align}e_{H}(\hat{W}, N[x,y]-\hat{W})+e_{G_{2}}(u, N[x,y]-\hat{W})&\geq r+|\{x,y\}|+2(\delta_{+}-1)-2(k_{y}+k_{x})-r\notag\\
 &\geq 2\delta_{+}-2(k_{y}+k_{x}).\notag
 \end{align} 
Applying this lower bound to (\ref{eq:blueNieghborsWI}), we prove our claim \[|N_{G_{2}}(u)\cap N[x,y]|\geq \delta_{+}\eta+\delta_{+}(k_{y}+k_{x}) +2\delta_{+}-2(k_{y}+k_{x})-1\geq \eta\delta_{+}+2\delta_{+}-1.\qedhere\]
\end{proof}
\end{claim}

\begin{claim}\label{cl:NoYNeighInWI}$N_{H}(y)\cap \hat{W}=\emptyset$
\begin{proof}If there is a $u\in N_{H}(y)\cap \hat{W}$, then Claim~\ref{cl:WIG2Nei} says $u$ is adjacent in $G_{2}$ to at least $\eta+3$ vertices in $N_{G}[x,y]$. However, this contradicts Claim~\ref{cl:eta2}. 
\end{proof}
\end{claim}

\begin{claim}\label{cl:uNeigWI}If $\epsilon=0$ and $u\in N_{H}(y)\cap (Y-X)$, then  $e_{G}(u,\hat{W})\geq \delta_{+}-1$.
\begin{proof}Since $\epsilon=0$, $W_{G}(x,y)=\hat{W}$. Suppose there is a $u'\in N_{H}(u)$ not in $\hat{W}$. By Claim~\ref{cl:setup} and Claim~\ref{cl:setup2}, $u'\in Y-X$ and every vertex in $\hat{W}$ is adjacent in $G$ to $u'$. Since $N_{H}(w)\subseteq \hat{W}$ for every $w\in I-L$, $u$ must be adjacent in $G$ to $\delta_{+}$ vertices in $\hat{W}$. Therefore, $e_{G}(u,\hat{W})\geq e_{H}(u,\hat{W}) \geq \delta_{+}-1$.
\end{proof}
\end{claim}

\begin{claim}\label{cl:N(x)WI}If $\Delta_{1}-deg_{H}(x)+\epsilon+k_{x}+k_{y}=0$, then $e_{H}((X-Y)\cup \{x\}, \hat{W})\geq \delta_{+}$.
\begin{proof}Since $\epsilon=0$ and $\delta_{+}\geq 2$, Claim~\ref{cl:NoYNeighInWI} says that $N_{H}(y)-\{x\}\subseteq Y-X$. We choose a $u\in N_{H}(y)\cap (Y-X)$. Claim~\ref{cl:setup} says that $N_{H}(v)\subseteq \hat{W}$ for any $v\in N_{H}(x)\cap (X-Y)$. Since $\Delta_{1}=deg_{H}(x)$ and $\Delta_{1}=\delta_{+}+1$, $|N_{H}(x)-\{y\}|=\delta_{+}$. Therefore,
\begin{align}
    e_{H}((X-Y)\cup \{x\}, \hat{W})&=e_{H}((X-Y)\cup \{x\}, \hat{W})+e_{H}(x,\hat{W})\notag\\
    &\geq (\delta_{+}-1)(\delta_{+}-|N_{H}(x)\cap \hat{W}|)+|N_{H}(x)\cap \hat{W}|\geq \delta_{+}.\notag \qedhere
\end{align}
\end{proof}
\end{claim}

\begin{claim}\label{cl:speacialwinI-L}There exists a $w\in I-L$ not adjacent in $G$ to any vertex in $N_{H}(y)$.
\begin{proof}We need to show that $|I-L|-e_{G_{2}}(N_{H}(y), I-L)>0$. In Claim~\ref{cl:NoYNeighInWI} we showed that $N_{H}(y)\cap \hat{W}=\emptyset$. Thus, no vertex in $I-L$ is adjacent in $H$ to a vertex in $N_{H}(y)$. By (\ref{eq:boundsAdjustedFurther}), $|\hat{W}|\geq \delta_{+}(\Delta_{1}-deg_{H}(x)+\eta+\epsilon+k_{y}+k_{x})+q$. Combining this with (\ref{eq:I-L})  we deduce that 
\begin{align}
    |I-L|&\geq |\hat{W}|+\Delta_{1}-deg_{H}(x)+\eta+\epsilon+k_{y}+k_{x}\notag\\
    &\geq \delta_{+}(\Delta_{1}-deg_{H}(x)+\eta+\epsilon+k_{y}+k_{x})+q+\Delta_{1}-deg_{H}(x)+\eta+\epsilon+k_{y}+k_{x}\notag\\
    &\geq (\delta_{+}+1)(\Delta_{1}-deg_{H}(x)+\eta+\epsilon+k_{y}+k_{x})+q.\notag
\end{align}
Recall that every vertex in $N_{H}(y)$ is adjacent in $G_{2}$ to at most $\eta+2$ vertices in $V-(I\cap L)$. Thus, \[e_{G_{2}}(N_{H}(y), I-L)\leq (\delta_{+}-1)(\eta+2)-e_{G_{2}}(N_{H}(y),\hat{W}).\] 
Suppose $\Delta_{1}-deg_{H}(x)+\epsilon+k_{x}+k_{y}=0$. This implies \[|I-L|\geq \eta(\delta_{+}+1)+q\geq\eta(\delta_{+}+1)+e_{H}(\hat{W},N[x,y]-\hat{W}).\] Let $u\in N_{H}(y)\cap (Y-X)$. By Claim~\ref{cl:uNeigWI} and Claim~\ref{cl:N(x)WI}, we deduce that \[e_{H}(\hat{W},N[x,y]-\hat{W})\geq e_{H}((X-Y)\cup \{x\}, \hat{W})+e_{G}(u, \hat{W})-e_{G_{2}}(u,\hat{W})\geq 2\delta_{+}-1-e_{G_{2}}(u,\hat{W}).\] Since $e_{G_{2}}(N_{H}(y),\hat{W})\geq e_{G_{2}}(u, \hat{W})$, we see that \[e_{G_{2}}(N_{H}(y), I-L)\leq (\delta_{+}-1)(\eta+2)-e_{G_{2}}(u,\hat{W}).\] Thus, 
\begin{align}
    |I-L|-e_{G_{2}}(N_{H}(y), I-L)&\geq \eta(\delta_{+}+1)+e_{H}(\hat{W},N[x,y]-\hat{W})-e_{G_{2}}(N_{H}(y), I-L)\notag\\
    &\geq \eta(\delta_{+}+1)+2\delta_{+}-1-e_{G_{2}}(u,\hat{W})-((\delta_{+}-1)(\eta+2)-e_{G_{2}}(u,\hat{W}))\notag\\
    &\geq 2\eta+1\geq 1.\notag
\end{align}

We are left with the case $\Delta_{1}-deg_{H}(x)+\epsilon+k_{x}+k_{y}\geq 1$. Claim~\ref{cl:epAnddelta} says $q\geq \delta_{+}-1$. Thus, 
\begin{align}
|I-L|-e_{G_{2}}(N_{H}(y), I-L)&\geq  (\eta+1)(\delta_{+}+1)+(\delta_{+}-1)- (\delta_{+}-1)(\eta+2)\notag \\
&=2\eta+2\geq 2\notag.
\end{align}This finishes the proof of Claim~\ref{cl:speacialwinI-L}.
\end{proof}
\end{claim}
We choose a $w\in I-L$ that is not adjacent to vertices in $N_{H}(y)$. Let $z\in N_{H}(w)$, and let $L_{1}\subseteq L$ be the vertices adjacent in $G_{2}$ to exactly one vertex in $N_{H}(y)$. Suppose $z$ is adjacent to every vertex in $L_{1}(H)\cap I$. Claim~\ref{cl:WIG2Nei} says that $z$ is adjacent to at most $\Delta_{2}-(\eta\delta_{+}+\delta_{+}+1)$ vertices in $I\cap L$. Thus, $|L_{1}(H)\cap I|\leq \Delta_{2}-(\eta\delta_{+}+\delta_{+}+1)$. By Claim~\ref{cl:MainObservation}, \[2|I\cap L|-|L_{1}(H)\cap I|\leq e_{G_{2}}(N_{H}(y),I\cap L).\] Since $|I\cap L|=(\delta_{+}-1)(\Delta_{2}-2)-\eta$, we may deduce that 
\begin{align}
    e_{G_{2}}(N_{H}(y),I\cap L)&\geq 2((\delta_{+}-1)(\Delta_{2}-2)-\eta)-(\Delta_{2}-(\eta\delta_{+}+\delta_{+}+1)\notag\\
    &\geq \Delta_{2}(2\delta_{+}-3)+\eta(\delta_{+}-2)+5-3\delta_{+}.\notag
\end{align}
Since $(H,y)\in \Acute{\mathcal{H}}$, there has to be at least one vertex in $N_{H}(y)-\{x\}$ that is adjacent to at most $\Delta_{2}-2$ vertices in $I\cap L$. Therefore, $e_{G_{2}}(N_{H}(y),I\cap L)\leq (\delta_{+}-1)\Delta_{2}-2$.
Combining the lower and upper bounds of  $e_{G_{2}}(N_{H}(y),I\cap L)$ we see that \[(\delta_{+}-1)\Delta_{2}-2\geq \Delta_{2}(2\delta_{+}-3)+\eta(\delta_{+}-2)+5-3\delta_{+}.\]
Moving all terms onto the right hand side and then collecting them we see that
\[0\geq \Delta_{2}(\delta_{+}-2)+\eta(\delta_{+}-2)+7-3\delta_{+}.\]
Since $\Delta_{2}\geq \delta_{+}+1$ and $\delta_{+}\geq 2$, we deduce the contradiction 
\begin{align}
    0&\geq (\delta_{+}+1)(\delta_{+}-2)+\eta(\delta_{+}-2)+7-3\delta_{+}\notag\\
    &=\eta(\delta_{+}-2)+\delta_{+}(\delta_{+}-4)+5\geq 1.\notag
\end{align} Therefore, there is an $l\in L_{1}\cap I$ that is not adjacent in $G$ to $z$.

Let $u$ be the sole neighbor of $l$ in $N_{H}(y)$. We exchange in $H$ the edges $wz$ and $yu$ with the non-edges $zy$ and $uw$ of $G$. Then in the resulting realization we interchange $y$ with $l$ to create a packing of $\pi(G_{1})$ with $G_{2}$. This final contradiction completes the proof of Theorem~\ref{theorem:BEC-half}.
\end{theoremlateproof}
\subsubsection{\label{sec:modify}Proof of Claim~\ref{cl:modify}}

\begin{claimlateproof}{cl:modify}We are given an $H\in \mathcal{R}(\pi(G_{1}))$ and a $y\in V$ with $deg_{H}(y)=\delta_{+}$. We define $G=H+ G_{2}$, $b=b(H, G_{2})$, and $g=g(H,G_{2})$. Let $S=N_{H}(y)$, $S_{0}=S_{0}(H,y)$, $S_{1}=S_{1}(H,y)$, and $S_{2}=S_{2}(H,y)$. We define $L=L_{0}\cup L_{1} \cup L_{2}$ to be a partition of $L$ such that every vertex in $L_{2}$ has at least two neighbors in $S$ and for $i\in \{0,1\}$, every vertex in $L_{i}$ has $i$ neighbors in $S$. By Claim~\ref{cl:MainObservation}, $L_{0}=\emptyset$. For any $u\in V$ and $i\in \{0,1,2\}$, we let $\Gamma(u)=N_{G_{2}}(u)\cap L$ and $\Gamma_{i}(u)=N_{G_{2}}(u)\cap L_{i}$. In addition, we need to define the set $\Gamma'_{1}(u)=\Gamma_{1}(u)\cup \{y\}$.

Suppose there is a smallest number $t$ in $\{0,1\}$ and an $l\in L(H)\cup \{y\}$ that satisfies (\ref{eq:modify}). In search of a contradiction, we assume there is an $s\in N_{G_{2}}(l)\cap S_{t'}$ for some $t'\leq t$ that violates Claim~\ref{cl:modify}.  

Our proof of Claim~\ref{cl:modify} is made easier if $l\neq y$. So, suppose $l=y$ and there is an $l'\in \Gamma(s)$ that satisfies (\ref{eq:modify}). If we prove Claim~\ref{cl:modify} holds for $l'$ in $H$ with parameters $y$, $t$, $t'$, and $s$, then there are edges $\{u_{1}u_{1}',\ldots, u_{q}u_{q}'\}$ that satisfy conditions \ref{claim2con1}, \ref{claim2con2}, and \ref{claim2con3} where $u_{1}=y$ and $u_{1}'=s$. Consequently, these same edges and non-edges prove the claim for when $y=l$. Therefore, if $l=y$, then we will assume every vertex in $\Gamma(s)$ does not satisfy (\ref{eq:modify}), and thus, every vertex in $\Gamma(s)$ is adjacent in $G_{2}$ to vertices in $S_{2}$. Observe, that if $\delta_{+}=1$, then $l\neq y$ since $\{s\}=S$ and $L$ is not empty.

\begin{subclaim}\label{scl:basicObservation}For $u\in S$ and $z\in \Gamma_{1}'(u)$, if $z'\in (L\cup \{y\})-N_{G_{2}}(z)$, then $u$ is adjacent in $H$ to some vertex in $N_{G_{2}}(z')$.
\begin{proof}If $u$ is not adjacent in $H$ to some vertex in $N_{G_{2}}(z')$, then we may interchange $y$ with $z$ in $H$ to construct an $H'\in \mathcal{R}(\pi(G_{1}))$ with $b(H',G_{2})=1$. We may then interchange $u$ with $z'$ in $H'$ to construct a realization of $\pi(G_{1})$ that packs with $G_{2}$. This contradiction implies $u$ is adjacent in $H$ to some vertex in $N_{G_{2}}(z')$.
\end{proof}
\end{subclaim}

\begin{subclaim}\label{scl:wandl}If there is a set of edges $\{w_{0}w'_{0},\ldots, w_{q}w'_{q}\}$ of $H$ with $w_{0}=y$ and $w'_{q}\in V-L-N_{H}[y]$ such that $\{w'_{0}w_{1},\ldots, w'_{q-1}w_{q}\}$ are non-edges of $G$, then every vertex in $\Gamma'_{1}(w'_{0})$ is adjacent in $G_{2}$ to $w'_{q}$.
\begin{proof}Suppose there is a $z\in \Gamma'_{1}(w'_{0})$ not adjacent in $G_{2}$ to $w'_{q}$. We may exchange in $H$ the edges $\{w_{0}w'_{0},\ldots, w_{q}w'_{q}\}$ with the non-edges $\{w'_{0}w_{1},\ldots, w'_{q}w_{0}\}$ to construct a new realization $H'\in \mathcal{R}(\pi(G_{1}))$. Since $N_{H'}(y)\cap N_{G_{2}}(z)=\emptyset$, we may interchange $y$ with $z$ in $H'$ to construct a realization of $G_{1}$ that packs with $G_{2}$. Therefore, no such $z$ exists.
\end{proof}
\end{subclaim}

Since $L_{0}=\emptyset$ and some vertex in $S$ is adjacent in $G_{2}$ to $y$, $|L|\leq |S|\Delta_{2}-|\{y\}|=\delta_{+}\Delta_{2}-1$. Therefore, there is a natural number $\lambda$ such that  $|L|=\delta_{+}\Delta_{2}-\lambda$. We now establish a general lower bound for $\lambda$.

We observe
\begin{align}
    e_{G_{2}}(S,L)&=e_{G_{2}}(S,L_{1})+e_{G_{2}}(S,L_{2})\notag \\
    &= \left(\sum_{u\in S}d_{G_{2}}(u)\right)-(b+2e_{G_{2}}(S,S)+e_{G_{2}}(S,V-S-L-\{y\})).\notag
\end{align}
Furthermore, we know \[|L|=|L_{1}|+|L_{2}|=e_{G_{2}}(S,L_{1})+|L_{2}|.\] When we combine the last two equations and then solve for $|L|$, we establish \[|L|= \left(\sum_{u\in S}d_{G_{2}}(u)\right)-(b+e_{G_{2}}(S,L_{2})-|L_{2}|+2e_{G_{2}}(S,S)+e_{G_{2}}(S,V-S-L-\{y\})).\] Thus, \begin{equation}\label{eq:lambda}
\lambda=\delta_{+}\Delta_{2}-\left(\sum_{u\in S}d_{G_{2}}(u)\right)+ b+e_{G_{2}}(S,L_{2})-|L_{2}|+2e_{G_{2}}(S,S)+e_{G_{2}}(S,V-S-L-\{y\}).
\end{equation}
Note that $e_{G_{2}}(S,L_{2})-|L_{2}|\geq |L_{2}|$ and \[2e_{G_{2}}(S,S)+e_{G_{2}}(S,V-S-L-\{y\})\geq 2|S_{2}|+|S_{1}|=|S_{2}|+|S|-|S_{0}|.\]
We can bound (\ref{eq:lambda}) from below with
\begin{align}
    \lambda&\geq b+|S_{2}|+\delta_{+}-|S_{0}|+e_{G_{2}}(S,L_{2})-|L_{2}|\label{eq:lowerlambda}\\
    &\geq b+|S_{2}|+\delta_{+}-|S_{0}|+|L_{2}|.\label{eq:lowerlambdafurther}
\end{align}

Using (\ref{eq:main}) and $|L|$ we establish that $|V-L|\geq (\Delta_{1}+1)(\Delta_{2}+1)-m-g-(\delta_{+}\Delta_{2}-\lambda)$. This can be rewritten as \begin{equation}\label{eq:V-L}|V-L|\geq \Delta_{2}(\Delta_{1}-\delta_{+})-g +\Delta_{2}+\Delta_{1}+1 - m + \lambda.\end{equation}

\begin{subclaim}\label{scl:missingNeighbor}$V-L-N_{H}(u)-N_{G_{2}}(u)\neq \emptyset$ for every $u\in S_{1}\cup S_{0}$.
\begin{proof}Since $|V-L|=|G^{+}|\geq \Delta_{1}+2$, the claim is easily established when $u\in S_{0}$. Suppose $u\in S_{1}$. Since $|S_{0}|\leq \delta_{+}-1$, (\ref{eq:lowerlambdafurther}) implies that $\lambda\geq 2$. To finish the proof, it is enough to show $|V-L|\geq \Delta_{1}+3$ since $u$ is adjacent in $G$ to at most $\Delta_{1}+1$ vertices in $V-L$.  If $\Delta_{1} > \delta_{+}$, then (\ref{eq:V-L}) says $|V-L|\geq \Delta_{1}+3$. We are left with the case $\Delta_{1}=\delta_{+}$. Let $u'\in N_{G_{2}}(u)\cap (V-L)$. If $u'\in V-L-N_{H}[y]$, then since $u'\in N_{G_{2}}(u)$, there must be a $w\in N_{H}(u')\cap (V-L-N_{H}[y])$. Thus, $|V-L|\geq |N_{H}[y]|+|\{u',w\}|\geq \Delta_{1}+3$. Finally, if $u'\in S$, then (\ref{eq:lowerlambdafurther}) shows that $\lambda\geq 3$, and therefore, $|V-L|\geq \Delta_{1}+3$ by (\ref{eq:V-L}).
\end{proof}
\end{subclaim}
\begin{subclaim}
    $\Delta_{2}\geq 2$.
    \begin{proof}By Claim~\ref{scl:missingNeighbor}, there is a $w_{1}\in V-L-N_{H}(s)-N_{G_{2}}(s)$. Thus, there is an $w_{1}'\in N_{H}(w_{1})-N_{H}[y]$. However, since $\{ys, w_{1}w_{1}'\}$ are edges of $H$ and $sw_{1}$ is a non-edge, we see by Claim~\ref{scl:wandl} that $l$ is adjacent in $G_{2}$ to both $s$ and $w_{1}'$. Thus, $\Delta_{2}\geq 2$.
    \end{proof}
\end{subclaim}

\begin{subclaim}\label{cl:lower}$\lambda\geq \text{min}\{\delta_{+}+1,\Delta_{2}+1\}$.
\begin{proof}

We let $J=V-L-N_{H}[y]$. Observe that $J$ is not empty since Claim~\ref{cl:minOrder} shows $|V-L|\geq \Delta_{1}+2$. By contradiction, suppose $\lambda\leq \text{min}\{\delta_{+},\Delta_{2}\}$. 
 
\begin{subsubclaim}\label{scl:lower}$\Gamma'_{1}(u)\neq \emptyset$ for every $u\in S$.
\begin{proof}
If there is a vertex $u\in S$ such that $\Gamma_{1}'(u)=\emptyset$, then $|L_{2}|\geq e_{G_{2}}(u,L_{2})$ and  $e_{G_{2}}(u,L_{2})+e_{G_{2}}(u,S)+e_{G_{2}}(u, V-S-L-\{y\})\geq d_{G_{2}}(u)$. However, since $e_{G_{2}}(S,L_{2})-|L_{2}|\geq |L_{2}|\geq e_{G_{2}}(u,L_{2})$, we see by (\ref{eq:lambda}) the contradiction
\begin{align}
    \lambda&\geq \Delta_{2}-d_{G_{2}}(u)+ b+e_{G_{2}}(S,L_{2})-|L_{2}|+2e_{G_{2}}(S,S)+e_{G_{2}}(S,V-S-L-\{y\})\notag\\
    &\geq \Delta_{2}-d_{G_{2}}(u)+b +|L_{2}|+2e_{G_{2}}(S,S)+e_{G_{2}}(S,V-S-L-\{y\})\notag\\
    &\geq \Delta_{2}-d_{G_{2}}(u)+b + e_{G_{2}}(u,L_{2})+e_{G_{2}}(u,S)+e_{G_{2}}(u, V-S-L-\{y\})\notag\\
    &\geq \Delta_{2}-d_{G_{2}}(u)+b+d_{G_{2}}(u)\geq \Delta_{2}+1.\notag \qedhere
\end{align}
\end{proof}
\end{subsubclaim}
Since $\delta_{+}\geq \lambda$, (\ref{eq:lowerlambdafurther}) implies $|S_{0}|\geq e_{G_{1}}(S,L_{2})-|L_{2}|+1\geq |L_{2}|+1$. Thus, by (\ref{eq:lowerlambda}), there must be a vertex in $S_{0}$ adjacent in $G_{2}$ to at least $\Delta_{2}-1$ vertices in $L_{1}$. Since $|\Gamma_{1}'(u)|\leq \Delta_{2}-1$ for all $u\in S-S_{0}$, we may choose an $x\in S_{0}$ that has the largest $|\Gamma'_{1}(x)|$ over all vertices in $S$. 

\begin{subsubclaim}\label{subsubcl:delta=1}$\delta_{+}=1$.
\begin{proof}By contradiction, suppose $\delta_{+}\geq 2$. We choose an $s'\in S-\{x\}$ that yields a maximum $|\Gamma'_{1}(s')|$ such that we give preference to such a vertex in $S_{0}-\{x\}$ if it exists.

We need to show that $|\Gamma_{1}'(s')\cup \Gamma_{1}'(x)|\geq \min\{\Delta_{2}+1,2(\Delta_{2}-1)\}$ and $|\Gamma_{1}'(x)|=\Delta_{2}$ when $s'\notin S_{0}$. If $|\Gamma_{1}'(x)|=\Delta_{2}$, then $|\Gamma_{1}'(s')\cup \Gamma_{1}'(x)|\geq \Delta_{2}+1$ by Claim~\ref{scl:lower}. Suppose $|\Gamma_{1}'(x)|=\Delta_{2}-1$, and therefore, $|\Gamma_{2}'(x)|=1$. Let $|\Gamma_{1}'(s')|=\Delta_{2}-i$ for some natural number $i\geq 1$. For all $u\in S_{0}-\{x\}$, our choice of $s'$ implies that $|\Gamma_{1}'(u)|\geq i+1$ when $s'\notin S_{0}$ and $|\Gamma_{1}'(u)|\geq i$ when $s'\in S_{0}$. Therefore, if $i\geq 2$ or $s'\notin S_{0}$, then we  can deduce the contradiction $e_{G_{1}}(S,L_{2})\geq 2(|S_{0}|-1)+|\Gamma_{2}'(x)|\geq 2|S_{0}|-1\geq |S_{0}|+|L_{2}|\geq e_{G_{1}}(S,L_{2})+1$. This implies $|\Gamma_{1}'(s')|=\Delta_{2}-1$ and $s'\in S_{0}$, and therefore, $|\Gamma_{1}'(s')\cup \Gamma_{1}'(x)|\geq 2(\Delta_{2}-1)$.

Suppose no vertex in $J$ is adjacent in $G_{2}$ to every vertex in $\Gamma_{1}'(s')\cup \Gamma_{1}'(x)$. By Claim~\ref{scl:missingNeighbor}, there is a $u\in V-L-N_{H}(s')-N_{G_{2}}(s')$, and therefore, a $u'\in N_{H}(u)\cap J$. Since $\{ys', uu'\}$ are edges of $H$ and $s'u$ is a non-edge of $G$, Claim~\ref{scl:wandl} says $u'$ is adjacent in $G_{2}$ to every vertex in $\Gamma'_{1}(s')$. We choose a $w_{0}\in V-L-N_{H}[x]$ such that we give preference to one in $J$ if it exists. Let $w'_{0}\in N_{H}(w_{0})\cap J$. Since $\{yx, w_{0}w'_{0}\}$ are edges of $H$ and $xw_{0}$ is a non-edge of $G$, Claim~\ref{scl:wandl} says $\Gamma'_{1}(x)\subseteq N_{G_{2}}(w'_{0})$. By our choice of $w_{0}$, if $w'_{0}\notin N_{H}(x)$, then $w_{0}\in J$. If $w_{0}=s'$, then $w'_{0}\in N_{H}(x)$, and therefore, since $\{ys', xw'_{0}\}$ are edges of $H$ and $xs'$ is a non-edge of $G$, Claim~\ref{scl:wandl} implies the contradiction $\Gamma_{1}'(s')\cup \Gamma_{1}'(x)\subseteq N_{G_{2}}(w'_{0})$. Thus, $w_{0}\neq s'$. If $w_{0}\notin N_{G}(s')$, then since $\{ys', w_{0}w'_{0}\}$ are edges of $H$ and $ s'w'_{0}$ is a non-edge of $G$, Claim~\ref{scl:wandl} implies the contradiction $\Gamma_{1}'(s')\cup \Gamma_{1}'(x) \subseteq N_{G_{2}}(w'_{0})$. If $w_{0}\in N_{H}(s')$, then since $\{yx,w_{0}s', uu'\}$ are edges of $H$ and $\{s'u, xw_{0}\}$ are non-edges of $G$, Claim~\ref{scl:wandl} implies the contradiction $\Gamma_{1}'(s')\cup \Gamma_{1}'(x) \subseteq N_{G_{2}}(u')$. Therefore, $w_{0}\in N_{G_{2}}(s')$, and thus, $s'\in S_{1}$ and as shown previously $|\Gamma_{1}'(x)|=\Delta_{2}$. If $w'_{0}\notin N_{H}(x)$, then $w_{0}\in J$, and since $\{yx,w_{0}w'_{0}\}$ are edges of $H$ and $xw'_{0}$ is a non-edge of $G$, we conclude that $\Gamma_{1}'(x)\subseteq N_{G_{2}}(w_{0})$ by Claim~\ref{scl:wandl}. However, this is a contradiction since $|\Gamma_{1}'(x)|=\Delta_{2}$ and $s'\in N_{G_{2}}(w_{0})$. Therefore, $w'_{0}\in N_{H}(x)$. Since $\Gamma'_{1}(x)\subseteq N_{G_{2}}(w'_{0})$ and $|\Gamma_{1}'(x)|=\Delta_{2}$, we deduce that $w'_{0}\notin N_{G_{2}}(s')$. If $w'_{0}$ is not adjacent in $H$ to $s'$, then since $\{ys',w'_{0}x, w_{0}w'_{0}\}$ are edges of $H$ and $\{s'w'_{0}, xw_{0}\}$ are non-edges of $G$, Claim~\ref{scl:wandl} implies the contradiction $\Gamma_{1}'(s')\cup \Gamma_{1}'(x) \subseteq N_{G_{2}}(w_{0}')$. Thus, $w'_{0}$ is adjacent in $H$ to $x$, $s'$, and $w_{0}$. Since $s'\in S_{1}$, we conclude that $|S_{0}|=1$ by our choice of $s'$, and therefore, by (\ref{eq:lowerlambdafurther}), $\lambda\geq \delta_{+}\geq 2$. Since $w_{0}\in N_{G_{2}}(s')$ and $w_{0}$ is not adjacent in $G$ to $x$, we deduce that $w_{0}$ is adjacent in $H$ to two vertices in $J$. Thus, $|J|\geq 2$. Therefore, if $\Delta_{1}=\delta_{+}$, then  $|V-L|\geq \Delta_{1}+3$. Moreover, if $\Delta_{1}>\delta_{+}$, then by (\ref{eq:V-L}), we deduce that $|V-L|\geq \Delta_{1}+3$. Since $|\Gamma_{1}'(w'_{0})|=\Delta_{2}$, there must be a $w_{1}\in V-L-\{y\}$ not adjacent in $G$ to $w'_{0}$. If $w_{1}\in N_{H}(s')$, then since $\{yx,w_{0}w_{0}', w_{1}s', uu'\}$ are edges of $H$ and $\{xw_{0}, w'_{0}w_{1}, s'u\}$ are non-edges of $G$, Claim~\ref{scl:wandl} implies the contradiction $\Gamma_{1}'(s')\cup \Gamma_{1}'(x) \subseteq N_{G_{2}}(u')$. Thus, $w_{1}\notin N_{H}(s')$. Since $s'$ is adjacent in $G_{2}$ to $w_{0}$ and $s'\in S_{1}$, we may conclude that $s'$ is not adjacent in $G$ to $w_{1}$. Thus, there must be a $w'_{1}\in N_{H}(w_{1})\cap J$. Since $\{yx,ys',w_{0}w'_{0}, w_{1}w'_{1}\}$ are edges of $H$ and $\{xw_{0}, w'_{0}w_{1}, s'w_{1}\}$ are non-edges of $G$, Claim~\ref{scl:wandl} implies the contradiction $\Gamma_{1}'(s')\cup \Gamma_{1}'(x) \subseteq N_{G_{2}}(w_{1})$. This last contradiction shows that there is a vertex $w\in J$ adjacent in $G_{2}$ to every vertex in $\Gamma_{1}'(s')\cup \Gamma_{1}'(x)$.

However, since \[\Delta_{2}\geq |N_{G_{2}}(w)|\geq |\Gamma_{1}'(s')|+|\Gamma_{1}'(x)|\geq \min\{\Delta_{2}+1,2(\Delta_{2}-1)\},\] it must be the case that $\Delta_{2}=2$, and therefore, $|\Gamma_{1}'(s')|=|\Gamma_{1}'(x)|=1$. Thus, there must be an $l'\in L_{2}$ adjacent in $G_{2}$ to $x$. However, since $\Delta_{2}\geq \lambda$, we learn from (\ref{eq:lowerlambdafurther}) that $\delta_{+}=2$, and therefore, $l'$ is adjacent in $G_{2}$ to $s'$. Since $\Gamma_{1}'(x)\cap N_{G_{2}}(l')=\emptyset$, Claim~\ref{scl:basicObservation} says $x$ is adjacent in $H$ to $s'$.

If there is a vertex $u\in J$ that is not adjacent in $G$ to neither $x$ nor $s'$, then there must be distinct vertices $u'$ and $u''$ in $N_{H}(u)\cap J$. However, since $\{yx,uu',uu''\}$ are edges of $H$ and $xu$ is a non-edge of $G$,  Claim~\ref{scl:wandl} implies that every vertex in $\Gamma_{1}'(x)$ is adjacent in $G_{2}$ to the three vertices $\{x,u',u''\}$. This contradicts our claim that $\Delta_{2}=2$. Thus, there is some $x'\in \{x,s'\}$ that is adjacent in $H$ to at least $|J|/2$ vertices in $J$. Since $m=\Delta_{2}+1=3$, we know from (\ref{eq:V-L}) that $|J|=2\Delta_{1}-3$. However, since $x'$ is adjacent in $H$ to two vertices in $S\cup \{y\}$, we can deduce the contradiction \[|N_{H}(x')|\geq |J|/2+2\geq \frac{2\Delta_{1}-3}{2}+2> \Delta_{1}.\] This completes the proof of Claim~\ref{subsubcl:delta=1}.
\end{proof}  
\end{subsubclaim}

    We are left with the case $\delta_{+}=1$. If $\Delta_{1}=1$, then no vertex in $J$ is adjacent in $G$ to $x$. Thus, $x$ is not adjacent in $H$ to vertices in $N_{G_{2}}(z)-\{x\}$ for every $z\in \Gamma_{1}'(x)$. Therefore, by Claim~\ref{scl:basicObservation} $y$ must be adjacent in $G_{2}$ to the $\Delta_{2}$ vertices in $\Gamma_{1}'(x)\cup \{x\}$ and not adjacent in $G$ to vertices in $J$. Since $J$ is not empty, there is a $w\in J$, and therefore, there is a $w'\in N_{H}(w)\cap J$. However, by Claim~\ref{scl:wandl} there is a contradiction since $\{yx, ww'\}$ are edges of $H$ and $\{xw, w'y\}$ are non-edges of $G$. Therefore, $\Delta_{1}>\delta_{+}=1$.
 
    Let $C\subseteq N_{H}(x)\cap J$ be the set of vertices adjacent in $H$ to some vertex in $J-N_{H}(x)$. We define $D=(N_{H}(x)\cap J)-C$. For some $z\in \Gamma_{1}'(x)$, we let $R(z)$ be the vertices in $J-N_{H}(x)$ adjacent in $G$ to $z$ and let $F(z)=J-R(z)-N_{H}(x)$. By Claim~\ref{scl:wandl}, every vertex in $\Gamma'(x)$ is adjacent in $G_{2}$ to every vertex in $C$. Furthermore, if $u\in F(z)$ for some $z\in \Gamma'(x)$, then by Claim~\ref{scl:wandl} $N_{H}(u)\subseteq C$. We choose a $z\in \Gamma'_{1}(x)$, and let $R:=R(z)$ and $F:=F(z)$. Let $C'\subseteq C$ be the set of vertices adjacent in $H$ to vertices in $F(z)$.

    We now show that $F$ is not empty. Observe that
    \begin{align}
    |R\cup F|&\geq |V|-|L|-|N_{H}[x]|\notag\\
    &\geq (\Delta_{1}+1)(\Delta_{2}+1)-m-g-|L|-|N_{H}[x]|\notag\\
    &\geq (\Delta_{1}+1)(\Delta_{2}+1)-m-(\Delta_{2}-1)- (\Delta_{2}-1)-(\Delta_{1}+1)\notag\\
    &=(\Delta_{1}-1)\Delta_{2}+2-m.\notag
    \end{align}
    Since $\Delta_{1}\geq 2$, $|R\cup F|$ is not empty. Suppose that $F$ is empty.
    Since $R$ is not empty, $z$ is not adjacent in $G$ to some $l'\in \Gamma_{1}'(x)-\{z\}$ since $|\Gamma_{1}'(x)|=\Delta_{2}$, and therefore, Claim~\ref{scl:basicObservation} says $N_{G_{2}}(z)\cap (D\cup C)$ is not empty. Observe that $|R|\leq \Delta_{2}-1-|C|-|N_{G_{2}}(z)\cap D|$ since $C\subset N_{G_{2}}(z)$. However, this leads to the contradiction \begin{align}
        |F|&\geq (\Delta_{1}-1)\Delta_{2}+2-m-(\Delta_{2}-2-|C|-|N_{G_{2}}(z)\cap D|)\notag\\
        &=(\Delta_{1}-2)\Delta_{2}+3-m+|C|+|N_{G_{2}}(z)\cap D|\label{eq:FlowerBound}\\
        &\geq 1. \notag
    \end{align}
    Thus, $F$ is not empty, and therefore, $C'$ is not empty.
    
    
    

    For $w\in F$, if $N_{G_{2}}(w)\cap (C\cup D)=\emptyset$, then the realization of $\pi(G_{1})$ created by interchanging in $H$ $x$ with $w$ and then interchanging $y$ with $z$ in the resulting realization packs with $G_{2}$. This is a contradiction, and therefore, there is a $v\in N_{G_{2}}(w)\cap (C\cup D)$. Since every vertex in $C\cup \{x\}$ is adjacent in $G_{2}$ to the $\Delta_{2}$ vertices in $\Gamma_{1}'(x)$, it must be the case $v\in D$. Since $C'$ and $D$ are not empty, $\Delta_{1}\geq 3$.
    
    Suppose there is a $w'\in C'$ and $v'\in C'$ that are not adjacent in $G$. Let $w\in N_{H}(w')\cap F$ and $v\in N_{H}(v')\cap F$. Since $\{yx,vv', w'w\}$ are edges of $H$ and $\{xv, v'w', wy\}$ are non-edges of $G$, Claim~\ref{scl:wandl} implies the contradiction $w\in N_{G_{2}}(z)$. Therefore, $w'v'$ is an edge of $G$, and specifically, $w'v'$ must be an edge of $H$ since both $w'$ and $v'$ are adjacent in $G_{2}$ to the $\Delta_{2}$ vertices in $\Gamma_{1}'(x)$. Therefore, every vertex in $C'$ is adjacent in $H$ to $|C'|$ vertices in $C'\cup \{x\}$. 
    
    We now have $|F|\leq e_{H}(F,C')\leq |C'|(\Delta_{1}-|C'|)$. Combining (\ref{eq:FlowerBound}) with the upper bound of $e_{H}(F,C')$ we deduce that \[(\Delta_{1}-2)\Delta_{2}+3-m+|C'|\leq |F|\leq e_{H}(F,C')\leq |C'|(\Delta_{1}-|C'|).\] 
    
    Grouping all the $|C'|$ terms onto the right hand side and recognizing that $m\leq\Delta_{2}+1$ we deduce that \[(\Delta_{1}-3)\Delta_{2}+2\leq (\Delta_{1}-|C'|-1)|C'|.\] If $|C'|=1$, then since $\Delta_{1}\geq 3$, we deduce the contradiction \[\Delta_{1}-1\leq (\Delta_{1}-3)\Delta_{2}+2\leq \Delta_{1}-2.\] If $|C'|\geq 2$, then since $\Delta_{1}-|C'|-1\leq \Delta_{1}-3$ and $|C'|\leq \Delta_{2}-1$, we can deduce our final contradiction \[(\Delta_{1}-|C'|-1)|C'|+2<(\Delta_{1}-3)\Delta_{2}+2\leq (\Delta_{1}-|C'|-1)|C'|.\] This completes the proof of Claim~\ref{cl:lower}. \qedhere
\end{proof}
\end{subclaim}

Let $b'$ be the number of bad pairs incident with vertices in $S_{0}$.
\begin{subclaim}If $t=1$ and $S_{0}\neq \emptyset$, then
\begin{equation}\lambda\geq b-b'+|S_{2}|+\delta_{+}+(\Delta_{2}-1)|S_{0}|.\label{eq:lambdalarge}
\end{equation}
    \begin{proof}Since $t=1$, we know that $\Gamma_{1}(u)=\emptyset$ for all $u\in S_{0}$ and every vertex in $L_{2}$ is adjacent in $G_{2}$ to a vertex in $S-S_{0}$. This implies \[e_{G_{2}}(S,L_{2})= e_{G_{2}}(S-S_{0},L_{2})+ e_{G_{2}}(S_{0},L_{2})\geq |L_{2}|+|S_{0}|\Delta_{2}-b'.\] By (\ref{eq:lowerlambda}), we see that
\begin{align}
    \lambda&\geq b+|S_{2}|+\delta_{+}-|S_{0}|+e_{G_{2}}(S,L_{2})-|L_{2}|\notag\\
    &\geq b+|S_{2}|+\delta_{+}-|S_{0}|+|L_{2}|+\Delta_{2}|S_{0}|-b'-|L_{2}|\notag \\
    &\geq b-b'+|S_{2}|+\delta_{+}+(\Delta_{2}-1)|S_{0}|.\qedhere
\end{align}
   \end{proof} 
\end{subclaim}

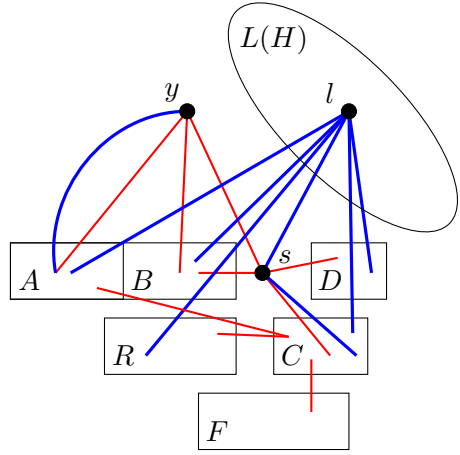
\begin{figure}[t]
    \centering
    \begin{tikzpicture}

\coordinate (s) at (3.35,2.35);
\coordinate (y) at (2.35,4.5);
\coordinate (l) at (4.5,4.5);
\coordinate (u) at (2.5,0);
\coordinate (v) at (4.5,.75);
\coordinate (uAdj) at (0,1.25);
\coordinate (vAdj) at (0,1.25);
\draw[] (u) rectangle (v); 
\draw[] (1.25,1) rectangle (3,1.75); 
\draw[] (0,2) rectangle (3,2.75); 
\draw[] (0,2) rectangle (1.5,2.75); 

\draw[] (3.5,1) rectangle (4.75,1.75); 
\draw[] (4,2.) rectangle (5,2.75);
 
\draw[color=red, thick] (y) -- (s);
\draw[color=red, thick] (y) -- (.6,2.35);
\draw[color=red, thick] (y) -- (2.25,2.35);

\draw[color=red, thick] (s) -- (4.35,2.55);
\draw[color=red, thick] (s) -- (2.5,2.35);
\draw[color=red, thick] (s) -- (4.25,1.25);
\draw[color=blue, very thick] (s) -- (4.6,1.25);
\draw[color=blue, very thick] (s) -- (l);

\draw[color=blue, very thick] (l) -- (4.55, 1.55);
\draw[color=blue, very thick] (l) -- (2.45, 2.5);
\draw[color=blue, very thick] (l) -- (4.8, 2.35);
\draw[color=blue, very thick] (l) -- (1.8, 1.25);
\draw[color=blue, very thick] (l) -- (.8, 2.35);

\draw[color=red, thick] (3.7,1.5) -- (2.75,1.54);
\draw[color=red, thick] (3.7,1.5) -- (1.15,2.15);
\draw[color=red, thick] (4,1.2) -- (4,.5);
\draw [blue, very thick]  (y) to[out=-180,in=100] (.6,2.35);
\draw[rotate=45] (6.25,0) ellipse (.8cm and 2cm);
  \draw[fill=black] (s) circle (.1cm);
  \draw[fill=black] (y) circle (.1cm);
\draw[fill=black] (l) circle (.1cm);
\node[] at (3.65,2.57)  {$s$};
\node[] at (2.15,4.75) {$y$};
\node[] at (4.25,4.75)  {$l$};

\node[] at (2.75,.25)  {$F$};
\node[] at (1.5,1.25)  {$R$};
\node[] at (.25,2.25)  {$A$};
\node[] at (1.75,2.25)  {$B$};
\node[] at (3.75,1.25)  {$C$};
\node[] at (4.25,2.25)  {$D$};
\node[] at (3.5,5.45)  {$L(H)$};

\end{tikzpicture}
    \caption{The edges of $H$ are colored red, and the edges of $G_{2}$ are a little thicker and colored blue. There is at least one bad pair between $A\cup B\cup \{s\}$ and $y$, and $s$ is adjacent in $G_{2}$ to at most one vertex in $B\cup C \cup D$. Moreover, $A\cup B\cup \{s\}=N_{H}(y)$, $B\cup C \cup D\subseteq N_{G}(s)$, $R\cup C\subseteq N_{G_{S}}(l)$, and $(A\cup R\cup F)\cap N_{H}(u)\neq \emptyset$ for every $u\in C$.}
    \label{fig:ABCSetup}
\end{figure}

Similar to the proof of Claim~\ref{cl:modify}, we distinguish some useful subsets in $V-L$ (See \ref{fig:ABCSetup}). First, we denote $B = S\cap N_{G}(s)$ and $A=S-B-\{s\}$. Next, we let $D\subseteq N_{G}(s)\cap (V-L-N_{H}[y])$ such that $N_{H}(u)\subseteq N_{G}[s]\cap (V-L)$ for all $u\in D$. With $D$ we can define $C=(N_{G}(s)\cap (V-L))-D-N_{H}[y]$. Finally, we let $R\subseteq V-L-(N_{G}[s]\cup A)$ be the set of vertices adjacent in $G_{2}$ to $l$, and let $F=V-L-(N_{G}[s]\cup A \cup R)$. 



Subtracting $|N_{G}[s]\cap (V-L)|=deg_{H}(s)+t'+1$ from (\ref{eq:V-L}) we see that
\begin{align}|A\cup R\cup F|&=|V-L|-|N_{G}[s]\cap(V-L)|\notag\\
&\geq \Delta_{2}(\Delta_{1}-\delta_{+})-g+\Delta_{1}+\Delta_{2}+1-m+\lambda-(deg_{H}(s)+t'+1)\notag\\
&=\Delta_{2}(\Delta_{1}-\delta_{+})-g+\Delta_{1}-deg_{H}(s)+\Delta_{2}+1-m+\lambda-1-t'.\label{eq:ARFFull}
\end{align}

We will also need the related bound 
\begin{equation}\label{eq:theRestofH}
    |R\cup C\cup D\cup F|\geq |V-L|-|N_{H}[y]|\geq (\Delta_{2}+1)(\Delta_{1}-\delta_{+})+\Delta_{2}-m-g+\lambda.
\end{equation}

If there is an edge $u_{2}v_{2}$ of $H$ such that $u_{2}\in A\cup R\cup F$ and $v_{2}\in C\cup R \cup F$, then the edges $\{ys, u_{2}v_{2}\}$ of $H$ and the non-edge $su_{2}$ of $G$ would satisfy Claim~\ref{cl:modify} if $v_{2}$ satisfies \ref{claim2con3}. Thus, $l$ is adjacent in $G_{2}$ to $v_{2}$, and $v_{2}$ must be adjacent in $G_{2}$ to at least  $\Delta_{2}-t'$ vertices in $L\cup \{y\}$. Therefore, $v_{2}\notin F$. This implies $N_{H}(u)\subseteq B\cup C$ for every $u\in F$, and $l$ must be adjacent in $G_{2}$ to every vertex in $R\cup C$. Thus, $|R\cup C|\leq \Delta_{2}-d$ where $d=|N_{G}(l)\cap S|$. Moreover, since every vertex in $C$ is adjacent in $H$ to a vertex in $A\cup R\cup F$, every vertex in $C$ is adjacent in $G_{2}$ to at least $\Delta_{2}-t'$ vertices in $L\cup \{y\}$.

\begin{subclaim}\label{subcl:dcupf}$D\cup F$ is not empty.
\begin{proof}Suppose $D\cup F$ is empty.  Utilizing (\ref{eq:theRestofH}) we may bound $|R\cup C|$ such that \[\Delta_{2}-d\geq|R\cup C|\geq (\Delta_{2}+1)(\Delta_{1}-\delta_{+})+\Delta_{2}-m-g+\lambda.\]
Rearranging, we construct
\begin{equation}\label{eq:RClower}
0\geq (\Delta_{2}+1)(\Delta_{1}-\delta_{+})-m-g+\lambda+d.
\end{equation}
If $\Delta_{1}=\delta_{+}$, then since $\lambda\geq\Delta_{1}+1\geq m$ by Claim~\ref{cl:lower}, $g=0$, and $d\geq 1$, there is a contradiction since the right hand side of (\ref{eq:RClower}) is at least one. Suppose $\Delta_{1}\neq \delta_{1}$, and thus, $g=\Delta_{2}-1$. If $\lambda\geq \Delta_{2}+1$, then $m=\Delta_{2}+1$, and therefore, there is a contradiction since \[0\geq (\Delta_{2}+1)(\Delta_{1}-\delta_{+})-(\Delta_{2}-1)+d\geq 3.\] If $\Delta_{2}+1>\lambda\geq\delta_{+}+1$, then $\Delta_{1}+1\geq m$, and we can deduce our final contradiction
\begin{align}
    0&\geq (\Delta_{2}+1)(\Delta_{1}-\delta_{+})-(\Delta_{1}+1)-(\Delta_{2}-1)+\delta_{+}+1+d\notag\\
    &= \Delta_{2}(\Delta_{1}-\delta_{+}-1)+1+d\geq 1.\notag \qedhere
\end{align}
\end{proof}
\end{subclaim}

\begin{subclaim}\label{subclaim:neighborhood}If there is a $v_{3}\in D$ not adjacent in $G_{2}$ to $l$, then $N_{H}(u)\subseteq B\cup C\cup \{y\}$ for all $u\in A\cup R\cup F$.
\begin{proof}If there is a $u_{2} \in A\cup R\cup F$ adjacent in $H$ to some $v_{2}\in A\cup R\cup F$, then the edges $\{ys,u_{2}v_{2},sv_{3}\}$ would satisfy Claim~\ref{cl:modify}.
\end{proof}
\end{subclaim}

    

\begin{subclaim}\label{subcl:AandRBound}$\Delta_{2}-d\geq |A|+|R|$.
\begin{proof}Since $\Delta_{2}-d\geq |R\cup C|$ we are done if $|C|\geq |A|$. Suppose $|A|>|C|$. Since $l$ is adjacent in $G_{2}$ to every vertex in $R\cup C$ and $d$ vertices in $S$, it is sufficient to show that $|N_{G_{2}}(l)\cap (C\cup D)|\geq |A|$.  If $l$ is adjacent in $G_{2}$ to every vertex in $D$, then \[|N_{G_{2}}(l)\cap (C\cup D)|=|N_{H}(s)\cap (C\cup D)|= |N_{H}(s)-B-\{y\}|\geq \delta_{+}-|B|-1=|A|.\] Suppose there is a $w\in D$ not adjacent in $G_{2}$ to $l$. Let $u\in A$. By Claim~\ref{subclaim:neighborhood}, we see that $N_{H}(u)-\{y\}\subseteq B\cup C$. Thus, $|C|\geq \delta_{+}-|B|-|\{y\}|=|A|$.
\end{proof}
\end{subclaim}

\begin{subclaim}\label{cl:FnotEmpty}$F$ is not empty.
\begin{proof}In search of a contradiction, we assume $F$ is empty.

\begin{subsubclaim}\label{cl:stats}The following is true:
\begin{itemize}
\item $t'=1$.
\item $\Delta_{1}\leq \delta_{+}+1$.
\item $|A\cup R|=\Delta_{2}-1$.
    \item $\lambda = m$ when $\Delta_{1}=\delta_{+}$.
    \item $\lambda = \Delta_{1}\leq \Delta_{2}$ when $\Delta_{1}=\delta_{+}+1$.
\end{itemize}
\begin{proof}By (\ref{eq:ARFFull}) and Claim~\ref{subcl:AandRBound}, we may construct \[\Delta_{2}-1\geq |R\cup A|\geq \Delta_{2}(\Delta_{1}-\delta_{+})+\Delta_{2}-m-g+\lambda-t'.\]

Suppose $\Delta_{1}\geq \delta_{+}+2$. Since $g\leq\Delta_{2}-1$, $m\leq \Delta_{2}+1$ and $t'\leq 1$ we may establish the contradiction \[|A\cup R|\geq 2\Delta_{2}+\Delta_{2}-(\Delta_{2}+1)-(\Delta_{2}-1)+\lambda-1\geq \Delta_{2}-1+ \lambda\geq \Delta_{2}.\] Thus, $\Delta_{1}\leq \delta_{+}+1$.

If $\Delta_{1}=\delta_{+}$, then $g=0$. Thus, $|A\cup R|\geq \Delta_{2}-t'+\lambda-m.$ This can only happen when $t'=1$, $\lambda=m$, and $|A\cup R|=\Delta_{2}-1$. If $\Delta_{1}=\delta_{+}+1$, then using $g=\Delta_{2}-1$ we see that \[|A\cup R|\geq \Delta_{2}+\Delta_{2}-m-(\Delta_{2}-1)+\lambda-t'\geq \Delta_{2}-t'+\lambda+1-m.\] If $\Delta_{2}<\Delta_{1}$, then $\lambda=m=\Delta_{2}+1$ and we may establish the contradiction $\Delta_{2}-1\geq |A\cup R|\geq \Delta_{2}+1-t'$. Thus,  $\Delta_{2}\geq\Delta_{1}=\lambda$, and therefore, $t'=1$ and $|A\cup R|=\Delta_{2}-1$.
\end{proof}
\end{subsubclaim}

\begin{subsubclaim}$l\neq y$.
\begin{proof}Suppose by contradiction that $y=l$. By our assumption on $y$, we know that $\Gamma(s)\subseteq L_{2}$ and every vertex in $\Gamma(s)$ is adjacent in $G_{2}$ to a vertex in $S_{2}$. Since $\Delta_{2}\geq 2$, $t=1$, and $s$ is adjacent in $G_{2}$ to $y$, we deduce that $|L_{2}|\geq \Delta_{2}-2$ and $S_{2}$ is not empty. Since $s\in S_{1}$, $|S_{0}|\leq \delta_{+}-2$. From (\ref{eq:lowerlambdafurther}), we deduce that \[\lambda\geq b+|S_{2}|+\delta_{+}-|S_{0}|+|L_{2}|\geq b+1+\delta_{+}-(\delta_{+}-2)+\Delta_{2}-2\geq \Delta_{2}+2.\] However, this contradicts Claim~\ref{cl:stats} since $\Delta_{2}+1\geq m\geq \lambda\geq \Delta_{2}+2$. Thus, $l\neq y$.
\end{proof}
\end{subsubclaim}

\begin{subsubclaim}\label{cl:lnotAdjacent} $|D-N_{G_{2}}(l)|\geq 2$ when $\Delta_{1}=\delta_{+}+1$ and $|D-N_{G_{2}}(l)|\geq 1$ otherwise.
\begin{proof}Let $i$ be the number of vertices in $D$ not adjacent to $l$. This implies $\Delta_{2}-1\geq |R\cup C\cup D|-i$ since $l$ is adjacent in $G_{2}$ to $s$ and not adjacent to $i$ vertices in $D$. On the other hand, with (\ref{eq:theRestofH}) we can deduce that \begin{equation}\label{eq:lnotAdjacent}
    \Delta_{2}-1\geq |R\cup C\cup D|-i\geq (\Delta_{2}+1)(\Delta_{1}-\delta_{+})+\Delta_{2}-m-g+\lambda-i.
\end{equation} Observe that when $\Delta_{1}=\delta_{+}$, the right hand side of (\ref{eq:lnotAdjacent}) is at least $\Delta_{2}-i$. Which is a contradiction when $i=0$. When $\Delta_{1}=\delta_{+}+1$, we know that $g=\Delta_{2}-1$ and $\lambda=m-1$. Thus, when $i\leq 1$ we observe the contradiction \[\Delta_{2}-1\geq |R\cup C\cup D|-i\geq \Delta_{2}+1+\Delta_{2}-m-(\Delta_{2}-1)+m-1\geq \Delta_{2}+1-i. \qedhere\]
\end{proof}
\end{subsubclaim}

Since $t'=1$, we know that $|N_{G}(s)\cap (V-L)|\leq \Delta_{1}+1$. This implies \[|B\cup C|=|N_{G}(s)\cap (V-L)-D-\{y\}|\leq \Delta_{1}+1-|\{y\}|-|D|=\Delta_{1}-|D|.\] Utilizing Claim~\ref{cl:stats} and Claim~\ref{cl:lnotAdjacent} we conclude that $|B\cup C|\leq \delta_{+}-1$. Since $l$ is not adjacent in $G$ to some vertex in $D$, we know by Claim~\ref{subclaim:neighborhood} that $N_{H}(u)\subseteq B\cup C\cup \{y\}$ for all $u\in A\cup R$. Since $N_{H}(u)\subseteq B\cup C$ for every $u\in R$ and $|B\cup C|\leq \delta_{+}-1$, we may conclude that $R$ is empty.  Therefore, $|A|=\Delta_{2}-1$, $|B\cup C|=\delta_{+}-1$, and
\[\Delta_{2}-1=|A|=\delta_{+}-1-|B|=|C|\leq |N_{G}(s)\cap (V-L)|\leq \Delta_{1}+1-|D|-|B|-1.\] Since $D$ is not empty, the last inequality yields a contradiction if $\Delta_{2}>\Delta_{1}$. Moreover, if  $\Delta_{1}=\delta_{+}+1$, then we can deduce a contradiction since Claim~\ref{cl:stats} implies $\Delta_{1}\leq \Delta_{2}$ and Claim~\ref{cl:lnotAdjacent} implies $|D|\geq 2$. Thus, $2\leq \Delta_{2}\leq \Delta_{1}=\delta_{+}$. By Claim~\ref{cl:stats}, (\ref{eq:lowerlambdafurther}), and (\ref{eq:lambdalarge}), we deduce that $\delta_{+}+1+|S_{2}|\leq \lambda=m=\Delta_{2}+1\leq \delta_{+}+1$. Therefore, $\delta_{+}=\Delta_{1}=\Delta_{2}$, and thus, $|D|=1$ by Claim~\ref{cl:lnotAdjacent}. This implies $B$ is empty since $\Delta_{2}-1\leq \Delta_{1}-|D|-|B|\leq \Delta_{1}-1-|B|$.

By Claim~\ref{cl:lnotAdjacent}, $D$ has a single vertex $u_{3}$ not adjacent to $l$. Since $B$ is empty and $|C\cup \{s\}|=\delta_{+}$, we learn from the definition of $D$ that $C\cup \{s\}= N_{H}(u_{3})$. Let $u_{4}\in C$ and $u_{5}\in A$. Recall that $C\subseteq N_{H}(u_{5})$. If there is a $v\in C$ not adjacent in $G$ to $u_{4}$, then the edges $\{ys, u_{5}u_{4}, vu_{3}\}$ in $H$ and the non-edges $\{su_{5},u_{4}v,u_{3}y\}$ of $G$ would satisfy Claim~\ref{cl:modify}. Thus, $u_{4}$ must be adjacent in $G$ to every vertex in $C-\{u_{4}\}\cup A\cup \{u_{3},s\}$. Thus, $|N_{G}(u_{4})\cap (V-L-\{y\})|\geq 2(\Delta_{2}-1)+1=2\Delta_{2}-1$. If $u_{4}$ is adjacent in $G_{2}$ to at most $\Delta_{2}-2$ vertices in $V-L$, then the edges $\{ys,u_{5}u_{4}\}$ of $H$ and the non-edge $su_{5}$ of $G$ would satisfy Claim~\ref{cl:modify}. Thus, $u_{4}$ is adjacent to at least $\Delta_{2}-1$ vertices in $L\cup \{y\}$. This implies $|N_{G}(u_{4})|\geq 2\Delta_{2}-1+\Delta_{2}-1=3\Delta_{2}-2$. Thus, $\Delta_{2}=2$ since $|N_{G}(u_{4})|\leq \Delta_{2}+\Delta_{1}\leq 2\Delta_{2}$ and $\Delta_{1}=\Delta_{2}$. Therefore, $A=\{u_{5}\}$, $C=\{u_{4}\}$, and $D=\{u_{3}\}$. We can conclude that $su_{4}$ is an edge of $G_{2}$ since $t'=1$ and $s$ is adjacent in $H$ to $y$ and $u_{3}$. Since $\lambda=3$, we may conclude that $L=\{l\}$. Therefore, $V=\{y,l,s,u_{3},u_{4},u_{5}\}$. Moreover, since $N_{G_{2}}(l)=\{s,u_{4}\}$ and $N_{G_{2}}(s)=\{l,u_{4}\}$, it must be the case that $u_{5}\in N_{G_{2}}(y)$. Therefore, by Claim~\ref{scl:basicObservation}, $u_{3}\in N_{G_{2}}(y)$. If $u_{3}u_{5}$ is an edge of $G_{2}$, then $G_{2}$ is the disjoint union of two three cycles $lsu_{4}$ and $yu_{3}u_{5}$. However, in this case $H$ and $G_{2}$ violate \ref{main:fivecycle} since $H$ is the disjoint union of the five cycle $ysu_{3}u_{4}u_{5}y$ and the isolate $l$. Therefore, $u_{3}u_{5}$ must be a non-edge of $G$. However, in this case, we may interchange $u_{3}$ with $y$ to create a realization of $\pi(G_{1})$ that packs with $G_{2}$. Therefore, $F$ must not be empty. This completes the proof of Claim~\ref{cl:FnotEmpty}.
\end{proof}
\end{subclaim}
Let $B'$ be the set of vertices in $B$ adjacent in $H$ to vertices in $A\cup R\cup F$, let $C'\subseteq C$ be the vertices adjacent in $H$ to vertices in $F$, and $D'\subseteq D$ be the set of vertices not adjacent in $G_{2}$ to $l$.

\begin{subclaim}\label{cl:wNieghbors}Let $w\in F$. For $u_{2}\in A\cup R\cup F$, every $u_{3}\in N_{H}(w)$ is adjacent in $G$ to every vertex in $N_{H}(u_{2})-\{y\}$. Moreover, $u_{3}$ is adjacent in $G$ to every vertex in $(B'\cup C)-\{u_{3}\}$.
\begin{proof}If there is a $v_{2}\in N_{H}(u_{2})-\{y\}$ not adjacent in $G$ to some $u_{3}\in N_{H}(w)$, then since $w$ satisfies \ref{claim2con3}, the edges $\{ys,u_{2}v_{2},u_{3}w\}$ of $H$ and the non-edges $\{su_{2},v_{2}u_{3},wy\}$ of $G$ would satisfy Claim~\ref{cl:modify}. Since every vertex in $B'\cup C$ is adjacent in $H$ to a vertex in $A\cup R\cup F$, it easily follows that every $u_{3}\in N_{H}(w)$ is adjacent to every vertex in $(B'\cup C)-\{u_{3}\}$.
\end{proof}
\end{subclaim}

Since $F$ is not empty, we now know that every $u\in C'$ is adjacent in $G$ to every vertex in $B'\cup (C-\{u\})$.

\begin{subclaim}\label{cl:neighofD'}Let $v_{3}\in D'$. If $v_{2}\in N_{H}(u_{2})$ for some $u_{2}\in F\cup R$, then $v_{2}$ is adjacent to every vertex in $N_{H}(v_{3})$.
\begin{proof}If $v_{2}$ is not adjacent to some $u_{3}\in N_{H}(v_{3})$, then the edges $\{ys, u_{2}v_{2}, u_{3}v_{3}\}$ of $H$ and the non-edges $\{su_{2}, v_{2}u_{3}, v_{3}y\}$ of $G$ satisfy Claim~\ref{cl:modify} since $v_{3}$ satisfies \ref{claim2con3}.
\end{proof}
\end{subclaim}

\begin{subclaim}\label{cl:delta_+special}If $t'=1$ and $\delta_{+}=1$, then $l\neq y$ and $\Delta_{2}\geq 3$. 
\begin{proof}If $\delta_{+}=1$, then $s$ and $y$ form a bad pair, and $S_{2}$ is empty. Therefore, $l\neq y$.  Since $t'=1$, we know that, in addition to $y$ and $l$, $s$ is adjacent in $G_{2}$ to some vertex not in $L\cup \{y\}$. Thus, $\Delta_{2}\geq 3$.
\end{proof}
\end{subclaim}

\begin{subclaim}\label{cl:D-dbound}Every vertex in $C$ is adjacent to at most $\Delta_{1}-\delta_{+}$ vertices in $A\cup R\cup F$.
\begin{proof}
Let $w\in F$. If there is a $u\in (B'\cup C)-N_{H}(w)$, then by Claim~\ref{cl:wNieghbors} every vertex in $C$ is adjacent to at least $\delta_{+}+1$ vertices in $N_{H}(w)\cup \{u,s\}$. Thus, every vertex in $C$ is adjacent to at most $\Delta_{1}+t'-(\delta_{+}+1)\leq \Delta_{1}-\delta_{+}$ vertices in $A\cup R\cup F$. Thus, we are left with the case that $N_{H}(w)=B'\cup C$ for every $w\in F$, and therefore, $C'=C$.

We choose a $v\in C'$ that is adjacent in $G$ to the most vertices in $A\cup R\cup F$. We suppose, by contradiction, that $v$ is adjacent in $G$ to at least $\Delta_{1}-\delta_{+}+1$ vertices in $A\cup R\cup F$. We know by Claim~\ref{cl:wNieghbors} that $v$ is adjacent to every vertex in $(B'\cup C-\{v\})\cup\{s\}$.  Thus, $v$ is adjacent to at most $\Delta_{1}+t'-|(B'\cup C-\{v\})\cup\{s\}|$ vertices in $A\cup R\cup F$. Since $|B'\cup C|\geq \delta_{+}$, it must be the case that $|B'\cup C|=\delta_{+}$ and $t=t'=1$. Thus, $v$ is adjacent in $G$ to exactly $\Delta_{1}-\delta_{+}+1$ vertices in $A\cup R\cup F$.

\begin{subsubclaim}\label{cl:minIntersection}Every vertex in $R\cup F$ is adjacent in $G$ to every vertex in $B'\cup C$, and every vertex in $A$ is adjacent in $G$ to at least $\delta_{+}-1$ vertices in $B'\cup C$.
\begin{proof}Let $u_{2}\in A\cup R\cup F$. If there exists a $v_{2}\in N_{H}(u_{2})\cap (A\cup R\cup F)$, then since $B'\cup C=N_{H}(w)$ for every $w\in F$, we know by Claim~\ref{cl:wNieghbors} that $u_{2}$ is adjacent in $G$ to every vertex in $B'\cup C$. Suppose $N_{H}(u_{2})\subseteq B'\cup C$. If $u_{2}\in A$, then since $u_{2}$ is adjacent in $H$ to $y$ we know that $|N_{H}(u_{2})\cap (B'\cup C)|\geq \delta_{+}-1$. If $u_{2}\in R\cup F$, then $N_{H}(u_{2})=B'\cup C$ since $|B'\cup C|=\delta_{+}$.  
\end{proof}
\end{subsubclaim}

\begin{subsubclaim}\label{cl:adjARF}$v$ is adjacent in $G$ to every vertex in $A\cup R\cup F$.
\begin{proof}By Claim~\ref{cl:minIntersection}, every vertex in $R\cup F$ is adjacent to every vertex in $B'\cup C$. If every vertex in $C'$ is not adjacent to every vertex in $A$, then since $|C|=\delta_{+}-|B'|\geq \delta_{+}-|B|=|A|+1$, we may conclude that some vertex $z\in A$ is not adjacent in $G$ to two vertices in $C$. However, this would contradict Claim~\ref{cl:minIntersection} since $|B'\cup C|=\delta_{+}$. Therefore, $v$ must be adjacent to every vertex in $A\cup R\cup F$.
\end{proof}
\end{subsubclaim}

We deduce from Claim~\ref{cl:adjARF} that $|A\cup R\cup F|=|N_{G}(v)\cap (A\cup R\cup F)|=\Delta_{1}-\delta_{+}+1$.

\begin{subsubclaim}\label{cl:delta2boundF}$\Delta_{2}\geq3$
\begin{proof}Suppose $\Delta_{2}=2$. This implies $|R\cup C|=1$ since $|N_{G_{2}}(l)\cap (R\cup C)|\leq \Delta_{2}-1$. Let $w\in F$. Since $N_{H}(w)= B'\cup C$ we see that $R$ is empty, $|C|=1$, and $|B'|=\delta_{+}-1$. Therefore, $B'=B$ and $A$ is empty. Since $\Delta_{2}=2$, $t'=1$, $L\neq \emptyset$, and $y$ forms a bad pair with some vertex in $S$, we may conclude that there is an $s'\in S-\{s\}$ adjacent in $G_{2}$ to some vertex in $L\cup \{y\}$. This implies $s'\notin S_{2}$ and $|N_{G}(s')\cap (V-L)|\leq \Delta_{1}+1$. From  Claim~\ref{cl:wNieghbors}, we see $s'$ is adjacent in $G$ to every vertex in $(B'-\{s'\})\cup C$. However, this gives us the contradiction \[\Delta_{1}+1\geq |N_{G}(s')\cap (V-L)|\geq |B'\cup C|-|\{s'\}|+|\{s,y\}|+|F|\geq \delta_{+}+1+\Delta_{1}-\delta_{+}+1\geq\Delta_{1}+2.\qedhere\]
\end{proof}
\end{subsubclaim}

\begin{subsubclaim}$\Delta_{1}>\delta_{+}$ 
\begin{proof}If $\Delta_{1}=\delta_{+}$, then $\lambda\geq m$ by Claim~\ref{cl:lower}. Combining $|A\cup R\cup F|=1$, $\Delta_{2}\geq 3$, and (\ref{eq:ARFFull}) with $g=0$ and $t'=1$ we establish the contradiction $1=|A\cup R\cup F|\geq \Delta_{1}-deg_{H}(s)+\Delta_{2}-t'\geq 2$.\qedhere
\end{proof}
\end{subsubclaim}

\begin{subsubclaim}\label{cl:Ddld1}$\Delta_{2}-3+\Delta_{1}-deg_{H}(s)+\lambda\leq \delta_{+}+1$.
\begin{proof}Suppose $\Delta_{2}-3+\Delta_{1}-deg_{H}(s)+\lambda\geq \delta_{+}+2$. Recall that $t'=1$ and $|A\cup R\cup F|=\Delta_{1}-\delta_{+}+1$. By (\ref{eq:ARFFull}), we see that \[\Delta_{1}-\delta_{+}+1= |A\cup R\cup F|\geq \Delta_{2}(\Delta_{1}-\delta_{+})-g+\Delta_{1}-deg_{H}(s)+\Delta_{2}-m+\lambda-1.\] Grouping all terms onto the right hand side and simplifying we see \[0\geq (\Delta_{2}-1)(\Delta_{1}-\delta_{+})-g+\Delta_{1}-deg_{H}(s)+\lambda-2+\Delta_{2}-m.\] Inserting our assumed lower bound $\Delta_{1}-deg_{H}(s)+\lambda\geq \delta_{+}+2-\Delta_{2}+3$ we see that \[0\geq (\Delta_{2}-1)(\Delta_{1}-\delta_{+})-g+\delta_{+}+3-m.\] If $\delta_{+}+1=\Delta_{1}$, then we deduce the contradiction \[0\geq \delta_{+}+3-m=\Delta_{1}+1-m+1\geq 1.\] If $\delta_{+}+2\leq \Delta_{1}$, then we can deduce the contradiction \[0\geq \Delta_{2}-1+\delta_{+}+3-m\geq \Delta_{2}+1-m+\delta_{+}+1\geq 1.\] Thus, $\Delta_{2}-3+\Delta_{1}-deg_{H}(s)+\lambda\leq \delta_{+}+1$.
\end{proof}
\end{subsubclaim}
Since $\Delta_{2}\geq 3$, Claim~\ref{cl:Ddld1} implies $\lambda\leq \delta_{+}+1$. Since $t=t'=1$, $\Delta_{2}\geq 3$, and $\lambda\leq \delta_{+}+1$, we know by (\ref{eq:lambdalarge}) that $S_{0}=\emptyset$. Therefore, by (\ref{eq:lowerlambdafurther}) we deduce that both $S_{2}$ and $L_{2}$ are empty. Thus, $S=S_{1}$, and therefore, $\lambda=\delta_{+}+1$. Moreover, Claim~\ref{cl:Ddld1} implies $\Delta_{2}=3$ and $\Delta_{1}=deg_{H}(s)$. Therefore, $|\Gamma_{1}'(u)|=2$ for every $u\in S$. Furthermore, $y\neq l$ since $S_{2}$ is empty and $\Delta_{2}=3$.

Since $\Delta_{1}>\delta_{+}$, we know that $|A\cup R\cup F|\geq 2$. Let $w\in F$. Suppose there is a $u\in A$. In this case $|B|\leq |S|-|\{s,u\}|\leq \delta_{+}-2$. Thus, $|C|\geq \delta_{+}-|B'|\geq \delta_{+}-|B|\geq 2$. Therefore, there is a $u'\in C\cap N_{H}(s)$. Since $\{yu,su'\}$ are edges of $H$ and $us$ is a non-edge of $G$, we know by Claim~\ref{scl:wandl} that $u'$ is adjacent in $G_{2}$ to every vertex in $\Gamma_{1}'(u)$. However, this implies $u'$ is adjacent in $G_{2}$ to the four vertices in $\Gamma_{1}'(u)\cup \Gamma_{1}'(s)$ since $u'\in C$. This contradicts $\Delta_{2}=3$. Therefore, $A=\emptyset$. 

Suppose there is a $u\in B-B'$. Since $|R\cup F|\geq 2$ and $u\in S_{1}$ there is an $x\in R\cup F$ not adjacent in $G$ to $u$. Let $x'\in N_{H}(x)\cap (R\cup C)$. Since $\{yu,ys,xx'\}$ are edges of $H$ and $\{sx,ux\}$ are non-edges of $G$, we see by Claim~\ref{scl:wandl} that $x$ is adjacent in $G_{2}$ to the four vertices in $\Gamma_{1}'(u)\cup \Gamma_{1}'(s)$. However, this is a contradiction since $\Delta_{2}=3$. Thus, $B=B'$. By Claim~\ref{cl:minIntersection}, every vertex in $u\in B'$ is adjacent in $G$ to every vertex in $R\cup F\cup B'\cup C \cup \{y,s\}-\{u\}$. However, this leads to the following contradiction \[\Delta_{1}+1\geq |R\cup F\cup B'\cup C \cup \{y,s\}-\{u\}|\geq \Delta_{1}-\delta_{+}+1+\delta_{+}-1+2=\Delta_{1}+2.\] Thus, $B$ is empty, and therefore, $\delta_{+}=1$ and $C=\{v\}$. 

Since $\Delta_{2}=3$ and $\{y,l\}\subseteq N_{G_{2}}(v)$, we deduce that $|R|\leq 1$, and therefore, $v$ is adjacent in $H$ to every vertex in $R\cup F$. Since $|A\cup R\cup F|=\Delta_{1}$, $v$ is adjacent in $G_{2}$  to $s$, $l$, and $y$. Furthermore, $v$ is not adjacent in $G$ to vertices in $D$.

We are left with the situation that $v$ is adjacent in $G$ to $\Delta_{1}$ vertices in $R\cup F$ and every vertex in $V-(R\cup F)-\{s\}$ is adjacent in $G$ to $s$. This implies $n\leq |N_{G}(v)\cap (R\cup F)|+|N_{G}[s]|\leq \Delta_{1}+\Delta_{1}+\Delta_{2}\leq 2\Delta_{1}+3$. On the other hand, (\ref{eq:main}) implies $n\geq 4(\Delta_{1}+1)-\Delta_{1}-1-2\geq 3\Delta_{1}+1$. This can only happen when $\Delta_{1}=2$.

Since $deg_{H}(s)=\Delta_{1}=2$ and $\{l,y,v\}=N_{G_{2}}(s)$, $s$ is adjacent in $H$ to exactly one vertex in $D$. Let $D=\{z\}$.

If $y$ is not adjacent in $G_{2}$ to $l$, then by Claim~\ref{scl:basicObservation}, both $y$ and $l$ must be adjacent in $G_{2}$ to $z$. Therefore, $R=\emptyset$ and $|F|\geq 2$. This implies there must be a $w\in F$ not adjacent in $G_{2}$ to $z$ nor $y$ and not adjacent in $H$ to  $l$ nor $y$. In this case, we may exchange in $H$ the edges $\{wv, zs\}$ with the non-edges $\{sw,vz\}$ to construct a realization $H'$ of $\pi(G_{1})$. We then in  $H'$ interchange $s$ with $w$ to construct a realization of $\pi(G_{1})$ that packs with $G_{2}$. Thus, it must be the case that every vertex in $C\cup B'$ is adjacent in $G$ to at most $\Delta_{1}-\delta_{+}$ vertices in $A\cup R\cup F$. This completes the proof of Claim~\ref{cl:D-dbound}.\qedhere
\end{proof}
\end{subclaim}

Since $F$ is not empty, we know that $e_{H}(C,A\cup R\cup F)\geq 1$. By Claim~\ref{cl:D-dbound}, \[|C|(\Delta_{1}-\delta_{+})\geq e_{H}(C,A\cup R\cup F).\] Clearly, $\Delta_{1}=\delta_{+}$ gives us a contradiction; hence, we conclude that $\Delta_{1}>\delta_{+}$. 

\begin{subclaim}\label{cl:ectoF} Every vertex in $A\cup R\cup F$ is adjacent in $H$ to a vertex in $C$. 
\begin{proof}
Let $Q\subseteq A\cup R\cup F$ be the largest set of vertices that are not adjacent in $H$ to the vertices in $C$. Let $w\in F$ and $v\in N_{H}(w)$. Since every vertex in $Q$ is adjacent in $H$ to some vertex in $A\cup R\cup F$, we know, by Claim~\ref{cl:wNieghbors} and our assumption on $Q$, that $v$ is adjacent in $G_{2}$ to every vertex $Q$. However, if $|Q|\geq 2$, then $v$ would be adjacent in $G_{2}$ to at most $\Delta_{2}-2$ vertices in $L\cup \{y\}$. However, this is a contradiction since $\{ys, wv\}$ would satisfy Claim~\ref{cl:modify}. Thus, $Q$ has a single vertex $u$, and therefore, there is a $u'\in N_{H}(u)\cap (A\cup R\cup F)$ not in $Q$. If there is an $x\in (A\cup R\cup F)-\{u\}$ not adjacent in $H$ to $v$, then there must be an $x'\in N_{H}(x)\cap C$. Since $\{ys, u'u$\} are edges of $H$ and $su'$ is a non-edge of $G$, we see by Claim~\ref{scl:wandl} that $u$ is adjacent in $G_{2}$ to the $\Delta_{2}-1$ vertices in $\Gamma_{1}'(s)$. Since $u$ is also adjacent in $G_{2}$ to $v$, it must not be adjacent in $G$ to $x'$. Therefore, since $\{ys, u'u, x'x, vw\}$ are edges of $H$ and $\{su',ux',xv\}$ are non-edges of $G$, Claim~\ref{scl:wandl} implies the contradiction that $l$ is adjacent in $G_{2}$ to $w$. Therefore, $v$ must be adjacent in $G$ to every vertex in $A\cup R\cup F$. Thus,  $\Delta_{1}-\delta_{+}\geq |A\cup R\cup F|$ by Claim~\ref{cl:D-dbound}. By Claim~\ref{cl:wNieghbors}, $v$ is adjacent in $G$ to the vertices in $B'\cup (C-\{v\})\cup\{s\}$. In particular, $v$ is adjacent in $H$ to $\delta_{+}$ vertices in $B'\cup (C-\{v\})\cup\{s\}$ since $v$ is adjacent in $G_{2}$ to $u$, $t'\leq 1$, and $|B'\cup (C-\{v\})\cup\{s\}|\geq \delta_{+}$. Since in addition $v$ is adjacent in $H$ to $w$ and $u'$, we see that $\Delta_{1}\geq \delta_{+}+2$. Since $t'=1$ and $g=\Delta_{2}-1$, (\ref{eq:ARFFull}) implies that \[\Delta_{1}-\delta_{+}\geq \Delta_{2}(\Delta_{1}-\delta_{+})-(\Delta_{2}-1)+\Delta_{1}-deg_{H}(s)+\Delta_{2}-m+\lambda-1.\] Combining all the terms onto the right hand side and simplifying we establish \[0\geq (\Delta_{2}-1)(\Delta_{1}-\delta_{+}-1)+\Delta_{1}-deg_{H}(s)+\Delta_{2}+\lambda-1-m.\] Since $\Delta_{2}\geq 2$, $\Delta_{1}\geq \delta_{+}+2$, and $\lambda\geq 2$, the right hand side is positive. With this contradiction we showed that every vertex in $A\cup R\cup F$ is adjacent in $H$ to a vertex in $C$. This ends the proof of Claim\ref{cl:ectoF}.
\end{proof}
\end{subclaim}

Since $N_{H}(w)\subseteq C$ for a given a $w\in F$, Claim~\ref{cl:ectoF} implies that  \[e_{H}(C,A\cup R\cup F)\geq |A\cup R\cup F|+ |N_{H}(w)\cap C|-1\geq |A\cup R\cup F|+\delta_{+}-1.\]

To finish the proof of Claim~\ref{cl:modify}, we first consider the case $\Delta_{2}\geq \Delta_{1}$, and therefore, $m=\Delta_{1}+1$ and $\lambda= \delta_{+}+i$ for some $i\geq 1$. We can reduce (\ref{eq:ARFFull}) in the following way \begin{align}
    |A\cup R\cup F|&\geq \Delta_{2}(\Delta_{1}-\delta_{+})-\Delta_{2}+1+\Delta_{1}-deg_{H}(s)+\Delta_{2}+1-m+\lambda-1-t'\notag \\
    &=\Delta_{2}(\Delta_{1}-\delta_{+})-\Delta_{2}+1+\Delta_{1}-deg_{H}(s)+\Delta_{2}+1-(\Delta_{1}+1)+\delta_{+}+i-1-t'\notag\\
    &=\Delta_{2}(\Delta_{1}-\delta_{+})+\Delta_{1}-deg_{H}(s)-\Delta_{1}+\delta_{+}+i-t'\notag\\
    &=(\Delta_{2}-1)(\Delta_{1}-\delta_{+})+\Delta_{1}-deg_{H}(s)+i-t'.\notag
\end{align}

Since $|C|\leq \Delta_{2}-1$ and no vertex in $C$ is adjacent in $G$ to more than $\Delta_{1}-\delta_{+}$ vertices in $A\cup R\cup F$, we conclude that \[\delta_{+}-1+|A\cup R\cup F|\leq e_{H}(C, A\cup R\cup F)\leq |C|(\Delta_{1}-\delta_{+})\leq (\Delta_{2}-1)(\Delta_{1}-\delta_{+}).\] This can only be true if $\delta_{+}=1$, $i=t'=1$, $\Delta_{1}=deg_{H}(s)$, $|C|=\Delta_{2}-1$, $|N_{H}(u)\cap F|=\Delta_{1}-\delta_{+}$ for all $u\in C$, and $|A\cup R\cup F|=(\Delta_{2}-1)(\Delta_{1}-\delta_{+})$. This implies $A=B=R=\emptyset$.

Claim~\ref{cl:delta_+special} implies $l\neq y$, $\Delta_{2}\geq 3$, and therefore, $|C|\geq 2$. Since $|C|=\Delta_{2}-1$ and $N_{G_{2}}(y)=\{s\}\cup C$, $y$ and $l$ are not adjacent in $G$. Since every vertex in $C$ is adjacent in $H$ to $\Delta_{1}-1$ vertices in $F$, we may conclude that no vertex in $C$ is adjacent in $H$ to two vertices in $C\cup D\cup \{s\}$. Since $t'=1$, there is exactly one vertex $v\in C\cup D$ that is adjacent in $G_{2}$ to $s$. If $v\in C$, then since $\Gamma'(s)\cup \{s\}=N_{G_{2}}(s)$, Claim~\ref{cl:wNieghbors} says $v$ is adjacent in $H$ to some $v'$ in $C$. However, this implies that $v'$ would be adjacent in $H$ to $v$ and $s$. This contradiction implies that $v$ must be in $D$.  Let $u$ and $u'$ be in $C$. Since $u$ and $u'$ are adjacent in $H$ to $s$, Claim~\ref{cl:wNieghbors} says that $uu'$ is an edge of $G_{2}$. Moreover, since $u$ and $u'$ are adjacent in $G_{2}$ to the $\Delta_{2}-1$ vertices in $\Gamma_{1}'(s)$ we may conclude that neither $u$ nor $u'$ are adjacent in $G$ to $v$ nor $v'$. Let $w\in N_{H}(u)\cap F$ and $w'\in N_{H}(u')\cap F$. Since $\{ys,wu,vv',u'w'\}$ are edges of $H$ and $\{sw,uv, v'u'\}$ are non-edges of $G$, Claim~\ref{scl:wandl} implies $l$ is adjacent in $G_{2}$ to $w'$. However, this is a contradiction. Therefore, $\Delta_{2}<\Delta_{1}$. 

In this case, $m=\Delta_{2}+1$. If $\lambda\geq\Delta_{2}+1$, then (\ref{eq:ARFFull}) says $|A\cup R\cup F|\geq \Delta_{2}(\Delta_{1}-\delta_{+}).$ However, this leads to the contradiction $(\Delta_{2}-1)(\Delta_{1}-\delta_{+})\geq e_{H}(C,A\cup R\cup F)\geq \Delta_{2}(\Delta_{1}-\delta_{+})$.
If $\Delta_{2}+1>\lambda\geq\delta_{+}+1$, then (\ref{eq:ARFFull}) implies $|A\cup R\cup F|\geq \Delta_{2}(\Delta_{1}-\delta_{+})-\Delta_{2}+\delta_{+}+1-t'$. This gives us
$(\Delta_{2}-1)(\Delta_{1}-\delta_{+})\geq e_{H}(C,A\cup R\cup F)\geq \Delta_{2}(\Delta_{1}-\delta_{+})-\Delta_{2}+\delta_{+}+1-t'$. After, moving all the terms on to the right hand side and simplifying, we deduce the following contradiction 
\[0\geq \Delta_{1}-\delta_{+}-\Delta_{2}+\delta_{+}+1-t'\geq \Delta_{1}-\Delta_{2}+1-t'\geq 1.\]
This completes the proof of Claim~\ref{cl:modify}. \qedhere
\end{claimlateproof}

\bibliographystyle{hsiam}

\begin{thebibliography}{10}

\bibitem{Aigner1993}
{\sc M.~Aigner and S.~Brandt}, {\em Embedding arbitrary graphs of maximum
  degree two}, Journal of the London Mathematical Society, s2-48 (1993),
  pp.~39--51, https://doi.org/10.1112/jlms/s2-48.1.39.

\bibitem{Alon1996}
{\sc N.~Alon and E.~Fischer}, {\em 2-factors in dense graphs}, Discrete
  Mathematics, 152 (1996), pp.~13--23,
  https://doi.org/10.1016/0012-365X(95)00242-O.

\bibitem{Bollobas1978}
{\sc B.~Bollob{\'{a}}s and S.~E. Eldridge}, {\em Packings of graphs and
  applications to computational complexity}, Journal of Combinatorial Theory,
  Series B, 25 (1978), pp.~105--124,
  https://doi.org/10.1016/0095-8956(78)90030-8.

\bibitem{Bollobas2008}
{\sc B.~Bollob{\'{a}}s, A.~Kostochka, and K.~Nakprasit}, {\em Packing
  d-degenerate graphs}, Journal of Combinatorial Theory, Series B, 98 (2008),
  pp.~85--94, https://doi.org/10.1016/j.jctb.2007.05.002.

\bibitem{Brooks1941}
{\sc R.~L. Brooks}, {\em On colouring the nodes of a network}, Mathematical
  Proceedings of the Cambridge Philosophical Society, 37 (1941), pp.~194--197,
  https://doi.org/10.1017/S030500410002168X.

\bibitem{Busch2012}
{\sc A.~H. Busch, M.~J. Ferrara, S.~G. Hartke, M.~S. Jacobson, H.~Kaul, and
  D.~B. West}, {\em Packing of graphic n-tuples}, Journal of Graph Theory, 70
  (2012), pp.~29--39, https://doi.org/10.1002/jgt.20598.

\bibitem{CamesVanBatenburg2018}
{\sc W.~Cames Van~Batenburg and R.~J. Kang}, {\em Packing graphs of bounded
  codegree}, Combinatorics, Probability and Computing, 27 (2018), pp.~725--740,
  https://doti.org/10.1017/S0963548318000032.

\bibitem{Catlin1976}
{\sc P.~A. Catlin}, {\em Embedding subgraphs and coloring graphs under extremal
  degree conditions}, PhD thesis, The Ohio State University, 1976.

\bibitem{Csaba2007}
{\sc B.~Csaba}, {\em On the {Bollob{\'{a}}s--Eldridge} conjecture for bipartite
  graphs}, Combinatorics, Probability and Computing, 16 (2007), pp.~661--691,
  https://doi.org/10.1017/S0963548307008395.

\bibitem{Csaba2003}
{\sc B.~Csaba, A.~Shokoufandeh, and E.~Szemer{\'e}di}, {\em Proof of a
  conjecture of {Bollob{\'a}s} and {Eldridge} for graphs of maximum degree
  three}, Combinatorica, 23 (2003), pp.~35--72.

\bibitem{Csaba2019}
{\sc B.~Csaba and B.~V{\'{a}}s{\'{a}}rhelyi}, {\em On embedding degree
  sequences}, Informatica,  (2019), pp.~23--31,
  https://doi.org/10.31449/inf.v43i1.2684.

\bibitem{Diemunsch2015}
{\sc J.~Diemunsch, M.~Ferrara, S.~Jahanbekam, and J.~M. Shook}, {\em Extremal
  theorems for degree sequence packing and the two-color discrete tomography
  problem}, SIAM Journal on Discrete Mathematics, 29 (2015), pp.~2088--2099,
  http://dx.doi.org/10.1137/140987912.

\bibitem{Diestel2016}
{\sc R.~Diestel}, {\em Graph theory}, vol.~173 of Graduate Texts in
  Mathematics, Springer-Verlag, Heidelberg, fifth~ed., Aug. 2016.

\bibitem{Eaton2000}
{\sc N.~Eaton}, {\em A near packing of two graphs}, Journal of Combinatorial
  Theory, Series B, 80 (2000), pp.~98--103,
  https://doi.org/10.1006/jctb.2000.1971.

\bibitem{Goldberg2011}
{\sc M.~Goldberg and M.~Magdon-Ismail}, {\em Embedding a forest in a graph},
  The electronic journal of combinatorics, 18 (2011),
  https://doi.org/10.37236/586.
\newblock No. P99.

\bibitem{Gollakota2020}
{\sc A.~Gollakota, W.~Hardt, and I.~Mikl\'{o}s}, {\em Packing tree degree
  sequences}, Graphs and Combinatorics, 36 (2020), pp.~779--801,
  https://doi.org/10.1007/s00373-020-02153-0.

\bibitem{Hajnal1970}
{\sc A.~Hajnal and E.~Szemer{\'{e}}di}, {\em Proof of a conjecture of{
  Erd{\H{o}}s}}, in Combinatorial Theory and Its Applications, Vol.{~}2 (Ed.
  P.{~}Erd{\H{o}}s, A.{~}R{\'{e}}nyi, and V.{~}T.{~}S{\'{o}}s), 1970,
  pp.~601--623.

\bibitem{Hakimi1962}
{\sc S.~L. Hakimi}, {\em On realizability of a set of integers as degrees of
  the vertices of a linear graph. {I}}, Journal of the Society for Industrial
  and Applied Mathematics, 10 (1962), pp.~496--506,
  https://doi.org/10.1137/0110037.

\bibitem{Havel1955}
{\sc V.~Havel}, {\em Pozn{\'{a}}mka o existenci kone{\v{c}}n{\'{y}}ch
  graf{\r{u}}}, {{\v{C}}}asopis pro p{\v{e}}stov{\'{a}}n{\'{i}} matematiky, 080
  (1955), pp.~477--480.

\bibitem{Kano1992}
{\sc M.~Kano and N.~Tokushige}, {\em Binding numbers and f-factors of graphs},
  Journal of Combinatorial Theory, Series B, 54 (1992), pp.~213--221,
  https://doi.org/10.1016/0095-8956(92)90053-Z.

\bibitem{Katerinis2000}
{\sc P.~Katerinis and N.~Tsikopoulos}, {\em Minimum degree and f-factors in
  graphs}, New Zealand Journal of Mathematics, 29 (2000), pp.~33--40.

\bibitem{Kaul2007}
{\sc H.~Kaul and A.~Kostochka}, {\em Extremal graphs for a graph packing
  theorem of {Sauer} and {Spencer}}, Combinatorics, Probability and Computing,
  16 (2007), pp.~409--416, https://doi.org/10.1017/S0963548306007929.

\bibitem{Kaul2008}
{\sc H.~Kaul, A.~Kostochka, and G.~Yu}, {\em On a graph packing conjecture by
  {Bollob{\'a}s}, {Eldridge} and {Catlin}}, Combinatorica, 28 (2008),
  pp.~469--485, https://doi.org/10.1007/s00493-008-2278-0.

\bibitem{Kaul2022}
{\sc H.~Kaul and B.~Reiniger}, {\em A generalization of the graph packing
  theorems of sauer-spencer and brandt}, Combinatorica, 42 (2022),
  pp.~1347--1356, https://doi.org/10.1007/s00493-022-4932-3.

\bibitem{Kierstead2008}
{\sc H.~A. Kierstead and A.~V. Kostochka}, {\em A short proof of the{
  Hajnal–Szemerédi} theorem on equitable colouring}, Combinatorics,
  Probability and Computing, 17 (2008), pp.~265--270,
  https://doi.org/10.1017/S096354830700861.

\bibitem{Kierstead2009}
{\sc H.~A. Kierstead, A.~V. Kostochka, and G.~Yu}, {\em Extremal graph packing
  problems: ore-type versus dirac-type}, London Mathematical Society Lecture
  Note Series, Cambridge University Press, 2009, pp.~113--136,
  https://doi.org/10.1017/CBO9781107325975.006.

\bibitem{Kostochka2007}
{\sc A.~Kostochka and G.~Yu}, {\em An {Ore}-type analogue of the
  {Sauer}-{Spencer} theorem}, Graphs and Combinatorics, 23 (2007),
  pp.~419--424, https://doi.org/10.1007/s00373-007-0732-1.

\bibitem{Plummer2007}
{\sc M.~D. Plummer}, {\em Graph factors and factorization: {1985–2003}: a
  survey}, Discrete Mathematics, 307 (2007), pp.~791--821,
  https://doi.org/10.1016/j.disc.2005.11.059.
\newblock Cycles and Colourings 2003.

\bibitem{Sauer1978}
{\sc N.~Sauer and J.~Spencer}, {\em Edge disjoint placement of graphs}, Journal
  of Combinatorial Theory, Series B, 25 (1978), pp.~295--302,
  https://doi.org/10.1016/0095-8956(78)90005-9.

\bibitem{Shook2022}
{\sc J.~M. Shook}, {\em On a conjecture that strengthens the $k$-factor case of
  kundu's $k$-factor theorem}, 2022, https://doi.org/10.48550/arXiv.2205.01645.

\bibitem{Tyshkevich2000}
{\sc R.~Tyshkevich}, {\em Decomposition of graphical sequences and unigraphs},
  Discrete Mathematics, 220 (2000), pp.~201--238,
  https://doi.org/10.1016/S0012-365X(99)00381-7.

\bibitem{Wang1994}
{\sc H.~Wang}, {\em Packing a forest with a graph.}, Australas. J Comb., 10
  (1994), pp.~205--210.

\bibitem{Yap1988}
{\sc H.~P. Yap}, {\em Packing of graphs---a survey}, Discrete Mathematics, 72
  (1988), pp.~395--404, https://doi.org/10.1016/0012-365X(88)90232-4.

\bibitem{Yin2016}
{\sc J.-H. Yin}, {\em A note on packing of graphic n-tuples}, Discrete
  Mathematics, 339 (2016), pp.~132--137,
  https://doi.org/10.1016/j.disc.2015.07.017.

\end{thebibliography}

\end{document}